%% file: ER_Annealed_v3.tex
\documentclass[]{article}

\usepackage{geometry}
\input{ER_preamble_biblio.tex}

\input{ER_preamble_main.tex}

\numberwithin{equation}{section}

\allowdisplaybreaks

\title{\sffamily Cutoff for Random Walk on Dynamical Erd\H{o}s--R\'enyi Graph}
\author{\sffamily Sam Olesker-Taylor\quad Perla Sousi}
\date{\vspace{-5ex}}

\begin{document}

\maketitle

\renewcommand{\abstractname}{\sffamily Abstract}

\begin{abstract}
	\noindent
	We consider dynamical percolation on the complete graph $K_n$, where each edge refreshes its state at rate $\mu \ll 1/n$, and is then declared open with probability $p = \lambda/n$ where $\lambda > 1$. We study a random walk on this dynamical environment which jumps at rate $1/n$ along every open edge. We show that the mixing time of the full system exhibits cutoff at $\tfrac32 \log n/\mu$. We do this by showing that the random walk component mixes faster than the environment process; along the way, we control the time it takes for the walk to become isolated.
\end{abstract}

\small
\begin{quote}
	\begin{description}[noitemsep, topsep=0pt, font={\mdseries\itshape}]
		\item [Keywords:]
		dynamical percolation, Erd\H{o}s--R\'enyi, random walk, mixing times, coupling
		
		\item [MSC 2020 subject classifications:]
		05C81, 60J27, 60K35, 60K37
	\end{description}
\end{quote}
\normalsize


\sffamily
\setcounter{tocdepth}{1}
\normalfont

\blfootnote{%
	Statistical Laboratory, University of Cambridge, UK%
\qquad%
	\url{sam.ot@posteo.co.uk} and \url{p.sousi@statslab.cam.ac.uk}%
\\%
	This research was supported by the EPSRC:%
\quad%
	SOT by Doctoral Training Grant \#1885554
and
	PS by EP/R022615/1%
\\%
	This article was originally published under the first named author's previous name, ``Sam Thomas''%
}


\romannumbering

\section{Introduction}
\label{sec:intro}

In this paper we consider a random walk on a dynamically evolving random graph. Fix an underlying {(undirected)} graph $G = (V,E)$.
Write $n = \abs V$.
The dynamics of the graph are that of \textit{dynamical percolation}: edges refresh {independently} at rate $\mu$; upon refreshing, the edge is declared open with probability $p$ and closed with probability $1-p$. We denote the state of the graph at time $t$ by $\eta_t \in \{0,1\}^E$: 0 corresponds to a closed edge and 1 to an open edge. The location of the random walker at time $t$ is denoted $X_t \in V$: it moves at rate 1; when its exponential clock rings, it chooses uniformly at random a neighbour $v$ of its current location, $x$ say, and jumps from $x$ to $v$ if and only if the edge connecting $x$ and $v$ is open (at this time), otherwise it remains in place.

Write $\pi_\RW$ for the invariant distribution of the nearest-neighbour simple random walk on $G$ (ie the degree-biased distribution) and $\pi_p$ for the product measure on $\{0,1\}^E$ with density $p$.
{The full process, $(X,\eta)$, is reversible with invariant distribution $\pi = \pi_\RW \times \pi_p$.}

We emphasise that the pair $(X_t,\eta_t)_{t\ge0}$ is Markovian, as is just the graph process $(\eta_t)_{t\ge0}$, while the location of the walker alone $\Xt$ is not: indeed, its transitions depend on the current graph. (Note that $(\eta_t)_{t\ge0}$ is a biased simple random walk on the hypercube $\{0,1\}^E$.)
For all our results, we take $p = \lambda/n$, and emphasise that~$\lambda$ is a fixed constant, while $n$ and $\mu = \mu_n$ vary.

This model was introduced by Peres, Stauffer and Steif in \cite{PSS:dyn-perc:sub}. They used the torus $\bz_n^d$ (with~$d$ fixed) as their underlying graph; in this paper we use the complete graph $K_n$ as our underlying graph. Hence from now on we take
\[
	V = \{1,...,n\}
\Quad{and}
	E = \brb{ (i,j) \mid i,j \in \{1,...,n\}, \, i \neq j }.
\]
Percolation on the complete graph gives precisely the Erd\H{o}s--R\'enyi graph, and hence the name `dynamical Erd\H{o}s--R\'enyi'; we denote the measure of an Erd\H{o}s--R\'enyi graph by $\pi_\ER$, and note that in this case $\pi_p = \pi_\ER$. Also, we denote the uniform measure on $\{1,...,n\}$ by $\pi_U$, and note that in this case $\pi_\RW = \pi_U$.

{Taking $p = \lambda/n$, for a constant $\lambda$, corresponds to the sparse regime for percolation, in which the expected degree of a vertex is order 1 (ie independent of $n$ in the limit)}.
{Since the majority of the degrees are order $1$, the walker takes steps on the timescale~$n$ (the majority of the time).}
We consider $\mu \ll 1/n$ so that the walk takes a large number of steps before seeing a local update to the graph; for bounded degree graphs, we could take $\mu \ll 1$. Our proofs actually require very slightly more, namely a polylogarithmic factor: we consider $\mu \ll \logn[-\alpha]/n$ for a fixed $\alpha > 0$; no attempt has been made to optimise this parameter.

As in \cite{PSS:dyn-perc:sub}, we look at the $\eps$-mixing time of the \textit{full system} $(X,\eta)$:
\[
	\tmix(\eps) = \inf\brb{ t \ge 0 \midb \max_{x_0,\eta_0} \, \tvb{\prb[x_0,\eta_0]{(X_t,\eta_t) \in \cdot} - \pi} \le \epsilon }.
\]
{When the leading order term of $\tmix(\eps)$ is independent of $\eps$, we say that there is \textit{cutoff}. The \textit{cutoff window} is given by the order of $\tmix(\eps) - \tmix(1-\eps)$, which will depend on $\eps$.}

In order to state our mixing result, we must first define some notation for iterated logarithm:
\[
	\text{set}
\quad
	\log_{(1)}\! n = \logn
\Quad{and define inductively}
	\log_{(m+1)}\! n = \log(\log_{(m)}\! n)
\Quad{for}
	m \ge 1.
\]
Our main result considers the \textit{supercritical regime} of percolation, ie has $p = \lambda/n$ where $\lambda > 1$ is a constant, and states that the full system $(X,\eta)$ then exhibits cutoff at time $\tfrac32(\logn)/\mu$ with cutoff window of smaller order than $(\log_{(M)}\! n)/\mu$ for all $M$.

\begin{thm}[Cutoff for Full System]
\label{res:main_thm:mixing_full}
	For all~$\lambda > 1$, all~$\eps \in (0,1)$, all~$M \in \bn$ and all~$n$ sufficiently large, for $p = \lambda/n$ and $\mu \le \logn[-20]/n$, we have
	\[
		\absb{ \tmix(\eps) - \tfrac32 (\logn)/\mu }
	\le
		(\log_{(M)}\! n)/\mu.
	\]
\end{thm}

We also consider the `mixing' of the \textit{random walk component}:
\[
	\tmix^\RW(\eps,\eta_0) = \sup\brb{ t \ge 0 \midb \max_{x_0} \, \tvb{ \prb[x_0,\eta_0]{ X_t \in \cdot } - \pi_U } \ge \eps }
\Quad{for}
	\eta_0 \in \{0,1\}^E.
\]
Since $X$ is not a Markov chain, we do not have a priori that the total variation distance from uniform is decreasing, hence we do not define the mixing time to be `the first time the total variation distance is below $\eps$', but rather `the last time the total variation distance is above $\eps$'. (Of course, for a Markov chain, these notions are the same.) Note that, trivially by projection, $\tmix^\RW(\eps,\eta_0) \le \tmix(\eps)$ for all~$\eps$ and all~$\eta_0$. We show that `the walk mixes faster than the environment' when the initial environment is `typical', in the following precise sense.

\begin{thm}[Mixing of Random Walk]
\label{res:main_thm:mixing_rw}
	For all~$\lambda > 1$, all~$\eps \in (0,1)$, all~$M \in \bn$ and all~$n$ sufficiently large, for $p = \lambda/n$ there exists a subset $H \subseteq \{0,1\}^E$ with $\pi_\ER(H) = 1 - \oh1$ so that, for all~$\eta_0 \in H$, for $\mu \le \logn[-20]/n$, we have
	\[
		\tmix^\RW(\eps,\eta_0) \le \log_{(M)}\! n / \mu;
	\]
	ie, if $\eta_0 \sim \pi_\ER$ then, for all~$\eps \in (0,1)$, all~$M \in \bn$ and all~$n$ sufficiently large, we have
	\[
		t_\mix^\RW(\eps,\eta_0) \le \log_{(M)}\! n / \mu
	\Quad{with probability}
		1 - \oh1.
	\]
\end{thm}

\begin{rmkt*}
\label{rmk:cutoff-intuition}
	From this intuitively it is clear why we get cutoff.
		First, let the environment mix.
	This is just a (biased) random walk on the hypercube $\{0,1\}^N$ where $N = \binom{n}{2}$, and {has cutoff at time $\tfrac12 \log(N/p) / \mu = (\tfrac32 \logn + \Oh1)/\mu$. We prove this in \Cref{res:hyp}; cf \cite[Example 12.19]{LPW:markov-mixing}, where the unbiased case is considered.}
	At this time, the graph is `approximately' Erd\H{o}s--R\'enyi, and so likely to be in the set~$H$ (from \Cref{res:main_thm:mixing_rw}).
	Finally we let the walk mix.
	This takes time little-$o$ of the mixing of the environment. {This is indeed the heuristic that we use, but one has to be careful due to correlations between the walk and environment.}
\end{rmkt*}

	%
Above we are allowed to choose $X_0$ dependent on $\eta_0$ (and vice versa). In \S\ref{sec:iiie}, we consider drawing $\eta_0$ according to $\pi_\ER$, and then choosing $X_0$ \emph{independently} of $\eta_0$. By symmetry, we may assume $X_0 = 1$. We then look at the mixing time of the walk on this (evolving) graph:
\[
	\tmix^\RW(\eps) = \sup\brb{ t \ge 0 \midb \tvb{ \prb[1,\ER]{ X_t \in \cdot } - \pi_U } \ge \eps },
\]
where $\pr[x_0,\ER]{ \cdot } = \sum_{\eta_0} \pr[x_0,\eta_0]{ \cdot } \pi_\ER(\eta_0)$ averages the initial environment with respect to the Erd\H{o}s--R\'enyi measure. We prove a sharp, up to constants, result on this mixing time. The result does not require us to consider the supercritical regime, ie $\lambda > 1$, but allows any~$\lambda \in (0,\infty)$, including the critical case $\lambda = 1$. In \S\ref{sec:iiie}, we prove the following.

\begin{thm}
\label{res:iiie:statement}
	For all~$\lambda \in (0,\infty)$, there exists a constant $C$ so that, for all~$\eps \in (0,1)$ and all~$n$ sufficiently large, $p = \lambda/n$, we have the following bounds on the mixing time:
	\begin{alignat*}{2}
		\tmix^\RW(\eps)
	&\ge
		\tfrac1\mu \cdot \tfrac1{2\lambda} \logeps
	\quad&&\text{if}\quad
		\eps \in (0,e^{-3\lambda} \wedge 1)
	\Quad{for any}
		\mu;
	\\
		\tmix^\RW(\eps)
	&\le
		\tfrac1\mu \cdot C \logeps
	\quad&&\text{if}\quad
		\eps \in (0,\tfrac14)
	\Quad{when}
		\mu \le \tfrac23(1+\lambda)^{-1}/n.
	\end{alignat*}
\end{thm}


Note that once all edges of the graph have been refreshed, the graph has the Erd\H{o}s--R\'enyi measure, and is independent of $\eta_0$. This is the coupon-collector problem, and takes time concentrated at $\log N \approx 2 \logn$. At first glance, then, it appears that the idea of the above remark along with \Cref{res:iiie:statement} can be applied to easily give pre-cutoff at $2 \logn$.
However, this is not the case: to prove the above statement, it is crucial that $X_0$ is chosen \emph{independently} of~$\eta_0$; we then exploit symmetry. At time $2\logn$, the walk and environment are correlated, and so the argument does not apply; hence the need for \Cref{res:main_thm:mixing_rw}.
	%

\medskip

The above statement \emph{suggests} that the upper bound in \Cref{res:main_thm:mixing_rw} is probably not sharp. However, the main interest in \Cref{res:main_thm:mixing_rw} compared with \Cref{res:main_thm:mixing_full} is not the specific upper bound, but the fact that `most' $\eta_0$ have $\tmix^\RW(\eps,\eta_0) \ll \tmix(\eps)$, ie that the walk mixes faster than the full system. This is what allows us to show cutoff.

\begin{rmkt*}
	For the rest of the paper, with the exception of \S\ref{sec:iiie}, we assume that $\mu \le \logn[-20]/n$; we shall not repeat this in the statement of every theorem.
\end{rmkt*}

\medskip

The model of dynamical percolation (without the random walker) was introduced by H\"aggstr\"om, Peres and Steif in \cite{HPS:dyn-perc}.
The model with the random walk was then introduced by Peres, Stauffer and Steif in \cite{PSS:dyn-perc:sub}, with underlying graph the torus $\bz_n^d$. They considered the \textit{subcritical regime}, ie $p < p_c(\bz^d)$, the critical probability for bond percolation on $\bz^d$, and obtained the correct order for the mixing of the full system, showing further that the order of the mixing of the walk is the same order as the mixing of the full system (in contrast to our model).
{While they give a very complete picture of the subcritical mixing, key to their proofs is that the percolation clusters are all small, and that there is no giant component (taking up a constant proportion of the vertices).}

The \textit{supercritical regime}, ie $p > p_c(\bz^d)$ where such a giant does exist, was then considered by Peres, Sousi and Steif in \cite{PSS:dyn-perc:sup}.
They considered the `quenched' case, where a `typical' environment process $\{\eta_t\}_{t\ge0}$ is fixed in advance, and the walker walks on this. They obtained the correct order for the mixing of the walk, up to polylogarithmic factors, but only in the regime where $\theta(p) > \tfrac12$, ie the probability that the component at 0 in $\bz^d$ is infinite is greater than $\tfrac12$.
The case $\theta(p) \le \tfrac12$ remains open.

Avena, G\"ulda\c{s}, van der Hofstad and den Hollander in \cite{AGHH:dynamic-cm,AGHH:dynamic-cm2} studied the mixing time of the non-backtracking random walk on a \textit{dynamical configuration model}. The \textit{configuration model} generates a random graph with a prescribed degree sequence, and the dynamics at every time step `rewire' uniformly at random a given proportion of the edges: {this rewiring involves cutting edges into two half-edges and then randomly repairing the half-edges}.

It is straightforward to see that in our model when the walker first crosses a refreshed edge (which is a randomised stopping time) it is then (almost) uniform. The authors of \cite{AGHH:dynamic-cm,AGHH:dynamic-cm2} considered an analogous time for the nearest-neighbour simple random walk, namely the first time the walk crosses a rewired edge, and showed that the distribution of the walk at this time is (almost) its invariant distribution. In both cases, however, these times are not sufficient to show mixing: it is not the case that the walk `remains close to invariant'; for example, there is significant probability that the walk will cross back over the same edge to its previous location.
{It is possible that a more refined analysis of a related stopping time---eg the first time the walk crosses a rewired edge and then `escapes', not recrossing this edge again (for a long time)---would work. This approach is not taken in \cite{AGHH:dynamic-cm,AGHH:dynamic-cm2}, though, and is left open.
Rather, to resolve this the authors consider the non-backtracking random walk
which, along with the locally tree-like structure of the configuration model, removes this `crossing back' issue.
They then show sharp asymptotics for the mixing time of this non-backtracking random walk, using the aforementioned stopping time.}

\smallskip

{In contrast to the above examples, in our work we show cutoff for the full process $(X,\eta)$ for the entire supercritical regime, ie consider $p = \lambda/n$ for any constant $\lambda > 1$, and obtain the correct order of the mixing of the walk, up to an iterated log factor.
Furthermore, our methods adapt immediately to the subcritical regime, ie $p = \lambda/n$ with $\lambda < 1$. We have not considered the details for this, but with a few concentration results on the structure of a subcritical Erd\H{o}s--R\'enyi graph, a similar mixing result will follow.}

\smallskip

Fountoulakis and Reed in \cite{FR:giant-mixing} and Benjamini, Kozma and Wormald in \cite{BKW:giant-mixing} studied the mixing time of the nearest-neighbour simple random walk on the giant component of a supercritical Erd\H{o}s--R\'enyi random graph, without any graph dynamics: they prove that the mixing time is order exactly~$\logn[2]$. Fountoulakis and Reed carefully studied the ratio between the size of the edge boundary of a set and the set itself, using a variation of the Lov\'asz-Kannan integral, which they developed in \cite{FR:bottleneck-mixing}.
Benjamini, Kozma and Wormald used a more geometric approach, defining a stripping process to analyse the (2-)\textit{core} and the \textit{kernel} of the graph; they show the kernel is a (type of) expander, and describe the \textit{decorations} attached to the kernel.

The two works above consider mixing from the worst-case starting point. Berestycki, Lubetzky, Peres and Sly in \cite{BLPS:giant-mixing} consider mixing when the starting point is chosen according to the invariant distribution. They show then that the mixing time is actually order $\logn$, obtaining the correct constant and also showing cutoff; contrast this with order $\logn[2]$ for the worst-case.

\section{Outline of Proof and Preliminaries}

\subsection{Outline of Proofs}

We now give a brief, informal outline of the proofs of the main results. First consider the following scenario: suppose a walker is isolated at vertex $u$ {(ie the walk is at the vertex $u$ which is an isolated vertex in the current graph)}, and suppose it becomes non-isolated by the edge $(u,v)$ opening, where $v$ was isolated immediately before $(u,v)$ opened. Now the pair $\{u,v\}$ is a component of the graph. Because $\mu \ll 1/n$, we see that the walker takes a large number of steps before this edge closes. If it closes before any other edge incident to $\{u,v\}$ opens (which has order 1 probability, by counting open/closed edges), then the walker is approximately uniformly distributed on $\{u,v\}$. So it has approximately `done a lazy simple random walker step'.

This motivates the following coupling. First wait for the two walkers to be simultaneously isolated \emph{in the same environment} of two full systems. Then to couple we want to imitate the standard coupling of the lazy simple random walk on the complete graph; we do this by considering the event that when the walkers become non-isolated they connect to a vertex that was isolated immediately prior. We give a precise definition of the coupling that we use in \S\ref{sec:coup}.

In order to find the time it takes for two walkers to be simultaneously isolated in the same environment, we first consider how long it takes one walker to become isolated. To find this time, we observe that a walker can only become isolated if it is at a degree~1 vertex and this vertex becomes isolated prior to the walk's leaving it. This motivates looking at the rate at which degree~1 vertices are hit. To do this, we compare the number of degree~1 vertices hit by a walker on the dynamic graph and the same quantity for a walker on a static graph; we then apply a Chernoff-style bound due to Gillman \cite{G:gillman-chernoff}. We give the precise details of this in \S\ref{sec:isol-times}.

A key element in studying the isolation time is to control how long the walk remains in the giant once it has entered. We show that since the graph updates slowly, as $\mu n \ll \logn[-19]$, the walk does not see updates to the graph for some while; this time is long enough for the walk to become approximately uniform on the giant prior to seeing a change. We can then use structure results on the Erd\H{o}s--R\'enyi graph to see how `near the core' of the giant the walk is.

We then use the fact that an Erd\H{o}s--R\'enyi graph with one vertex conditioned to be isolated has the distribution of an Erd\H{o}s--R\'enyi graph on $n-1$ vertices (with edge-probability $p$) union an isolated vertex: this allows us to say that `conditioning on one walker's being isolated has almost no affect on the other walker', which will allow us to treat the walkers as almost independent. We give the precise details of this in \S\ref{sec:isol2}.

\subsection{Notation and Terminology}

For functions $f$ and $g$, we write
$f(n) \lesssim g(n)$, or $f(n) = \Oh{g(n)}$,
if there exists a positive constant $C$ so that $f(n) \le C g(n)$ for all~$n$; write
$f(n) \gtrsim g(n)$, or $f(n) = \Om{g(n)}$, if $g(n) \lesssim f(n)$.
We write
$f(n) \asymp g(n)$, or $f(n) = \Th{g(n)}$,
if we have both $f(n) \lesssim g(n)$ and $g(n) \lesssim f(n)$. Write
$f(n) \ll g(n)$, or $f(n) = \oh{g(n)}$,
if $f(n)/g(n) \to 0$ as $\ninf$; write
$f(n) \gg g(n)$, or $f(n) = \om{g(n)}$,
if $g(n) \ll f(n)$.

For random variables $X$ and $Y$, we write $X \preccurlyeq Y$ if $Y$ stochastically dominates $X$ (from above), ie if $\pr{X \ge z} \le \pr{Y \ge z}$ for all~$z$; write $X \succcurlyeq Y$ if $Y \preccurlyeq X$.

For a real number $r > 0$, we write $\me(r)$ for the exponential random variable with rate $r$.

For real numbers $\alpha$ and $\beta$, we write $\alpha \wedge \beta = \min\{\alpha,\beta\}$.


\subsection{Good Graphs and Erd\H{o}s--R\'enyi Structure Results}
\label{sec:good-graphs}

In this section we state some results on Erd\H{o}s--R\'enyi graphs. Since a realisation can be \emph{any} graph, we want to describe explicitly what we shall mean by a \textit{good} graph.
We now define some notation for graph properties.
For the moment, let $c_*$ and $C_*$ be any two constants.

\begin{notationt} \label{notation}
Let $n \in \bn$; in the following, we suppress the $n$-dependence.
For a graph $G = (V,E)$ with $V = \{1, ..., n\}$, we use the following notation.
\begin{enumerate}
	\item 
	For $x \in V$, write $d(x)$ for the degree of $x$ (in $G$).
	
	
	\item 
	Write $\mg$ for the (set of vertices in the) largest component, and call it the \textit{giant}; if there is a tie, choose the option that includes the smallest labelled vertex.
	
	\item
	For $x \in V$, call an edge a \textit{removal edge} for $x$ if its removal breaks the component of $x$ in two, leaving $x$ in the smaller component (breaking ties as above). Write $\mr(x)$ for the set of removal edges for $x$, and write $R(x) = |\mr(x)|$.
	
	\item 
	For $M \in \bn$, write $\mw^M$ for the vertices of the giant with at most $C_* \log_{(M)}\! n$ removal edges, ie
	\[
		\mw^M = \brb{ x \in \mg \mid R(x) \le C_* \log_{(M)}\! n }.
	\]
	
	\item 
	Write $\gamma$ for the spectral gap and $\Phi$ for the isoperimetric constant, ie
	\[
		\Phi
	=
		\min\brbb{ \frac{|\partial S|}{d(S)} \midbb {S \subseteq V, \,} S \neq \emptyset, \, d(S) \le |E| }
	\]
	where $d(S) = \sum_{x \in S} d(x)$, $\partial S = \bra{ (x,y) \in E \mid x \in S, y \notin S }$; also write $\gamma$ for the spectral gap of the transition matrix of the nearest-neighbour simple random walk on $G$ (ie gap between 1 and the second largest eigenvalue).
\end{enumerate}
If the graph $G = \zeta$, then we use subscript $\zeta$, eg writing $\mg_\zeta$ for the giant of $\zeta$. If the graph $G = \eta_t$, then we add a subscript $t$, eg writing $\mg_t$, $R_t(x)$ or $d_t(x)$.
\end{notationt}

Let $\omega_*$ be any function of $n$.
{We now define what we mean by a \textit{good graph}.}

\begin{defnt}[Good Graph] \label{defn:good-graph}
	Let $n \in \bn$; in the following, we suppress the $n$-dependence.
	We say that a graph $G = (V,E)$ with $V = \{1,...,n\}$ is \textit{good}, and write $G \in \sg$, if it has a unique component $\mg$ with $|\mg| \ge C_* \logn$, which we call the \textit{giant}, which satisfies the following properties.
	\begin{enumerate}
		\item \label{gg:size}
		\emph{Size.}
		We have $|\mg| \ge c_* n$.
		
		
		\item \label{gg:max-deg}
		\emph{Maximum degree.}
		The maximum degree of $\mg$ is at most $C_* \logn$.
		
		
		\item \label{gg:num-edges}
		\emph{Number of edges.}
		There are at most $C_* n$ edges in $\mg$.
		
		\item \label{gg:num-deg1-giant}
		\emph{Number of degree~1 vertices.}
		The number of degree~1 vertices in $\mg$ is at least $c_* n$.
		
		\item \label{gg:removal-edges}
		\emph{Removal edges.}
		We have $R(x) \le C_* \logn$ for all $x \in \mg$.
		
		\item \label{gg:vtx-far-from-core}
		\emph{Vertices far from the core.}
		For all~$2 \le M \le \omega_*(n)$, the proportion of vertices $x$ of $\mg$ with $R(x) \ge C_* \log_{(M)}\! n$ is at most $(\log_{(M-1)}\! n)^{-4}$, ie $|\mg \setminus \mw^M|/|\mg| \le (\log_{(M-1)}\! n)^{-4}$.
		
		\item \label{gg:giant-expansion}
		\emph{Expansion properties.}
		We have $\Phi_\mg \ge c_* \logn[-2]$ and $\gamma_\mg \ge c_* \logn[-4]$.
	\qedhere
	\end{enumerate}
\end{defnt}

{This concept of \textit{good} depends on the choice of constants $c_*$ and $C_*$ and function $\omega_*$. The next proposition says that we can choose these parameters suitably so that an Erd\H{o}s--R\'enyi graph is overwhelmingly likely to be good.}
Write $\sg_t = \bra{ \eta_t \in \sg }$, and also
\[
	\Ggood{s,t}
=
	\brb{ \eta_u \in \sg \, \forall  u \in [s,t] }
=
	\cap_{s \le u \le t} \sg_u.
\qedhere
\]

\begin{prop}
\label{res:ER_good}
	There exist positive constants $c_*$ and $C_*$ and a function $\omega_*(n) \to \infty$ so that, for $G \sim \pi_\ER$, we have
	\[
		\prb{G \notin \sg}
	=
		\Ohb{ n^{-9} }
	\Quad{and}
		\prb{ \tfrac1n |\{x \in V \mid d(x) = 0\}| \le c_* }
	=
		\Ohb{ n^{-9} }.
	\]
\end{prop}

The claims in this proposition are fairly standard, but usually in the literature the proved decay rate is only $\oh1$, whereas we desire {the quantitative} $\Oh{n^{-9}}$.
As such we give the proof of this proposition, but defer it until the appendix.

\medskip

For the remainder of the paper, we select $(\omega_*,c_*,C_*)$ as guaranteed by this proposition and fix them permanently; whenever $\omega_*$, $c_*$ or $C_*$ is written below, it will refer to these constants.

We now consider our graph dynamics. We want not only the starting graph to be good, but we want it to remain good for a long time.

\begin{defnt}
\label{defn:very-good-graph}
	We make the following definitions:
	\begin{align*}
		\sh &= \brb{ \eta_0 \in \{0,1\}^E \midb \prb[\eta_0]{ \sg[0,1/\mu]^c } \le n^{-1} };
	\\
		\sh[s,t] &= \brb{ \eta_u \in \sh \, \forall  u \in [s,t] };
	\\
		H &= \brb{ \eta_0 \in \{0,1\}^E \midb \prb[\eta_0]{ \sh[0,n/\mu]^c } \le n^{-1} }.
	\end{align*}
	
	Further, if we are considering two environment processes, $\eta$ and $\xi$ say, then we (abuse notation slightly and) use the same notation, eg
	\(
		\sh[s,t] = \bra{ \eta_u,\xi_u \in \sh \, \forall  u \in [s,t] }.
	\)
\end{defnt}


We have the following result on `how good' an Erd\H{o}s--R\'enyi graph is.

\begin{prop}
\label{res:ER_very-good}
	For all~$t \le n/\mu$, if $\eta_0 \sim \pi_\ER$, then we have
	\[
		\prb{ \Hbad{0,t} }
	=
		\Ohb{ n^{-3} },
	\Quad{and hence}
		\pi_\ER(H) = 1 - \oh1.
	\]
	Moreover, these still hold if we add the condition that at least a proportion $c_*$ of the vertices are isolated to the definition of a good graph.
\end{prop}

We state a large deviations result on the Poisson distribution, which we use on a number of occasions throughout the paper.

\begin{lem}
	\label{res:poisson}
	We have the following bounds, valid for all~$\lambda > 0$ and all~$\eps \in (0,1)$:
	\[
		\prb{ \Po(\lambda) \ge (1+\eps)\lambda }
	\le
		\expb{ - \tfrac12 \lambda \eps^2 (1-\tfrac13\eps) }
	\Quad{and}
		\prb{ \Po(\lambda) \le (1-\eps)\lambda }
	\le
		\expb{ - \tfrac12 \lambda \eps^2 }.
	\]
\end{lem}

\begin{Proof}[Proof of \Cref{res:ER_very-good}]
	Since $\pi_\ER$ is the invariant distribution for the environment, using the concentration of the Poisson distribution
	and \Cref{res:ER_good} we find that
	\[
		\prb{ \sg[0,1/\mu]^c }
	=
		\prb{ \exists \, t \le 1/\mu \ST \eta_t \notin \sg }
	\le
		n^2 \cdot \Ohb{n^{-9}} + \expb{ - \tfrac1{12} n^2 }
	=
		\Ohb{ n^{-7} }.
	\]
	We now restrict the implicit sum in $\bp$ from $\eta_0 \in \{0,1\}^E$ to $\eta_0 \in \sh^c$:
	\begin{align*}
		\prb{ \sg[0,1/\mu]^c }
	&= \textstyle
		\sum_{\eta_0 \in \{0,1\}^E}
		\prb[\eta_0]{ \sg[0,1/\mu]^c }
		\pi_\ER(\eta_0)
	\\&\ge \textstyle
		\sum_{\eta_0 \in \sh^c}
		\prb[\eta_0]{ \sg[0,1/\mu]^c }
		\pi_\ER(\eta_0)
	\ge
		n^{-1} \cdot \pi_\ER(\sh^c).
	\end{align*}
	Hence we deduce that $\pi_\ER(\sh^c) = \Oh{n^{-6}}$. Repeating the same argument, we find that
	\[
		\pi_\ER(H^c) = \Ohb{n^{-2}}.
	\qedhere
	\]
\end{Proof}

\section{Hitting and Exit Times of the Giant}
\label{sec:hit/exit-giant}

In this section we study the hitting time of the giant, and how long the walk remains in the giant given that it starts there.
Write the following for the hitting and exit times of the giant:
\[
	\tau_\mg = \inf\brb{ t \ge 0 \mid X_t \in \mg_t }
\Quad{and}
	\tau_\mg' = \inf\brb{ t \ge 0 \mid X_t \notin \mg_t }.
\]


%
\begin{lem}[Hitting the Giant]
\label{res:hitgiant}
	There exists a positive constant $c$ so that, for all $n$ sufficiently large and all~$(x_0,\eta_0)$, we have
	\[
		\prb[x_0,\eta_0]{ \tau_\mg \le 1/\mu }
	\ge
		c - \prb[\eta_0]{ \Gbad{0,1/\mu} }.
	\]
\end{lem}

\begin{Proof}
	Write $\mU$ for the first time an edge incident to the walker refreshes and opens. By the memoryless property, $\mU \sim \me(\lambda \mu (1-1/n))$.
	When such an edge opens, it connects to the giant with probability $|\mg_\mU|/(n-1)$. Write $\theta_t = |\mg_t|/n$. We then have
	\begin{align*}
		\prb[x_0,\eta_0]{ \tau_\mg \le 1/\mu }
	&\ge
		\prb[x_0,\eta_0]{ \tau_\mg \le \mU, \, \mU \le 1/\mu, \, \theta_{\mU^-} \ge c_* }
	\\&\ge
		\prb[x_0,\eta_0]{ \tau_\mg \le \mU \mid \mU \le 1/\mu, \, \theta_{\mU^-} \ge c_* }
	\cdot
		\prb[x_0,\eta_0]{ \mU \le 1/\mu, \, \theta_t \ge c_* \, \forall t \le 1/\mu }
	\\&\ge
		c_* \rbb{ \prb[x_0,\eta_0]{ \mU \le 1/\mu } - \prb[\eta_0]{ \brb{ \theta_t \ge c_* \, \forall t \le 1/\mu }^c } }.
	\end{align*}
	The proof is completed by noting that $c_* \le 1$ and $U \sim \me(\lambda\mu(1-1/n))$ so $\pr[x_0,\eta_0]{\mU \le 1/\mu} \asymp 1$.
\end{Proof}

We now consider how long the walker remains in the giant once it enters. Recall the definition of $R$ and $\mw^M$ from \S\ref{sec:good-graphs}. Since the number of removal edges satisfies $R(x) \le C_* \logn$ for all~$x \in \mg$ when the graph is good, while the graph is good a trivial bound is $\tau_\mg' \succcurlyeq \me(C_* \mu \logn)$. We do a more careful analysis which shows that, for all~$M \in \bn$, `most of the time' $X_t \in \mw_t^M$, ie satisfies $R(X_t) \le C_* \log_{(M)}\! n$; this is because $|\mw^M|/|\mg| = 1 - \oh1$ for a good graph. The precise statement that we prove is as follows.

\begin{prop}[Exit Time from the Giant] \label{res:leavegiant}
	There exists a constant $C$ so that,
	for all~$M \in \bn$, all~$n$ sufficiently large, all~$t$ with $\logn[-5] \le \mu t \le \tfrac1{10}$ and all~$(x_0,\eta_0)$ with $\eta_0 \in \sh$ and $x_0 \in \mg_0$, we have
	\[
		\prb[x_0,\eta_0]{ \tau_\mg' \le t }
	\le
		C \mu t \log_{(M)}\! n.
	\]
\end{prop}

Since $\mu n \ll 1$, the `majority of the time' the walker takes a step before any edge incident to its location changes state. This motivates looking at a random walker moving on a \emph{static} graph, ie one without graph-dynamics. We call such a walk the \textit{static walk}, and the original walk (on the dynamic graph) the \textit{dynamic walk}; we denote them by $\tilde X$ and $X$, respectively.

Consider starting the two walks together. Observe that until the static walk encounters an edge that is in a different state to its original, the two walks have the same distribution, and hence we can couple them to be the same (until this time) as follows. Give $X$ and $\tilde X$ the same jump clock. When this clock rings, at time $t$ say, both walks choose the same vertex; $\tilde X$ performs the jump if and only if the connecting edge is present in $\eta_0$, while $X$ performs the jump if and only if the connecting edge is present in $\eta_t$. We call this the \textit{static-dynamic coupling}.

We define the \textit{set of edges seen by the walker} (in an interval $[s,t]$) as the set of all edges (open or closed) that are incident to the walker at some time (in an interval $[s,t]$).
Until $X$ sees an edge which is in a different state to its original, we can couple it with $\tilde X$, as described above.

\medskip

Label the edges of the (complete) graph $e_1, ..., e_N$, where $N = \binom{n}{2}$, in an arbitrary ordering (eg lexicographically). At time $t$, write $\so_t = \{o_t^1, o_t^2, ...\}$ for the (ordered) set of open edges in $\eta_t$ and $\sc_t = \{c_t^1, c_t^2, ...\}$ for the set of closed ones.
We say that an edge is a \textit{bridge} for a component if its removal splits the component into two (disconnected) parts. Recall also from \Cref{notation} the definition of a removal edge, and of $\mr$.

Fix $\rho = \logn[11]$. We now define an edge-set process $\mP = (\mP_t)_{t\ge0}$, {with $\mP_t \subseteq E$ for all $t$.}

\begin{defnt}[Set Processes]
{We define the edge-set process $\mP$ inductively. Throughout the definition, we assume that the graph is good, ie we define $(\mP_s)_{s \le t}$ on the event $\sg[0,t]$; recall the definition of good and $\sg[0,t]$ from \Cref{defn:good-graph} and the display after it.}

Suppose we have defined the set process $\mP$ up until time $t$, ie have defined $(\mP_s)_{s \in [0,t)}$. Let $V'$ be the most recent `update time' for the process $(\mP_s)_{s \in [0,t)}$ in the following sense:
\[
	V'
=
	\sup\brb{ s \in [0,t) \mid \text{an edge of $\mP_s$ changes state at time $s$} }.
\]
(If $t = 0$, then we take $V' = 0$. What follows is also used for the base case of the induction.)
For $r \ge 0$, let $\ma_{r+t}'$ (respectively $\mb_{r+t}'$) be the set of open (respectively closed) edges seen by $X$ in $[V',t+r]$. For $s \ge 0$, write $\mr_s = \mr_s(X_s)$. If $|\ma_{r+t}' \cup \mr_{r+t}| \le \rho$ and $|\mb_{r+t}'| \le \rho n$, then set
\[
	\ma_{r+t} = \ma_{r+t}' \cup \mr_{r+t} \cup \brb{o_{r+t}^1, ..., o_{r+t}^{a_{r+t}}}
\Quad{and}
	\mb_{r+t} = \mb_{r+t}' \cup \brb{c_{r+t}^1, ..., c_{r+t}^{b_{r+t}}},
\label{eq:set-def-good}
\nt
\]
where $a_{r+t}$ and $b_{r+t}$ are so that $|\ma_{r+t}| = \rho$ and $|\mb_{r+t}| = \rho n$; otherwise, set
\[
	\ma_{r+t} = \bra{o_{r+t}^1, ..., o_{r+t}^{\rho}}
\Quad{and}
	\mb_{r+t} = \bra{c_{r+t}^1, ..., c_{r+t}^{\rho n}}.
\label{eq:set-def-bad}
\nt
\]
Let $V$ be the first update time for the process $(\ma_{r+t} \cup \mb_{r+t})_{r>0}$:
\[
	V = \inf\brb{ r \ge 0 \mid \text{an edge of $\ma_{r+t} \cup \mb_{r+t}$ changes state at time $r+t$} } + t.
\]
Define $\mP_{r+t} = \ma_{r+t} \cup \mb_{r+t}$ for $r \in [0,V)$.

If at time $s$ we use \eqref{eq:set-def-good} to define $\ma_s$ and $\mb_s$, then we say that the set definitions \textit{succeeded} at time $s$, and write $\ms_s$ for this event; if we used \eqref{eq:set-def-bad}, then we say they have \textit{failed}.

Finally, for $t > 0$ define the event $\ms[0,t) = (\cap_{s<t} \ms_s) \cap \sg[0,t)$.
\end{defnt}

\begin{defnt}
\label{defn:update-times}
	Set $S = n\logn[8]$ and $U_0 = 0$. We say that $\mP$ \textit{updates} (at time $t$) when an edge in the set $\mP$ refreshes and changes state (at time $t$). We inductively define the sequence $U_1, U_2, ...$: for all~$k \ge 1$, let $V_k$ be the first time after $U_{k-1}$ that $\mP$ updates; set $U_k = V_k \wedge (U_{k-1}+S)$.
\end{defnt}

Work on the event $\ms[0,t)$, and fix $s < t$. Then $\ma_s$ is a collection of open edges of size $\rho$ and $\mb_s$ is a collection of closed edges of size $\rho n$. Hence the set $\mP_s = \ma_s \cup \mb_s$ updates at rate $\kappa \mu$ where
\[
	\kappa
=
	(1-p) \rho + p \rho n
=
	(1-\lambda/n) \rho + \lambda \rho
=
	(1 + \lambda - \lambda/n) \rho;
\Quad{note that}
	\kappa \asymp \rho = \logn[11].
\label{eq:kappa-def}
\nt
\]
By the memoryless property, $U_k - U_{k-1} \sim^\iid \me(\kappa\mu) \wedge S$. Observe also that the walk may only leave the giant when the set $\mP$ updates, in particular only at one of the times $U_1, U_2, ...$, but note that not all of these times are caused by updates: some are caused because of the threshold $S$.

We first look at the probability of the event $\{ \tau_\mg' \le U_k \} \cap \ms[0,U_k)$.

\begin{lem}
\label{res:leavegiant:Uk-U1}
	For all~$n$, all~$k$ and all~$(x_0,\eta_0)$ with $x_0 \in \mg_0$, we have
	\[
		\prb[x_0,\eta_0]{ \tau_\mg' \le U_k, \, \ms[0,U_k) }
	\le
		k \max_{\eta_0, x_0 \in \mg_0} \prb[x_0,\eta_0]{ \tau_\mg' = U_1, \, \ms[0,U_1) }.
	\]
\end{lem}

\begin{Proof}
	Consider any~$(x_0,\eta_0)$ satisfying $x_0 \in \mg_0$. By the union bound and the strong Markov property (applied at time $U_{j-1}$ for the $j$-th term of the sum), we find that
	\GAP{3}
	\begin{align*}
	&	\prb[x_0,\eta_0]{ \tau_\mg' \le U_k, \, \ms[0,U_k) }
	\le	\textstyle
	\sum_{j=1}^k
		\prb[x_0,\eta_0]{ \tau_\mg' = U_j, \, \ms[0,U_j) }
	\\&\gap
	\le
		\max_{\eta_0, x_0 \in \mg_0} \prb{ \tau_\mg' = U_1, \, \ms[0,U_1) }
	\cdot
	{ \textstyle \sum_{j=1}^k 
		\prb[x_0,\eta_0]{ \tau_\mg' > U_{j-1}, \, \ms[0,U_{j-1}) } }.
	\end{align*}
	Upper bounding the sum by $k$ completes the proof.
\end{Proof}

We now determine the \textit{uniform} mixing time of the static walk on a good giant. For a Markov chain $Z$ with transition matrix $P$ and invariant distribution $\pi$, the \textit{uniform mixing time} is
\[
	\tunif(\eps,Z)
=
	\inf\brb{ t \ge 0 \mid \max_{x,y} \absb{1 - p_t(x,y)/\pi(y)} \le \eps }.
\]

\begin{lem}
\label{res:leavegiant:unif-mix-static}
	Let $G$ be a good graph, and let $\mg$ be its (unique) giant. Consider the static walk, denoted $\tilde X$, on the giant. Write $\tunif(\eps,\tilde X)$ for the $\eps$-uniform mixing time of $\tilde X$ (on $\mg$). Then
	\[
		\tunif\rbb{\tfrac18,\tilde X} \lesssim n \logn[6].
	\]
\end{lem}

\begin{Proof}
To prove this lemma, we compare $\tilde X$ with a `sped-up' version. Consider a walk $Z$ on the giant $\mg$ of a good graph $\zeta$. Write $m = |\mg|$; so $m \asymp n$. Write $d(z)$ for the degree of $z$ in $\zeta$, and $d^*$ for the maximum degree; note that $d^* \le C_* \logn$. Associate to a vertex $z \in \mg$ the following set:
\[
	\mv_z = \mn_z \cup \{ 1,...,k_z \} \setminus \{z\}
\Quad{where}
	\text{$k_z$ is such that} \
	V_z = |\mv_z| = 2 C_* \logn,
\]
where $\mn_z = \bra{z' \in \mg \mid \zeta(z,z') = 1}$ is the (open) neighbourhood of $z$ (in $\zeta$). (This is possible since $d^* \le C_*\logn$.) Give $Z$ a rate 1 jump clock: when this clock rings, if $Z$ is at $z$ then a vertex $z'$ is chosen uniformly at random from $\mv_z$ and $Z$ moves (from $z$) to $z'$ if and only if the edge $(z,z')$ is open, ie $\zeta(z,z') = 1$.
Observe that $Z$ is the same as the static walk $\tilde X$, except that it is sped-up by a factor $n/(2C_*\logn)$. Hence the mixing times are in ratio $n/(2C_*\logn)$, for both total variation and uniform mixing. We now calculate the uniform mixing time $\tmix(\tfrac18,Z)$.

Since $V_z = |\mv_z| = 2C_*\logn \ge 2 d^*$ for all~$z \in \mg$, we see that the chain $Z$ is \textit{lazy} in the sense that if we discretise by its rate 1 jump clock then the resulting discrete-time chain is \textit{lazy}, ie $p(z,z) \ge \tfrac12$ for all~$z \in \mg$. Moreover,
\[
	\pi_Z(z) = 1/|\mg| = 1/m
\Quad{and}
	p(z,z') = \frac1{2C_*\logn} \oneb{ \zeta(z,z') = 1 }.
\]
Hence $Z$ is reversible. It is then known that
\[
	\tunif(\epsilon,Z)
\lesssim
	\Phi_*^{-2} \rbb{ \log(1/\pimin) + \logeps },
\]
where $\Phi_* = \inf\bra{ \Phi_S \mid \pi_Z(S) \le \tfrac12 }$ and $\Phi_S = \sum_{x \in A, y \in B} \pi_Z(x) p_Z(x,y) / \pi_Z(S)$; for a proof of this, see \cite{JS:bottleneck-mixing}.
For any set $S \subseteq \mg$, we have
\[
	\Phi_S
=
	\frac1{2C_*\logn} \cdot \frac{|\partial S|}{|S|}
\ge
	\frac1{2C_*\logn} \cdot \frac{|\partial S|}{d(S)}
=
	\frac1{2C_*\logn} \cdot \Phi_S',
\]
where the prime ($\prime$) denotes that we are considering the corresponding quantity for the nearest-neighbour discrete-time random walk. But we know that $\Phi_*' \gtrsim \logn[-2]$ since the graph is good, and hence $\Phi_* \gtrsim \logn[-3]$. Hence
\[
	\tunif(\tfrac18,Z) \lesssim \logn[7],
\Quad{and hence}
	\tunif\rbb{\tfrac18,\tilde X} \lesssim n \logn[6].
\qedhere
\]
\end{Proof}

We now use this mixing of the static walk along with our static-dynamic coupling to determine where the dynamic walk is at the update times of $\mP$.

\begin{lem}
\label{res:leavegiant:U1-q}
	There exists a constant $C$, so that, for all~$M$, all~$n$ sufficiently large, we have
	\[
		\max_{\eta_0, x_0 \in \mg_0}
		\prb[x_0,\eta_0]{ \tau_\mg' = U_1, \, \ms[0,U_1) }
	\le
		C \mu S \log_{(M)}\! n.
	\]
\end{lem}

\begin{Proof}
For this whole proof, we only consider the first update time $U_1$; as such, we drop the 1 from the subscript, just writing $U$. Also, we write $\tilde X$ for the static walk on $\eta_0$ (as above).

For the walk to leave the giant, we need the time $U $ to be triggered by an update to $\mP$, ie we need $U < S$.
If this is the case, then the walk leaves the giant if and only if the update was caused by the closing of one of the removal edges which, given $R_{U^-}(X_U)$, has probability $R_{U^-}(X_U)/\kappa$;
write $R_U = R_{U^-}(X_U)$.
Hence
\[
	\prb[x_0,\eta_0]{ \tau_\mg' = U, \, \ms[0,U) }
=
	\tfrac1\kappa
	\exb[x_0,\eta_0]{ R_U \one{U < S} \one{\sg[0,U)} }.
\label{eq:leavegiant:prob-R}
\nt
\]
We now set $T = n \logn[7]$; so $S = T \logn \gg T$.
We decompose according to $\{U < T\}$ or $\{U \ge T\}$. When $U < T$ we use the trivial bound $R_U \le C_* \logn$ (which holds whenever the graph is good):
\[
	\exb[x_0,\eta_0]{ R_U \one{U < T} \one{ \sg[0,U) } }
\le
	C_* \logn \cdot \prb[x_0,\eta_0]{ U < T, \, \sg[0,U) }.
\]
Since $\mw_0^{k+1} \subseteq \mw_0^k$, for all~$M$ and all~$n$ sufficiently large, we have
\GAP{5}
\begin{align*}
&	\exb[x_0,\eta_0]{ R_U \one{T \le U < S} \one{ \sg[0,U) } }
=
	\exb[x_0,\eta_0]{ R_U \one{T \le U < S} \oneb{ X_U \in \mw_0^M } \one{ \sg[0,U) } }
\\&\gap \textstyle
+	\sum_{k=1}^{M-1}
	\exb[x_0,\eta_0]{ R_U \one{T \le U < S} \oneb{ X_U \in \mw_0^k \setminus \mw_0^{k+1} } \one{ \sg[0,U) } }.
\label{eq:leavegiant:W-expansion}
\nt
\end{align*}
When $X_U \in \mw_0^k$, we have (by definition) $R_U \le C_* \log_{(k)}\! n$. Hence we have
\GAP{1}
\begin{align*}
&	\exb[x_0,\eta_0]{ R_U \one{T \le U < S} \one{ \sg[0,U) } }
\\&\gap
\le
	C_* \prb{ U < S, \, \sg[0,U) } \cdot \rbb{ \textstyle
		\log_{(M)}\! n + \sum_{k=1}^{M-1} \log_{(k)}\! n \cdot \prb{ X_U \notin \mw_0^{k+1} \mid T \le U < S, \, \sg[0,U) } }.
\end{align*}

What is crucial is that, on the event that the graph is good, the update \emph{times} are independent of the evolution of the walk: since $\mP$ always, regardless of the number of edges seen by the walker, contains precisely $\rho$ open edges and $\rho n$ closed edges, the update rate is always $\kappa\mu$. Thus an equivalent way of realising $(\mP_t)_{t \in [0,U)}$ is the following.
	Define the processes $(\ma_r)_{r\ge0}$ and $(\mb_r)_{r\ge0}$ as in (\ref{eq:set-def-good}, \ref{eq:set-def-bad}), taking $t = 0$.
	Then sample independently $V \sim \me(\kappa\mu)$.
	At time $V$ with probability $q = (1-p)/(1-p+\lambda)$ choose an edge uniformly at random from $\ma_V$ and change its state from open to closed, and with probability $1-q$ choose an edge uniformly at random from $\mb_V$ and change its state from closed to open.
	Then set $\mP_r = \ma_r \cup \mb_r$ for all~$r \in [0,V)$.
	Finally, set $U = V \wedge S$.

It remains to calculate this final probability, of $X_U \notin \mw_0^{k+1}$.
We want to couple $X_U$ with $\tilde X_U$, as we can then apply the (uniform) mixing result \Cref{res:leavegiant:unif-mix-static} to obtain good control over its location.
However, we can only do this under certain conditions; sufficient conditions are that none of the edges of $\mP_{U^-}$ have changed \emph{throughout the entire interval $[0,U)$}.
(Note that an edge could change state \emph{before} it is added to the set process $\mP$.)
Write $\mc$ for this sufficient condition.
Then
\begin{align*}
&	\prb[x_0,\eta_0]{ X_U \notin \mw_0^{k+1} \mid T \le U < S, \, \sg[0,U) }
\\&\qquad
\le
	\prb[x_0,\eta_0]{ \tilde X_U \notin \mw_0^{k+1} \mid T \le U < S }
+	\prb[x_0,\eta_0]{ \mc^c \mid T \le U < S, \, \sg[0,U) },
\end{align*}
since conditioning on $\sg[0,U)$ has no effect on the static walk.

Since $T \gg n \logn[6]$, which is the uniform mixing time of the static walk on a good giant (\Cref{res:leavegiant:unif-mix-static}), if $U \ge T$ then $\tilde X_U$ has (uniformly) mixed and so, since the invariant distribution of the static walk is uniform (on the giant), for all~$k \le M$, we have
\[
	\prb[x_0,\eta_0]{ \tilde X_U \notin \mw_0^{k+1} \mid T \le U < S }
\le
	\tfrac32 \, \absb{ \mg_0 \setminus \mw_0^{k+1} } / \abs{\mg_0}
\le
	\tfrac32 (\log_{(k)}\! n)^{-4},
\]
with the final inequality holding by definition of a good graph; here we have used crucially that the path $(\tilde X_t)_{t\ge0}$ is independent of $U$.
Also, since the update rate of $\mP$ is always $\kappa \mu \asymp \mu \logn[11]$ and we run for time $U \le S = n \logn[8]$, we find that
\begin{align*}
	\prb[x_0,\eta_0]{ \mc^c \mid T \le U < S, \, \sg[0,U) }
&
\le
\frac{
	\prb{ \me(\kappa\mu) < S }
}{
	\prb[x_0,\eta_0]{ T \le U < S, \, \sg[0,U) }
}
\le
	\kappa S \mu \cdot \rbb{1 + \oh1 }
\ll
	(\log_{(k)}\! n)^{-4},
\end{align*}
by the assumption that $\mu \le \logn[-20]/n$ and the fact that $\eta_0 \in \sh$.
Together, these give
\[
	\prb{ X_U \notin \mw_0^{k+1} \mid T \le U < S, \, \sg[0,U) }
\le
	2 (\log_{(k)}\! n)^{-4}.
\label{eq:leavegiant:bad-vts}
\nt
\]

Also, for any~$s \ge 0$, we have
\[
	\prb[x_0,\eta_0]{ U < s, \, \sg[0,U) }
\le
	\prb{ \me(\kappa\mu) \le s }
=
	1 - e^{-\kappa \mu s}
\le
	\kappa \mu s.
\]
Hence combining these inequalities, for all~$M$ and all~$n$ for sufficiently large, we have
\[
	\prb[x_0,\eta_0]{ \tau_\mg' = U_1, \, \ms[0,U_1) }
\le
	2 C_* \mu S \log_{(M)}\! n.
\qedhere
\]
\end{Proof}

Let $K$ be the (random) index given by $U_K \le t < U_{K+1}$.
{Note that
\[
	\prb[x_0,\eta_0]{ \tau_\mg' \le t, \, K \le k-1, \, \ms[0,U_k) }
\le
	\prb[x_0,\eta_0]{ \tau_\mg' \le U_k, \, \ms[0,U_k) },
\]
by monotonicity of $t \mapsto \{\tau_\mg' \le t\}$, along with the fact that $t < U_{K+1}$.}
Hence
\begin{align*}
	\prb[x_0,\eta_0]{ \tau_\mg' \le t }
&\le
	\prb[x_0,\eta_0]{ \tau_\mg' \le U_k, \, \ms[0,U_k) }
+	\prb[x_0,\eta_0]{ K \ge k, \, \ms[0,U_k) }
\\&\qquad
+	\prb[x_0,\eta_0]{ \ms[0,U_k)^c, \, \sg[0,U_k) }
+	\prb[x_0,\eta_0]{ \sg[0,U_k)^c }.
\end{align*}
We have already dealt with the first term in the previous lemmas; we now just need to show that the three `remainder' terms are sufficiently small. We do this now.

\begin{lem}
\label{res:leavegiant:K-tail}
	For all $n$ sufficiently large, all~$t$ with $\mu t \ge \logn[-5]$ and all~$(x_0,\eta_0)$ with $x_0 \in \mg_0$, for $k = \ceil{5t/S}$, we have
	\[
		\prb[x_0,\eta_0]{ K \ge k, \, \ms[0,U_k) }
	\le
		n^{-5}.
	\]
\end{lem}

\begin{Proof}
If we were to not have the thresholding by $S$, then we would have $K \sim \Po(\kappa \mu t)$. However, we do have the thresholding.
Set $\tilde U_0 = 0$, and inductively define $\tilde U_j$, for $j = 1,2,...$, by
\[
	\tilde U_j - \tilde U_{j-1} = S \cdot \oneb{ U_j - U_{j-1} = S },
\Quad{and set}
	\tilde K = \inf\brb{ k \ge 0 \mid \tilde U_k \le t < \tilde U_{k+1} }.
\]
We have $\tilde U_j \le U_j$ for all~$j \ge 0$, and thus $\tilde K \ge K$.

Recall that when the set definitions succeed, $\{U_j - U_{j-1}\}_{j\ge1}$ is a collection of iid random variables, and are independent of the starting point $(x_0,\eta_0)$.
{Recalling from \Cref{defn:update-times} that $S = n\logn[8]$ and from \eqref{eq:kappa-def} that $\kappa \asymp \logn[11]$, note that}
\[
	\prb{ \me(\kappa\mu) \ge S } = e^{-\kappa\mu S}
=
	1 - \oh1
\ge
	\tfrac12,
\]
by the assumption $\mu n \ll \logn[-19]$.
Also let us write $k' = \ceil{t/S}$; since $\mu t \ge \logn[-5]$ and $\mu n \ll \logn[-14]$, we have $t/S \gg \logn$, and so $k' \gg 1$ and $k' \le 2t/S$.
Then, on the event that the set definitions succeed, we have $\tilde K \preccurlyeq \Po(4t/S)$, since $\pr{ \me(\kappa\mu) \ge S } \ge \tfrac12$. Hence
\[
	\prb[x_0,\eta_0]{ K \ge k, \, \ms[0,U_k) }
\le
	\prb{ \Po(4t/S) \ge 5t/S }
\le
	\expb{ - \tfrac1{10} t/S },
\]
by Poisson concentration.
Since $t/S \gg \logn$, we deduce our lemma.
\end{Proof}


\begin{lem}
\label{res:leavegiant:setdefs-general}
	For all~$n$ sufficiently large, all~$k$ and all~$(x_0,\eta_0)$ with $x_0 \in \mg_0$, we have
	\[
		\prb[x_0,\eta_0]{ \ms[0,U_k)^c, \, \sg[0,U_k) }
	\le
		k \cdot \expb{ - \tfrac13 C_* \logn[9] }.
	\]
\end{lem}

\begin{Proof}
By the union bound, we have
\[
	\prb[x_0,\eta_0]{ \ms[0,U_k)^c, \, \sg[0,U_k) }
\le
	k \cdot \max_{x_0, \eta_0} \prb[x_0,\eta_0]{ \ms[0,U_1)^c, \, \sg[0,U_1) }.
\]
Since we work on the event that the graph is good, we have at most $C_* \logn$ removal edges for each vertex; we also have that there are $\Om{n}$ open edges and $\Om{n^2}$ closed edges. Hence the only part that can `go wrong' in the definitions is if the number of open or closed edges seen since the last update is too high. However, the maximum degree is at most $C_* \logn$ and we walk for a time at most $S = n \logn[8]$, so Poisson concentration will tell us that we do not see too many.

Consider a static walk $\tilde X$ on a good graph $\eta_0$, starting from $x_0 \in \mg_0$ and run for a time $S = n \logn[8]$. Let $\alpha$ be the number of open edges seen in this time, and $\beta$ the number of closed. Write $N$ for the number of steps taken; by Poisson thinning, we have $N \preccurlyeq \Po(C_* S \logn/n)$, and $S \logn/n = \logn[9]$. On the event $N \le 2 C_* \logn[9]$, we have $\alpha \le 2 C_*^2 \logn[10] \ll \rho$ and $\beta \le 2 C_* n \logn[9] \le \rho n$, as required for the set definitions to succeed. Hence
\[
	\prb[x_0,\eta_0]{ \ms[0,U_1)^c, \, \sg[0,U_1) }
\le
	\prb{ \Po(C_* \logn[9]) > 2 C_* \logn[9] }
\le
	\expb{ - \tfrac13 C_* \logn[9] },
\]
by Poisson concentration.
The result now follows from the union bound given above.
\end{Proof}


\begin{cor}
	For all~$n$ sufficiently large, all~$t$ with $\mu t \ge \logn[-5]$ and all~$(x_0,\eta_0)$ with $x_0 \in \mg_0$, for $\mu \ge n^{-8}$ and $k = \ceil{5t/S}$, we have
	\[
		\prb[x_0,\eta_0]{ \ms[0,U_k)^c, \, \sg[0,U_k) }
	\le
		\mu t \cdot n^{-5}.
	\]
\end{cor}

\begin{Proof}
	As in \Cref{res:leavegiant:K-tail} for $\ceil{t/S}$, we have $k \le 6 t / S$. Hence
	\[
		k \cdot \expb{ - \tfrac13 C_* \logn[9] }
	\le
		6 S^{-1} \cdot \mu t \cdot n^8 \expb{ - \tfrac13 C_* \logn[9] }.
	\le
		\mu t \cdot n^{-5}.
	\qedhere
	\]
\end{Proof}

We now have all the ingredients to prove \Cref{res:leavegiant} for the case $\mu \ge n^{-8}$.

\begin{Proof}[Proof of \Cref{res:leavegiant} \textnormal{(\emph{when $\mu \ge n^{-8}$})}]
Fix $M$.
Combining the above results, we have, for $k = \ceil{5t/S}$, recalling that $U_k \le kS \le 6t$ for the times $t$ we are considering, that
\begin{align*}
	\prb[x_0,\eta_0]{ \tau_\mg' \le t }
&\le
	\prb[x_0,\eta_0]{ \tau_\mg' \le U_k, \, \sg[0,U_k) }
+	\prb[x_0,\eta_0]{ K \ge \ceil{5t/S}, \, \ms[0,U_k) }
\\&\quad
+	\prb[x_0,\eta_0]{ \ms[0,U_k)^c, \, \sg[0,U_k) }
+	\prb[x_0,\eta_0]{ \sg[0,6t]^c }
\label{eq:leavegiant:final}
\nt
\\&\le
	C \mu t \log_{(M)}\! n
+	n^{-5}
+	\mu t n^{-5}
+	n^{-1}
\le
	2 C \mu t \log_{(M)}\! n
\end{align*}
since $\eta_0 \in \sh$, $\mu \ge n^{-8}$ and $\logn[-5] \le \mu t \le \tfrac1{10}$.
\end{Proof}

It remains to prove the proposition in the case $\mu \le n^{-8}$. In this case, `almost always' the static walk mixes on the entire giant before any of the graph even refreshes; this will make this proof easier. The general idea will be very similar, particularly to \Cref{res:leavegiant:U1-q}.

\begin{Proof}[Proof of \Cref{res:leavegiant} \textnormal{(\emph{when $\mu \le n^{-8}$})}]
Fix $M$.
For this proof, let $U_1, U_2, ...$ be the refresh times of the graph; let $U_0 = 0$. Note then that $U_j-U_{j-1} \sim^\iid \me(\mu N)$ where $N = \binom{n}{2} \le n^2$.

We are now interested in the probability that $\tau_\mg' = U_1$; as previously, drop the subscript 1.
We have $U \sim \me(\mu N)$, independent of $X$.


Suppose $\eta_0 \in \sg$ and $x_0 \in \mg_0$. Then, similarly to in \eqref{eq:leavegiant:prob-R}, we have
\begin{align*}
	\prb[x_0,\eta_0]{ \tau_\mg' = U, \, \sg[0,U) }
=
	\tfrac1N \exb[x_0,\eta_0]{ R_U \one{U \ge n^2} \one{\sg[0,U)} }
+	\prb[x_0,\eta_0]{ U \le n^{-2} }.
\end{align*}
We know that $\pr{U \le n^2} = \pr{\me(\mu N) \le n^2} \le \mu n^2 N \le n^{-4}$. Similarly to in \eqref{eq:leavegiant:W-expansion}, we have
\GAP{5}
\begin{align*}
&	\exb[x_0,\eta_0]{ R_U \one{U \ge n^2} \one{ \sg[0,U) } }
=
	\exb[x_0,\eta_0]{ R_U \one{U \ge n^2} \oneb{ X_U \in \mw_0^M } \one{ \sg[0,U) } }
\\&\gap \textstyle
+	\sum_{k=1}^{M-1}
	\exb[x_0,\eta_0]{ R_U \one{U \ge n^2} \oneb{ X_U \in \mw_0^k \setminus \mw_0^{k+1} } \one{ \sg[0,U) } }.
\end{align*}
When $X_U \in \mw_0^k$, we have (by definition) $R_U \le C_* \log_{(k)}\! n$. Hence we have
\GAP{1}
\begin{align*}
	\exb[x_0,\eta_0]{ R_U \one{U \ge n^2} \one{ \sg[0,U) } }
\le \textstyle
	\log_{(M)}\! n + \sum_{k=1}^{M-1} \log_{(k)}\! n \cdot \prb{ X_U \notin \mw_0^{k+1} \mid U \ge n^2, \, \sg[0,U) }.
\end{align*}
Using our static-dynamic coupling, we may couple $X_t = \tilde X_t$ for all $t \le U = U_1$, where $\tilde X$ is the static walk, since the $U_j$ are the refresh times of the entire graph.
Hence, as in \eqref{eq:leavegiant:bad-vts}, but without needing to consider the condition $\mc$, we use the (uniform) mixing of the static walk to obtain
\[
	\prb[x_0,\eta_0]{ X_U \notin \mw_0^{k+1} \mid U \ge n^2, \, \sg[0,U) }
\le
	2 \, \absb{ \mg_0 \setminus \mw_0^{k+1} } / \abs{\mg_0}
\le
	2 (\log_{(k)}\! n)^{-4}.
\]
Hence combining these inequalities, for all~$M$ and all~$n$ sufficiently large, we have
\[
	\max_{x_0,\eta_0}
	\prb[x_0,\eta_0]{ \tau_\mg' = U_1, \, \sg[0,U_1) }
\le
	2 C_* \tfrac1N \log_{(M)}\! n.
\]

As in \Cref{res:leavegiant:Uk-U1} (except replacing $\sg$ by $\ms$), for all~$k \in \bn$, we have
\[
	\prb[x_0,\eta_0]{ \tau_\mg' \le U_k, \, \sg[0,U_k) }
\le
	k \max_{x_0,\eta_0}
	\prb[x_0,\eta_0]{ \tau_\mg' = U_1, \, \sg[0,U_1) }
\le
	2 C_* k \tfrac1N \log_{(M)}\! n.
\]

Observe that $U_k \sim \Gamma(k,\mu N)$; let $K$ be the (random) index given by $U_K \le t < U_{K+1}$, and observe then that $K \sim \Po(\mu t N)$. Set $k = \ceil{2 \mu t N}$. By the same arguments as used in \eqref{eq:leavegiant:final} (except without the `set-definitions' term) we have
\begin{align*}
	\prb[x_0,\eta_0]{ \tau_\mg' \le t }
&
\le
	\prb[x_0,\eta_0]{ \tau_\mg' \le U_k, \, \sg[0,U_k) }
+	\prb{ K \ge 2 \mu t N }
+	\prb{U_k \ge 3k/(\mu N)}
+	\prb[\eta_0]{ \sg[0,7t]^c }
\\&
\le
	2 C_* \mu t \log_{(M)}\! n
+	\expb{- \tfrac13 \mu t N}
+	\expb{- 2 \mu t N}
+	\prb[x_0,\eta_0]{ \sg[0,7t]^c }.
\end{align*}
We now recall that we restricted consideration of $t$ to satisfy $\logn[-5] \le \mu t \le \tfrac1{10}$; note then that $7t \le 1/\mu$. We also consider only $\eta_0 \in \sh$; as such, the graph remainder term in the final line above is at most $1/n$. Since $N \asymp n^2$, we see that the first term dominates, leaving us with
\[
	\prb[x_0,\eta_0]{ \tau_\mg' \le t }
\le
	C \mu t \log_{(M)}\! n
\quad
	\text{for a constant $C$}.
\qedhere
\]
\end{Proof}

\section{Isolation Times}
\label{sec:isol-times}

In this section we prove two main results on isolation times. They will involve, respectively, a single random walker on a dynamical environment and two independent random walkers on the same dynamical environment.
When considering just one walk, we write $\tauisol$ for the isolation time; eg for a dynamical percolation system $(Z,\zeta)$ we write
\[
	\tauisol^Z = \inf\brb{ t \ge 0 \mid d^\zeta_t(Z_t) = 0 }.
\]
When the context is clear, we omit the superscript, just writing $\tauisol$; similarly, when the context is clear we write $d$ for the degree, rather than $d^\zeta$. When we consider two walks, $X$ and $Y$, on the same system, $\eta$, we use superscript $X$ or $Y$ to indicate which walk we are referring to: define
\[
	\tauisol^X = \inf\brb{ t \ge 0 \mid d_t(X_t) = 0 }
\Quad{and}
	\tauisol^Y = \inf\brb{ t \ge 0 \mid d_t(Y_t) = 0 }.
\]
{Recall from \Cref{defn:very-good-graph} that the event $\Hgood{0,t}$ guarantees that the graph is good up until time~$t$.}

\begin{thm}[Single-Walker Isolation Time]
\label{res:1:tauX>t;H}
	For all~$M \in \bn$, all~$n$ sufficiently large and all pairs $(x_0,\eta_0)$, we have
	\[
		\prb[x_0,\eta_0]{ \tauisol >  t, \, \Hgood{0,t} }
	\le
		2 \expb{ - \mu t / \log_{(M)}\! n }.
	\]
	Moreover, if $\mu t \ge 3$, then we may remove the pre-factor of 2.
\end{thm}

Once we have proved this, we shall be able to use a type of concentration result to prove a bound on the isolation time of two independent random walkers on the same environment. When we are considering this, we write $\bp_{x_0,y_0,\eta_0}$ for the measure. For two walks $X$ and $Y$ on the same (dynamical) environment $\eta$, write
\[
	\tau = \inf\brb{ t \ge 0 \mid d_t(X_t) = 0 = d_t(Y_t) }.
\]

\begin{thm}[Dual-Walker Isolation Time]
\label{res:2:tau>t;H}
	For all~$M \in \bn$, all~$n$ sufficiently large and all triples $(x_0,y_0,\eta_0)$, we have
	\[
		\prb[x_0,y_0,\eta_0]{ \tau > t, \, \Hgood{0,t} }
	\le
		2 \expb{ - \mu t / \log_{(M)}\! n }.
	\]
\end{thm}

\subsection{Single-Walker Isolation Time}

In this section we prove \Cref{res:1:tauX>t;H} on the isolation time of a walk $X$ on a dynamical environment $\eta$.
In order to find the isolation time, we wait until the walk joins the giant and then look at becoming isolated from there. It is easier to consider the giant, rather than subcritical components, because we are able to use concentration results on the structure of the giant.

We first state the proposition on isolation from the giant, and then show how to conclude \Cref{res:1:tauX>t;H} from it; we then prove the proposition to finish. Throughout, $M$ is a positive integer.

\begin{prop}[Isolation from the Giant]
\label{res:sup:tauX<t}
	There exists a positive constant $c$ so that, for all~$M$, all~$n$ sufficiently large and all~$(x_0,\eta_0)$ with $\eta_0 \in \sh$ and $x_0 \in \mg_0$, we have
	\[
		\prbb[x_0,\eta_0]{ \tauisol \le \frac1{\mu \log_{(M)}\! n} }
	\ge
		c \cdot \frac1{\log_{(M)}\! n}.
	\]
\end{prop}

\begin{Proof}[Proof of \Cref{res:1:tauX>t;H}]
Observe that this trivially holds (for all~$n$ large enough) if $\mu t \le 3$. By monotonicity, replacing $M$ by $M-1$, it suffices to prove an upper bound of $\exp{ - c \mu t / \log_{(M)}\! n }$ for a positive constant $c$ when $\mu t \ge 3$.

Fix $M$.
For this proof, rescale time so that $\mu = 1$.
We prove this theorem by performing independent experiments.
Note that if $x_0 \in \mg_0$ then $\tau_\mg = 0$, and otherwise we apply \Cref{res:hitgiant}. By direct calculation, we have
\GAP{3}
\begin{align*}
&	\prb[x_0,\eta_0]{ \tauisol \le 2, \, \Hgood{0,1} }
\ge
	\prb[x_0,\eta_0]{ \tauisol \le 2, \, \tau_\mg \le 1, \, \Hgood{0,1} }
\\&\gap
\ge
	\prb[x_0,\eta_0]{ \tauisol - \tau_\mg \le 1/\log_{(M)}\! n, \, \tau_\mg \le 1, \Hgood{0,1} }
\\&\gap
= \textstyle
	\sum_{x_0',\eta_0'}
	\prb[x_0',\eta_0']{ \tauisol \le 1/\log_{(M)}\! n }
\cdot
	\prb[x_0,\eta_0]{ X(\tau_\mg) = x_0', \eta(\tau_\mg) = \eta_0', \, \tau_\mg \le 1, \, \Hgood{0,1} }
\\&\gap
\ge
	c (\log_{(M)}\! n)^{-1} \cdot	\rbB{
		\prb[x_0,\eta_0]{ \tau_\mg \le 1 } - \prb[x_0,\eta_0]{ \Hbad{0,1} } },
\end{align*}
for a positive constant $c$,
where for the final inequality we used that on the event $\Hgood{0,1}$ we have $\eta_0' \in \sh$, and hence we may apply \Cref{res:sup:tauX<t}. Now applying \Cref{res:hitgiant}, we obtain
\[
	\prb[x_0,\eta_0]{ \tauisol \le 2, \, \Hgood{0,1} }
\ge
	\tfrac12 c c_1 / \log_{(M)}\! n
-	\prb[x_0,\eta_0]{ \Hbad{0,1} },
\]
with the positive constant $c_1$ coming from \Cref{res:hitgiant}, noting that $\prb[x_0,\eta_0]{ \Gbad{0,1} } = \oh1$ since $\eta_0 \in \sh$. Rearranging this, we obtain, for a positive constant $c$, that
\[
	\prb[x_0,\eta_0]{ \tauisol > 2, \, \Hgood{0,2} }
\le
	\prb[x_0,\eta_0]{ \tauisol > 2, \, \Hgood{0,1} }
\le
	\expb{ - c / \log_{(M)}\! n}.
\]
Hence, for any~$k \in \bn$, applying the strong Markov property ($k-1$ times), we obtain
\[
	\prb[x_0,\eta_0]{ \tauisol > 2k, \, \Hgood{0,2k} }
\le
	\max_{x_0',\eta_0'}
	\prb[x_0',\eta_0']{ \tauisol > 2, \, \Hgood{0,2} }^k
\le
	\expb{ - c k / \log_{(M)}\! n}.
\]
This completes the proof.
	%
\end{Proof}

It remains to prove \Cref{res:sup:tauX<t}. We do this via a sequence of lemmas, using the following rough methodology. Observe that the \emph{only} way for the walk to become isolated is to be at a degree~1 vertex and for the one open edge to close before any closed incident edges open or the walker leaves the vertex. This motivates looking at the rate at which the walk hits degree~1 vertices.

Since the walk is on a dynamically evolving graph, even though when we require the graph to be good this includes that the giant has a lot of degree~1 vertices, the location of these degree~1 vertices is changing. This makes using averaging properties (like a law of large numbers) difficult. However, since we take steps at rate at least $1/n$ (when non-isolated) and $\mu n \ll 1/\logn$ (the order of the maximum degree), we see that the vast majority of the time the walker takes a step before any edge incident to its location changes state. This motivates looking at the rate at which a walk (with the same walk-dynamics) hits degree~1 vertices \emph{on a static (good) graph}, and then relating this quantity to the relevant quantity for the walk on the dynamic graph. In \S\ref{sec:hit/exit-giant} we referred to this as the \textit{static walk} and the original as the \textit{dynamic walk}, denoting them by $\tilde{X}$ and $X$, respectively.

With this motivation in mind, we first collect some results regarding a walk with our dynamics on a static graph. To do this, we use a Chernoff-style bound on the number of visits to a set, which is due to Gillman \cite{G:gillman-chernoff}. It applies to \textit{discrete-time} random walks. We do not apply it to a discretisation of our continuous chain, but to the jump chain of the walk (on a static graph). We state it in a general form; an even more general form is given in \cite[Theorem 2.1]{G:gillman-chernoff}.

\begin{thm}[Gillman \cite{G:gillman-chernoff}]
\label{res:sup:gillman}
	Consider the discrete-time random walk on a weighted, connected graph $G = (V,E)$ with any initial distribution. Let $\pi$ be the unique invariant distribution. Let $A \subseteq V$, and let $N_m$ be the number of visits to $A$ in $m$ steps. Write $\gamma$ for the spectral gap. Then
	\[
		\prb{ \absb{ N_m - m \pi(A) } \ge R }
	\le
		3 \pimin^{-1/2} \expb{-\tfrac1{20} \gamma R^2 / m}
	\Quad{for any}
		\eps \in [0,m].
	\]
\end{thm}

We now apply this to a walk on a good (static) giant.


\begin{lem}
\label{res:sup:hit-deg1_gillman}
	Consider the discrete-time nearest-neighbour simple random walk on a graph $G$, and write $N_m$ for the number of visits to the set of degree~1 vertices in $m$ steps.
	There exists a positive constant $c$ so that, for all~$n$ sufficiently large and all~$m \ge \logn[6]$, if the graph is good, ie $G \in \sg$, and the walk starts from its giant, then we have
	\[
		\prb{N_m \le c m} \le n^{-1}.
	\]
\end{lem}

\begin{Proof}
Note that the invariant measure of this walk, which we denote $\pi'$, is given by $\pi_i = d_i/d_\mg$, where $d_\mg = \sum_{i \in \mg} d_i$. Since $d_i \ge 1$ for all~$i \in \mg$, we have $\pimin' \ge 1/d_\mg$. Now, trivially we have that $d_\mg \le d_G$, where $d_G = \sum_{i \in G} d_i$, and $d_G \le 2 C_* n$ by Definition \ref{defn:good-graph}\ref{gg:num-edges}. Hence $1/\pimin' \le d_\mg \le 2 C_* n$.

Let $A = \bra{ x \in \mg \mid d(x) = 1 }$ be the set of degree~1 vertices in the giant. Definition \ref{defn:good-graph}\ref{gg:num-deg1-giant} tells us that $|A| \ge c_* n$. Thus, since $d(x) \ge 1$ for all~$x \in \mg$ and $d_\mg \le 2C_*n$, we have that $\pi'(A) \ge c_*/(2C_*)$; let $c = c_*/(4C_*)$ so that $\pi'(A) \ge 2c$.

Recall from \Cref{defn:good-graph}\ref{gg:giant-expansion} that the spectral gap $\gamma$ of a good giant satisfies $\gamma \ge c_* \logn[-4]$. We now take $R = \tfrac12 \pi'(A) m \le m$ in \Cref{res:sup:gillman} to obtain
\[
	\prb{ \absb{ N_m - m \pi'(A) } \ge \tfrac12 \pi'(A) m }
\le
	3 \sqrt{2 c_*} \cdot n^{-1/2} \expb{ - \tfrac1{20} c_* \logn[-4] \cdot \tfrac14 \pi'(A)^2 m }.
\]
Since $\pi'(A) \ge 2c$, taking $m \ge \logn[6]$ gives super-polynomial decay, completing the proof.
\end{Proof}

We now make rigorous the motivation given at the start of this section in the following lemma.

\begin{lem}
\label{res:sup:tauX<s}
	There exists a positive constant $q$ so that, for all~$n$ sufficiently large and all~$(x_0,\eta_0)$ with $\eta_0 \in \sg$ and $x_0 \in \mg_0$, for $s = n \logn[6]$, we have
	\[
		\prb[x_0,\eta_0]{ \tauisol \le s }
	\ge
		q \mu s.
	\]
\end{lem}

\begin{rmkt*}
	Observe that, as a best-case scenario, if the walker were always at a degree~1 vertex until it becomes isolated, then the isolation time would simply be the time it takes for that one edge to close, which is $\me( \mu(1-p) )$. Thus, for any~$(x_0,\eta_0)$ with $d_0(x_0) \neq 0$, we have
	\[
		\prb[x_0,\eta_0]{ \tauisol \le s } \le \mu s / (1-p).
	\]
	Hence, for this $s$, such a result as \Cref{res:sup:tauX<s} is best-possible up to constants.
\end{rmkt*}

\begin{Proof}[Proof of \Cref{res:sup:tauX<s}]
Fix a pair $(x_0,\eta_0)$ with $\eta_0 \in \sg$ and $x_0 \in \mg_0$; for this proof, drop it from the notation, writing $\pr{\cdot}$ in place of $\pr[x_0,\eta_0]{\cdot}$.

\Cref{res:sup:hit-deg1_gillman} tells us the rate at which the static walk hits degree~1 vertices (with high probability). In order to transfer this result to our dynamic walk, we define a coupling between the two walks; this was given in \S\ref{sec:hit/exit-giant}, but we recall it precisely here. Write $\tilde{X}$ for the static walk, walking on the static graph $\eta_0$. Set $X_0 = \tilde{X}_0 = x_0$. Give $X$ and $\tilde{X}$ the same jump clock. When the clock rings, at time $t$ say, both walks choose the same vertex; $\tilde{X}$ performs the jump if and only if the edge is present in $\eta_0$, while $X$ performs the jump if and only if the edge is present in $\eta_t$. We call this the \textit{static-dynamic coupling}.

We now define $\tilde{T}_i$ to be the $i$-th time that (the static walk) $\tilde{X}$ hits a degree~1 vertex: set $\tilde{T}_0 = \tilde{T}_0' = 0$ and define inductively, for $i \ge 1$,
\[
	\tilde{T}_i = \inf\brb{ t \ge \tilde{T}_{i-1}' \midb d_t\rbr{\tilde{X}_t} = 1 }
\Quad{and}
	\tilde{T}_i' = \inf\brb{ t \ge \tilde{T}_i \midb \tilde{X}(t) \neq \tilde{X}\rbr{\tilde{T}_i} }.
\]
Since the jump rate of $\tilde{X}$ is always at least $1/n$, by standard Poisson concentration it takes at least $s/(2n)$ steps in time $s-n$ with probability $1 - \oh{1}$. Along with \Cref{res:sup:hit-deg1_gillman} this says that
\[
	\prb{ \tilde{T}_k \le s-n }
=
	1 - \oh{1}
\Quad{for}
	k = s/(8n).
\label{eq:sup:Tk<s;gg}
\nt
\]
For $i \ge 1$ define the event that all the (open or closed) edges incident to a vertex that the static walk visited remain in the same state between visits to degree~1 vertices:
\[
	\se_i
=
	\brb{ \text{neighbourhood of path of static walk did not change in $[\tilde{T}_{i-1}, \tilde{T}_i]$} }.
\]
Similarly, for $u \ge v \ge 0$ define
\[
	\se_{u,v} = \brb{ \text{neighbourhood of path of static walk did not change in $[u,v]$} }.
\]
Note that, by definition, on the event $\{\tilde{T}_k \le s\}$ we have $\se_1 \caps \se_k \supseteq \se_{0,s}$.

Write $N$ for the number of steps taken by $\tilde{X}$ in time $s$. Since $\eta_0 \in \sg$, the maximum degree is at most $C_* \logn$. Hence, by Poisson thinning, $N \preccurlyeq L \sim \Po(C_* s \logn / n)$. Let $\alpha$ and $\beta$ be the total number of open and closed edges, respectively, that are adjacent to the path of the (static) walk by time $s$. When $N \le 2 C_* s \logn / n$, we have $\alpha \le 2 C_*^2 s \logn[2]/n$ and $\beta \le 2 C_* s \logn$. Hence
\begin{align*}
	\prb{ \se_{0,s}^c }
&\le
	\prb{ \me\rbb{ \alpha(1-p)\mu + \beta p \mu } \le s, \, N \le 2 C_* s \logn/n }
+
	\prb{ L > 2 C_* s \logn/n }
\\&\le
	C s^2 \logn[2] \mu / n
+
	\exp{-c s \logn / n},
\end{align*}
for positive constants $c$ and $C$, by Poisson concentration. Hence, since $\mu n \ll \logn[-14]$, we have
\[
	\prb{ \se_{0,s} }
=
	1 - \oh{1}
\Quad{when}
	s = n \logn[6].
\label{eq:sup:static_coupling_prob;gg}
\nt
\]

Note that $X$ can only become isolated when it is at a degree~1 vertex immediately prior. We use our static-dynamic coupling to lower bound:
\GAP{1}
\[
\begin{aligned}
&	\prb{ \tauisol \le s }
\ge
	\prB{ \textstyle \bigcup_{i=1}^k \brb{ \tauisol \in \bigl[ \tilde{T}_i, \tilde{T}_i' \bigr), \, \tilde{T}_i \le s - n, \, \tauisol - \tilde{T}_i \le n, \, \cap_{j \le i} \se_j } }
\\&\gap \textstyle
=
	\sum_{i=1}^k
	\prb{ \tauisol \in \bigl[ \tilde{T}_i, \tilde{T}_i' \bigr), \, \tauisol - \tilde{T}_i \le n \midb \tilde{T}_i \le s - n, \, \cap_{j \le i} \se_j }
\cdot
	\prb{ \tilde{T}_i \le s - n, \, \cap_{j \le i} \se_j }.
\end{aligned}
\label{eq:sup:disjoint-union}
\nt
\]

Using the static-dynamic coupling on the event $\se_1 \caps \se_i$, we see that if the unique open edge adjacent to $X$ at time $\tilde{T}_i$ closes before anything else opens or $X$ jumps, then $X$ becomes isolated during $[\tilde{T}_i, \tilde{T}_i')$. Writing $\me_1$, $\me_2$ and $\me_3$ for independent exponential random variables, we have
\GAP{3}
\begin{align*}
&	\prb{ \tauisol \in \bigl[ \tilde{T}_i, \tilde{T}_i' \bigr), \, \tauisol - \tilde{T}_i \le n \midb \tilde{T}_i \le s - n, \, \cap_{j \le j} \se_j }
\\&\gap
\ge
	\prb{ \me_1((1-p)\mu) < \min\bra{\me_2(p(n-1)\mu), \me_3(1/n)}, \, \me_1((1-p)\mu) \le n }
\asymp
	\mu n,
\end{align*}
since $\mu n \ll 1$, by comparing rates.
Using this in \eqref{eq:sup:disjoint-union} along with \eqref{eq:sup:Tk<s;gg} and \eqref{eq:sup:static_coupling_prob;gg} we obtain
\begin{align*}
	\prb{ \tauisol \le s }
&
\gtrsim
	\mu n k \,
	\prb{ \tilde{T}_k \le s - n, \, \se_1 \caps \se_k }
\\&
\ge
	\mu n k \rbb{1 - \prb{ \se_{0,s}^c } - \prb{ \tilde{T}_k > s-n } }
\asymp
	\mu n k.
\end{align*}
Since $k = s/(8n)$, this concludes the proof.
	%
\end{Proof}

We now use this to prove our isolation result \Cref{res:sup:tauX<t}.


\begin{Proof}[Proof of \Cref{res:sup:tauX<t}]
In this proof, we use the following shorthand:
\[
	\prb[\mg]{ \cdot }
=
	\min_{\eta_0 \in \sg, x_0 \in \mg_0}
	\prb[x_0,\eta_0]{ \cdot }
\Quad{and}
	\prgb{ \cdot }
=
	\max_{\eta_0 \in \sg, x_0 \in \mg_0}
	\prb[x_0,\eta_0]{ \cdot }.
\]
Consider an initial pair $(x_0,\eta_0)$ with $\eta_0 \in \sh$ and $x_0 \in \mg_0$. For any~$s \in \br$ and $r \in \bn$, using the Markov property we have
\GAP{0.5}
\begin{align*}
&	\prb[x_0,\eta_0]{ \tauisol \in \bigl(sr,s(r+1)\bigr], \, \Ggood{0,sr} }
\\&\gap
\ge
	\prb[x_0,\eta_0]{ \tauisol \le s(r+1) \midb \tauisol > sr, \, \XinG{sr}, \, \Ggood{0,sr} }
\,	\prb[x_0,\eta_0]{ \tauisol > sr, \, \XinG{sr}, \, \Ggood{0,sr} }
\\&\gap
\ge
	\prb[\mg]{ \tauisol \le s }
\cdot
	\rbB{
		\prb[x_0,\eta_0]{ \tauisol > sr, \, \Ggood{0,sr} }
	-	\prb[x_0,\eta_0]{ \existsXnotinG{sr} }
	}.
\end{align*}
Hence we have
\begin{align*}
&	\prb[x_0,\eta_0]{ \tauisol > s(r+1), \, \Ggood{0,sr} }
\\&\gap
=
	\prb[x_0,\eta_0]{ \Ggood{0,sr} }
-	\prb[x_0,\eta_0]{ \tauisol \le sr, \, \Ggood{0,sr} }
-	\prb[x_0,\eta_0]{ \tauisol \in \bigl(sr, s(r+1)\bigr], \, \Ggood{0,sr} }
\\&\gap
\le
	\prb[x_0,\eta_0]{ \Ggood{0,sr} }
-	\prb[x_0,\eta_0]{ \tauisol \le sr, \, \Ggood{0,sr} }
\\&\gap
\qquad -
	\prb[\mg]{ \tauisol \le s }
\cdot
	\rbB{
		\prb[x_0,\eta_0]{ \tauisol > sr, \, \Ggood{0,sr} }
	-	\prb[x_0,\eta_0]{ \existsXnotinG{sr} } }
\\&\gap
=
	\prb[x_0,\eta_0]{ \tauisol > sr, \, \Ggood{0,sr} }
\cdot
	\prgb{ \tauisol > s }
\\&\gap
\qquad +
	\prb[\mg]{ \tauisol \le s }
\cdot
	\prb[x_0,\eta_0]{ \existsXnotinG{sr} }.
\end{align*}
Hence, upon iterating, we obtain
\begin{align*}
	\bp_{x_0,\eta_0}\bigl( \tauisol &> sr, \, \Ggood{0,sr} \bigr)
\le
	\prb[x_0,\eta_0]{ \tauisol > sr, \, \Ggood{0,s(r-1)} }
\\&
\le
	\prgb{ \tauisol > s }^{\!r}
+
	r \cdot
	\prb[\mg]{ \tauisol \le s }
\,	\prb[x_0,\eta_0]{ \existsXnotinG{sr} }.
\label{eq:sup:tauX>sr;gg}
\nt
\end{align*}
Observe that, by the memoryless property, we have
\[
	\prb[\mg]{\tauisol \le s } \le \mu s.
\label{eq:sup:tauX<s|upper}
\nt
\]
We now set $s = n\logn[6]$, $t = \gamma (\log_{(M)}\! n)^{-1}/\mu$ for a constant $\gamma$, to be chosen later, and $r = \floor{t/s}$; note then that $\tfrac23 t \le rs \le t$ as $\mu n \ll \logn[-5]$. Since $\eta_0 \in \sh$, we may apply \Cref{res:leavegiant} for this $t$ to obtain a constant $C$ so that
\[
	\prb[x_0,\eta_0]{ \existsXnotinG{sr} }
\le
	C \mu t \log_{(M)}\! n.
\label{eq:sup:giant_stay-in}
\nt
\]
Using \eqref{eq:sup:tauX<s|upper} and \eqref{eq:sup:giant_stay-in} along with \Cref{res:sup:tauX<s} in \eqref{eq:sup:tauX>sr;gg}, we find that
\[
	\prb[x_0,\eta_0]{ \tauisol > t, \, \Ggood{0,t} }
\le
	\rbr{1 - q \mu s}^r + C (\mu t)^2 \log_{(M)}\! n
\le
	1 - \tfrac13 q \mu t + C (\mu t)^2 \log_{(M)}\! n,	
\]
valid for any~$(x_0,\eta_0)$ with $\eta_0 \in \sh$ and $x_0 \in \mg_0$.
We then take $\gamma = q/(6C)$ and obtain
\[
	\prb[x_0,\eta_0]{ \tauisol > t, \, \Ggood{0,t} }
\le
	1 - \frac{q^2}{36C} \cdot \frac1{\log_{(M)}\! n}.
\]
Since we can take $C \ge 1$ and $q \le 1$, we then have
\GAP{3}
\begin{align*}
&	\prbb[x_0,\eta_0]{ \tauisol > \frac1{C \mu \log_{(M)}\! n}, \, \Ggoodbb{0,\frac1{C \mu \log_{(M)}\! n}} }
\\&\gap
\le
	\prbb[x_0,\eta_0]{ \tauisol > \frac{q}{3 C \mu \log_{(M)}\! n}, \, \Ggoodbb{0,\frac{q}{3 C \mu \log_{(M)}\! n}} }.
\end{align*}
Hence there exists a positive constant $c$ so that
\[
	\prbb[x_0,\eta_0]{ \tauisol > \frac1{\mu \log_{(M)}\! n}, \, \Ggoodbb{0,\frac1{\mu \log_{(M)}\! n}} }
\le
	1 - c \cdot \frac1{\log_{(M)}\! n}.
\]
Finally, $\eta_0 \in \sh$, so $\pr[\eta_0]{ \Gbad{0,(\log_{(M)}\! n)^{-1}/\mu} } \le n^{-1}$, and the result follows.
\end{Proof}

\begin{rmkt*}
	Observe that, as in the remark after \Cref{res:sup:tauX<s}, for this time-scale the result of \Cref{res:sup:tauX<t} is best-possible, up to constants.
\end{rmkt*}

\subsection{Dual-Walker (Joint) Isolation Time}
\label{sec:isol2}

In this section we prove \Cref{res:2:tau>t;H} on the joint isolation time of two walkers on a single dynamical environment. We start by introducing some more notation. Consider two walks $X$~and~$Y$, which start from $x_0$ and $y_0$ respectively, walking \emph{independently} on the same environment $\eta$.
Let
\[
	\tau = \inf\brb{ t \ge 0 \mid d_t(X_t) = 0 = d_t(Y_t) }.
\]
Let $\tau_0^X = \tau_0^Y = \hat\tau_0^X = 0$, and for $k \ge 1$ define inductively
\begin{gather*}
	\tau_k^X 		= \inf\brb{ t \ge \hat\tau_{k-1}^X	\mid d_t(X_t) = 0 },
\quad
	\hat\tau_k^X	= \inf\brb{ t \ge \tau_k^X \mid d_t(X_t) \neq 0 }
\\
\text{and} \quad
	\tau_k^Y		= \inf\brb{ t \ge \tau_k^X \mid d_t(Y_t) = 0 }.
\end{gather*}

We prove a result on the joint-isolation time of two walks, $X$ and $Y$, walking independently on the same (dynamic) environment~$\eta$. For the probability measure associated to this system~$(X,Y,\eta)$, when it is started from $(x_0,y_0,\eta_0)$, we write $\bp_{x_0,y_0,\eta_0}$.


In order to prove the dual-walker isolation result, we first state two lemmas that we use. We prove the theorem using the lemmas, then prove the lemmas. Throughout, $M$ is a positive integer.

\begin{lem} \label{res:2:term1}
	There exists a positive constant $c_1$ so that, for all~$M$, all~$n$ sufficiently large and all~$(x_0,y_0,\eta_0)$, we have
	\[
		\prb[x_0,y_0,\eta_0]{
			\tau_1^Y > \hat\tau_1^X, \, \Hgoodb{\tau_1^X, \hat\tau_1^X} }
	\le
		\expb{- c_1 / \log_{(M)}\! n }.
	\]
\end{lem}

\begin{lem} \label{res:2:term2}
	There exists a positive constant $c_2$ so that, for all~$M$, all~$n$ sufficiently large and all~$(x_0, \eta_0)$, we have
	\[
		\prb[x_0,\eta_0]{ \tau_{K+1}^X > t, \, \Hgood{0,t} }
	\le
		\expb{-\tfrac23 K}
	\Quad{when}
		K = \floorb{ c_2 \mu t / \log_{(M)}\! n }.
	\]
\end{lem}

\begin{Proof}[Proof of \Cref{res:2:tau>t;H}]
By monotonicity, replacing $M$ by $M-1$, it suffices to find a positive constant $c$ so that the probability is upper bounded by $2 \exp{ - c \mu t / (\log_{(M)}\! n)^2}$ for a positive constant $c$. Hence we may assume that $\mu t \ge (\log_{(M)}\! n)^2$, as otherwise the result trivially holds.

For any~$t \ge 0$ and for $K = \floor{ c_2 \mu t / \log_{(M)}\! n }$, using \Cref{res:2:term2} we have
\[
	\prb[x_0,y_0,\eta_0]{ \tau > t, \, \Hgood{0,t} }
\le
	\prb[x_0,y_0,\eta_0]{ \tau > t, \, \tau_{K+1}^X \le t, \, \Hgood{0,t} }
+	\expb{-\tfrac23 K}.
\label{eq:2:tau0>t_sum}
\nt
\]
Since $\tau_{K+1}^X \le t$ implies $\hat\tau_K^X \le t$, on the event $\{\tau_{K+1}^X \le t\}$ we have $\Hgood{0,t} \subseteq \Hgood{0, \hat\tau_K^X}$. We then use the strong Markov property at time $\hat\tau_1^X$ to iterate:
\GAP{0.5}
\begin{align*}
&	\prb[x_0,y_0,\eta_0]{ \tau > t, \, \tau_{K+1}^X \le t, \, \Hgood{0,t} }
\le
	\prb[x_0,y_0,\eta_0]{
		\cap_{k=1}^{K}
		\brb{ \tau_k^Y > \hat\tau_k^X },
	\,
		\Hgoodb{0, \hat\tau_K^X} }
\\&\gap
\le
	\prb[x_0,y_0,\eta_0]{
		\cap_{k=2}^{K}
		\brb{ \tau_k^Y > \hat\tau_k^X },
	\,
		\Hgoodb{\hat\tau_1^X, \hat\tau_K^X}
	\midb
		\tau_1^Y > \hat\tau_1^X,
	\,
		\Hgoodb{0, \hat\tau_1^X} }
\\&\gap\qquad
\cdot
	\prb[x_0,y_0,\eta_0]{
		\tau_1^Y > \hat\tau_1^X,
	\,
		\Hgoodb{\tau_1^X, \hat\tau_1^X} }
\\&\gap
\le
	\max_{x_0',y_0',\eta_0'}
	\prb[x_0',y_0',\eta_0']{
		\cap_{k=1}^{K-1}
		\brb{ \tau_k^Y > \hat\tau_k^X },
	\,
		\Hgoodb{0, \hat\tau_{K-1}^X} }
\cdot
	\prb[x_0,y_0,\eta_0]{
		\tau_1^Y > \hat\tau_1^X,
	\,
		\Hgoodb{\tau_1^X, \hat\tau_1^X} }
\\&\gap
\le \cdots \le
	\max_{x_0',y_0',\eta_0'}
	\prb[x_0',y_0',\eta_0']{ \tau_1^Y > \hat\tau_1^X, \, \Hgoodb{\tau_1^X, \hat\tau_1^X} }^K.
\label{eq:2:iteration}
\nt
\end{align*}
Since $\mu t \ge (\log_{(M)}\! n)^2$, we have $K \ge \tfrac12 c_2 \mu t / \log_{(M)}\! n$. Using \eqref{eq:2:iteration} and \Cref{res:2:term1} in \eqref{eq:2:tau0>t_sum}, we have
\GAP{5}
\begin{align*}
&	\prb[x_0,y_0,\eta_0]{ \tau > t, \, \Hgood{0,t} }
\le
	\max_{x_0',y_0',\eta_0'}
	\prb[x_0',y_0',\eta_0']{ \tau_1^Y > \hat\tau_1^X, \, \Hgoodb{\tau_1^X, \hat\tau_1^X} }^K
+	\expb{-\tfrac23 K}
\\&\gap
\le
	\expb{ - c_1 K / \log_{(M)}\! n }
+	\expb{ - \tfrac23 K }
\le
	2 \expb{ - \tfrac12 c_1 c_2 \mu t / (\log_{(M)}\! n)^2 }.
\qedhere
\end{align*}
\end{Proof}

It remains to prove \Cref{res:2:term1,res:2:term2}.

\begin{Proof}[Proof of \Cref{res:2:term1}]
By the strong Markov property, used in the same way as above, and recalling that $\tau_1^Y$ is the first time after $\tau_1^X$ that $Y$ becomes isolated, we have
\[
	\max_{x_0,y_0,\eta_0}
	\prb[x_0,y_0,\eta_0]{
		\tau_1^Y > \hat\tau_1^X, \, \Hgoodb{\tau_1^X, \hat\tau_1^X} }
\le
	\max_{\substack{x_0,y_0,\eta_0 \\ d_0(x_0) = 0}}
	\prb[x_0,y_0,\eta_0]{
		\tau_1^Y > \hat\tau_1^X, \, \Hgoodb{0, \hat\tau_1^X} }.
\label{eq:2:Yisol|Xisol}
\nt
\]

For the moment, we emphasise that our underlying graph has $n$ vertices: we do this by using super- and subscript $n$, eg $\bp^n$ and $\sg_n$. Recall \Cref{res:1:tauX>t;H}, which says that
\[
	\bp^n_{y_0,\eta_0}\rbb{ \tauisol^Y > t, \, \Hgood{0,t} }
\le
	\expb{ - \mu t / \log_{(M)}\! n }
\Quad{when}
	\mu t \ge 3.
\]
We wish to bound (the related quantity)
\[
	\bp^n_{x_0,y_0,\eta_0} \rbb{
		\tau_1^Y > \hat\tau_1^X, \, \Hgoodb{0, \hat\tau_1^X} \midb \hat\tau_1^X }.
\]
This is trivially 0 for $x_0 = y_0$; consider $x_0 \neq y_0$. To bound this, we observe that, conditional on the value of $\hat\tau_1^X = T$, this is conditioning the vertex $x_0$ to be isolated until time $\hat\tau_1^X = T$; the rest of the graph is unaffected.

Let $\eta^n$ be a dynamical environment on $n$ vertices, and let $x_0 \in \{1,...,n\}$ be a vertex. Define $\tilde{\eta}^n$ by conditioning on the event that the vertex $x_0$ is isolated until time $T$. Write $\hat\eta^n$ for the restriction of $\tilde{\eta}^n$ to $\{1,...,n\} \setminus \{x_0\}$. Observe then that $\hat\eta^n \sim \eta^{n-1}$ (where $\eta^{n-1}$ is a dynamical environment on $n-1$ vertices), up to relabelling of vertices. (In words, this says that if a vertex is conditioned to be isolated, then the rest of the graph behaves as a dynamical environment on $n-1$ vertices.) Note also that two edges do not update at the same time, so we cannot have $\tau_1^Y = \hat\tau_1^X$ (since the first requires an edge to close and the second an edge to open).

Hence, $(Y_t \mid t \le \hat\tau_1^X)$ is a walk on the environment $\hat\eta^n$, which has the distribution of $\eta^{n-1}$. Note that $Y$ may still pick the (conditioned to be isolated) vertex $x_0$ (with probability $1/(n-1)$), in which case it does not move; thus, under this conditioning, $(Y,\hat\eta^n)$ is simply a realisation of dynamical percolation on $n-1$ vertices, but with added laziness: when $Y$'s $\me(1)$ clock rings, with probability $1/(n-1)$ it does nothing; with the remaining probability, it performs the usual step.

Note that we can rescale $\mu$ to get rid of the laziness of $Y$. Indeed, the laziness has the effect of changing the walker's clock from rate 1 to rate $1-1/(n-1)$. As such, if we replace $\mu$ by $\mu' = \mu (1-1/(n-1))$, then the ratio of the rate edge-clocks to the rate of the walker-clock is $\mu$: we have simply slowed both down. We then speed up everything by a factor $1-1/(n-1)$. We apply previous results with $\mu$ replaced by $\mu'$. The restrictions on $\mu$ are satisfied by $\mu'$ also, since $\mu' \le \mu$.

For all~$m$, define $\sg_m'$ by replacing $c_*$ and $C_*$ in \Cref{defn:good-graph} by $\tfrac12 c_*$ and $2 C_*$, respectively; define $\sh'$ in terms of $\sg'$ as in \Cref{defn:very-good-graph}. We then have that if $\eta_0^n \in \sg_n$ and the vertex $x_0$ is isolated (in $\eta_0^n$), then $\tilde\eta_0^n = \eta_0^n - \{x_0\}$ defined by removing the vertex $x_0$ satisfies $\tilde\eta_0^n \in \sg_{n-1}'$ (for $n$ sufficiently large).
Hence we have the following inequality:
	let $(Z,\zeta)$ be a full system, independent of $X$ and $Y$, on $n-1$ vertices, and start it from $(Z_0,\zeta_0) = (y_0,\hat\eta_0^n)$; we then have
	\[
		\bp^n_{x_0,y_0,\eta_0^n} \rbb{
			\tau_1^Y > \hat\tau_1^X, \, \sh_n\sbb{0, \hat\tau_1^X} \midb \hat\tau_1^X }
	\le
		\bp^{n-1}_{y_0,\hat\eta_0^n} \rbb{
			\tau_1^Z > \hat\tau_1^X, \, \sh_{n-1}'\sbb{0, \hat\tau_1^X} \midb \hat\tau_1^X }.
	\]
Note that \Cref{res:1:tauX>t;H} still holds if we replace $\sh$ by $\sh'$ in its statement.
Combining all the above considerations, applying \Cref{res:1:tauX>t;H}, on the event $\{\mu \hat\tau_1^X \ge 3\}$ we have
\GAP{3}
\begin{align*}
&	\bp^n_{x_0,y_0,\eta_0^n} \rbb{
		\tau_1^Y > \hat\tau_1^X, \, \sh_n\sbb{0, \hat\tau_1^X} \midb \hat\tau_1^X }
\le
	\bp^{n-1}_{y_0,\hat\eta_0^n} \rbb{
		\tau_1^Z > \hat\tau_1^X, \, \sh_{n-1}'\sbb{0, \hat\tau_1^X} \midb \hat\tau_1^X }
\\&\gap
\le
	\expb{ - \mu \hat\tau_1^X \rbb{1-\tfrac1n} / \log_{(M)}\! (n-1) }
\le
	\expb{ - \tfrac12 \mu \hat\tau_1^X / \log_{(M)}\! n }.
\end{align*}
We now calculate the unconditioned value. Fix $(x_0,y_0,\eta_0)$ with $d_0(x_0) = 0$. We have
\GAP{2}
\begin{align*}
&	\prb[x_0,y_0,\eta_0]{
		\tau_1^Y > \hat\tau_1^X, \, \sh\sbb{0, \hat\tau_1^X} \midb \hat\tau_1^X }
\\&\gap
\le
	\exb[x_0,y_0,\eta_0]{
		\prb[x_0,y_0,\eta_0]{
			\tau_1^Y > \hat\tau_1^X, \, \sh\sbb{0, \hat\tau_1^X} \midb \hat\tau_1^X }
	\cdot
		\oneb{ \hat\tau_1^X \ge 3/\mu } }
+	\prb[x_0,\eta_0]{ \hat\tau_1^X < 3/\mu }
\\&\gap
\le
	\prb[x_0,\eta_0]{ \hat\tau_1^X \ge 3/\mu }
\cdot
	\expb{ - \tfrac12 \mu \rbr{3/\mu} / \log_{(M)}\! n }
+	\prb[x_0,\eta_0]{ \hat\tau_1^X < 3/\mu }
\\&\gap
\le
	1 - \prb[x_0,\eta_0]{ \hat\tau_1^X < 3/\mu } / \log_{(M)}\! n
\le
	\expb{ - \prb[x_0,\eta_0]{ \hat\tau_1^X < 3/\mu } / \log_{(M)}\! n }
\end{align*}
where we have used the inequality $\exp{-\tfrac32 x} \le 1 - x$, valid for sufficiently small $x$.
Since $d_0(x_0) = 0$, we have $\hat\tau_1^X \sim \me( \lambda \mu (1-1/n) )$,
and hence we have
\(
	\pr[x_0,\eta_0]{ \hat\tau_1^X \ge 3/\mu } \asymp 1.
\)
Substituting this into \eqref{eq:2:Yisol|Xisol} gives the required bound.
\end{Proof}

\begin{Proof}[Proof of \Cref{res:2:term2}]
We may assume that $K \ge 1$, otherwise the result is trivial.

For this lemma we only consider one walker, $X$; as such, we drop the $X$ superscripts. We define
\[
	\tauisol(s) = \inf\brb{ t \ge s \mid d_t(X_t) = 0 }
\Quad{and}
	\hat\tau_\isol(s) = \inf\brb{ t \ge \tauisol(s) \mid d_t(X_t) > 0 };
\]
also write $\tauisol = \tauisol(0)$ and $\hat\tau_\isol = \hat\tau_\isol(0)$.

For $k = 0,...,3K$, set $t_k = \tfrac{t}{3K}$ and $t_k' = t_k + \tfrac12 t / (3K)$; also, for $k = 1,...,3K$, write
\[
	\sh_k = \sh\sbb{t_{k-1},t_{k-1}'}
\Quad{and}
	\mj_k
=
	\brb{ \hat\tau_\isol(t_{k-1}) \le t_k }.
\]
If $\mj_k$ occurs then at some point in the interval $[t_{k-1},t_k]$ the walk is isolated and at a later point (in the same interval) is not. Observe that we have
\begin{align*}
	\brb{ \tau_{K+1} > t } \cap \Hgood{0,t}
&\subseteq
	\brB{ {\textstyle \sum_{k=1}^{3K}} \oneb{\mj_k} \le K } \cap \Hgood{0,t} 
\\
&=
	\brB{ {\textstyle \sum_{k=1}^{3K}} \oneb{\mj_k^c} \ge 2K } \cap \Hgood{0,t}
\subseteq
	\brB{ {\textstyle \sum_{k=1}^{3K}} \oneb{\mj_k^c \cap \sh_k} \ge 2K }.
\end{align*}
Write $J = \sum_{k=1}^{3K} \oneb{\mj_k^c \cap \sh_k}$. Note that by the Markov property we have $J \preccurlyeq \Bin(3K,q)$ where
\[
	q
=
	\max_{x_0,\eta_0}
	\prb[x_0,\eta_0]{ \mj_1^c \cap \sh_1 }.
\]
We shall show, for a suitable constant $c_2$ in the definition of $K$, that $q \le \tfrac13$, and then deduce that
\[
	\prb[x_0,\eta_0]{ \tau_{K+1} > t, \, \Hgood{0,t} }
\le
	\prb{ \Bin(3K,\tfrac13) \ge 2K }
\le
	\expb{- \tfrac23 K}.
\]

Observe that we have
\[
	\brB{ \hat\tau_\isol > \tfrac{t}{3K}, \, \tauisol \le \tfrac{t}{6K} }
\subseteq
	\brB{ \hat\tau_\isol - \tauisol > \tfrac{t}{6K} }.
\]
Thus we have, for any~$(x_0,\eta_0)$, that
\[
	\prb[x_0,\eta_0]{ \mj_1^c \cap \sh_1 }
\le
	\prb[x_0,\eta_0]{ \hat\tau_\isol - \tauisol > \tfrac{t}{6K} }
+	\prb[x_0,\eta_0]{ \tauisol > \tfrac{t}{6K}, \, \sh_1 }.
\]
The first term is simply
\[
	\prb{ \me\rbb{(n-1)\lambda\mu/n} > \tfrac{t}{6K} }
\le
	\prb{ \me(\mu) > \tfrac{t}{6K} }
=
	\expb{ - \tfrac16 \mu t / K }
\]
since there are $n-1$ edges that can open, and $\lambda > 1$. Applying \Cref{res:1:tauX>t;H}, we have
\[
	\prb[x_0,\eta_0]{ \tauisol > \tfrac{t}{6K}, \, \sh_1 }
\le
	2 \expb{ - \tfrac16 \mu t / (K \log_{(M)}\! n) }.
\]
Combining these two bounds, we then find that
\begin{align*}
	q
=
	\max_{x_0,\eta_0}
	\prb[x_0,\eta_0]{ \mj_1^c \cap \sh_1 }
&\le
	\expb{ - \tfrac16 \mu t / K }
+	2 \expb{ - \tfrac16 c^*_1 \mu t / (K \log_{(M)}\! n) }
\\&
\le
	3 \expb{ - \tfrac16 c^*_1 \mu t / (K \log_{(M)}\! n) }.
\end{align*}
Hence there exists a positive constant $c_2$ so that if $K = \floor{c_2 \mu t / \log_{(M)}\! n}$ then $q \le \tfrac13$.
\end{Proof}

\counterwithout{thm}{subsection}
\counterwithin{thm}{section}
\section{Coupling}
\label{sec:coup}

\subsection{Statement and Application of Coupling to Mixing}

For this section only, we call a graph \textit{good} if it satisfies the conditions of \Cref{defn:good-graph} and in addition the condition that at least a proportion $c_*$ of its vertices are isolated. Since this is an additional condition, the probability that a graph is good decreases, and hence all our isolation results (from \S\ref{sec:isol-times}) still hold with this extra condition. Recall also the definition of $H$ from \Cref{defn:very-good-graph}, and in particular that $\pi_\ER(H) = 1 - \oh1$.

%

In this section, $(X,\eta)$ and $(Y,\xi)$ are two realisations of the dynamical percolation system; we shall define a Markovian coupling of the two systems, and find a tail bound on the coupling time.
We look first at the case when the environments $\eta$ and $\xi$ start with $\eta_0 = \xi_0$.

\begin{prop}[Coupling Tail Bound]
\label{res:coup:same-env}
	There exists a Markovian coupling, which we denote by
	$\bp_{(x_0,\eta_0),(y_0,\eta_0)}$ when $(X,\eta)$ and $(Y,\xi)$ start from $(x_0,\eta_0)$ and $(y_0,\xi_0)$ respectively, so that, for all~$M \in \bn$, all~$n$ sufficiently large, all~$t$ and all~$(x_0,y_0,\eta_0)$, we have
	\[
		\prb[(x_0,\eta_0),(y_0,\eta_0)]{ (X_t,\eta_t) \neq (Y_t,\xi_t) }
	\le
		3 \expb{- \mu t / \log_{(M)}\! n}.
	\]
\end{prop}


From \Cref{res:coup:same-env} we are able to deduce the upper bounds in \Cref{res:main_thm:mixing_full,res:main_thm:mixing_rw}.

\begin{Proof}[Proof of \Cref{res:main_thm:mixing_rw}]
Observe that we have
\[
	\tvb{ \prb[x_0,\eta_0]{ X_t = \cdot } - \prb[y_0,\eta_0]{ Y_t = \cdot } }
\le
	\tvb{ \prb[x_0,\eta_0]{ X_t = \cdot, \eta_t = \cdot } - \prb[y_0,\eta_0]{ Y_t = \cdot, \xi_t = \cdot } }.
\]
While the walk component alone is not a Markov chain, the full system is. It is then standard to upper bound the total variation distance by the tail probability of the coalescence time:
\[
	\tvb{ \prb[x_0,\eta_0]{ X_t = \cdot, \eta_t = \cdot } - \prb[y_0,\eta_0]{ Y_t = \cdot, \xi_t = \cdot } }
\le
	\prb[(x_0,\eta_0),(y_0,\eta_0)]{ \tauc > t }.
\]
The coupling is coalescent, and so this tail bound is monotone in $t$. \Cref{res:coup:same-env} now implies that this is $\oh1$ when $t = \tfrac1\mu \log_{(M)}\! n$, since $\eta_0 \in H$ (and replacing $M$ with $M+1$ in the proposition).

Since the uniform distribution $\pi_\RW$ is invariant for our walk on \emph{any} graph,
\(
	\pr[\pi_\RW, \eta_0]{ Y_t = \cdot } = \pi_\RW
\)
for any $\eta_0$. Thus we obtain our result: for $t = \tfrac1\mu \log_{(M)}\! n$ and any $\eta_0 \in H$, we have
\begin{align*}
	\max_{x_0} \,
	\tvb{ \prb[x_0,\eta_0]{ X_t = \cdot } - \pi_\RW }
&
\le
	\max_{x_0,y_0} \,
	\tvb{ \prb[x_0,\eta_0]{ X_t = \cdot } - \prb[y_0,\eta_0]{ Y_t = \cdot } }
\\&
\le
	\max_{x_0,y_0} \,
	\prb[(x_0,\eta_0), (y_0,\eta_0)]{ \tauc > t }
=
	\oh1.
\qedhere
\end{align*}
\end{Proof}

In order to prove the mixing of the full system $(X,\eta)$, we also need to know the mixing of the environment by itself.
First recall that the environment process is simply a $p$-biased walk on the hypercube $\{0,1\}^N$, where $N = \binom n2$ and $p = \lambda/n$, and where each coordinate refreshes at rate $\mu$.
We state the result now, then prove it at the end of the subsection.

\begin{prop}[Hypercube Mixing]
\label{res:hyp}
	Consider the rate-1 $p$-biased random walk on the hypercube $\{0,1\}^N$, with $1/N \ll p \le \tfrac12$; denote it $\eta = (\eta_t)_{t\ge0}$, and its invariant distribution $\pi_p$.
	There is cutoff at $\tfrac12 \log(N/p)$ with window order 1:
	for all $\eps \in (0,1)$, there exists a constant $C_\eps$ so that
	\begin{alignat*}{3}
		\textstyle \max_{\eta_0} \,
		\tvb{ \prb[\eta_0]{ \eta_t = \cdot } - \pi_p }
	&\le
		\eps
	\quad&\text{if}\quad
		t &\ge \tfrac12 \log(N/p) + C_\eps
	\\
		\textstyle \min_{\eta_0} \,
		\tvb{ \prb[\eta_0]{ \eta_t = \cdot } - \pi_p }
	&\ge
		1-\eps
	\quad&\text{if}\quad
		t &\le \tfrac12 \log(N/p) - C_\eps.
	\end{alignat*}
\end{prop}

Since we work in continuous time, we can apply this directly when the refresh-rate is $\mu$.

\begin{Proof}[Proof of \Cref{res:main_thm:mixing_full}]
We consider first a lower bound on $\tmix(\eps)$.
Observe that, trivially,
\[
	\tvb{ \prb[x_0,\eta_0]{ (X_t,\eta_t) \in \cdot } - \pi_U \times \pi_\ER }
\ge
	\tvb{ \prb[\eta_0]{ \eta_t \in \cdot } - \pi_\ER }.
\]
Thus it suffices to only show that the environment has not mixed by time $t$.
This follows immediately from the lower bound in \cref{res:hyp}, since $\tfrac12 \log(N/p) = \tfrac32 \logn + \Th1$.

Now consider the upper bound.
Fix $(x_0,\eta_0)$ and $(y_0,\xi_0)$. By Chapman-Kolmogorov, we have
\[
	\prb[x_0,\eta_0]{ (X_{s+t},\eta_{s+t}) \in \cdot }
=
	\exb[x_0,\eta_0]{ \prb[X_s,\eta_s]{ (X_t,\eta_t) \in \cdot } }.
\]
Hence for any coupling $\bq$ of $\pr[x_0,\eta_0]{ (X_s,\eta_s) = \cdot }$ and $\pr[y_0,\xi_0]{ (Y_s,\xi_s) = \cdot }$ we have
\GAP{2}
\begin{align*}
&	\tvb{
		\prb[x_0,\eta_0]{ (X_{s+t},\eta_{s+t}) \in \cdot }
	-	\prb[y_0,\xi_0]{  (Y_{s+t},\xi_{s+t})  \in \cdot } }
\\&\gap
\le
	\bq\rbb{ \eta_s \neq \xi_s }
+	\prb{ \eta_s \notin H }
+	\max_{x_0',y_0',\eta_0' \in H}
	\tvb{
		\prb[x_0',\eta_0']{ (X_t,\eta_t) \in \cdot }
	-	\prb[y_0',\eta_0']{ (Y_t,\xi_t)  \in \cdot } }.
\end{align*}
In particular, consider the following such coupling $\bq$: fix $s \ge 0$, and couple $(\eta_s,\xi_s)$ using the optimal coupling when started from $(\eta_0,\xi_0)$; given~$\eta_s$ (and the fixed $\eta_0$), sample $X_s$ conditional on $\eta_s$ (and $\eta_0$); do similarly (and independently) for $Y_s$ with $\xi_s$ (and $\xi_0$). This then has
\[
	\bq\rbb{ \eta_s \neq \xi_s }
=
	\tvb{ \prb[\eta_0]{ \eta_s \in \cdot } - \prb[\xi_0]{ \xi_s \in \cdot } }.
\]

Now fix $M \in \bn$ and choose $s$ so that $\mu s = \tfrac32 \logn + \log_{(M+2)}\! n$, which has $\mu s \ge \tfrac12 \log(N/p) + \log_{(M+2)}\! n$.
Then
by the upper bound in \cref{res:hyp}
and the triangle inequality, we have
\[
	\bq\rbb{ \eta_s \neq \xi_s } = \oh1
\Quad{and}
	\prb{ \eta_s \notin H } \le \pi_\ER(H^c) + \oh1 = \oh1.
\]

Since our coupling $\bp_{\cdot,\cdot}$ from \Cref{res:coup:same-env} is Markovian and coalescent, we have
\[
	\tvb{
		\prb[x_0,\eta_0]{ (X_t,\eta_t) \in \cdot }
	-	\prb[y_0,\eta_0]{ (Y_t,\xi_t)  \in \cdot } }
\le
	\prb[(x_0,\eta_0),(y_0,\eta_0)]{ \tauc > t }.
\]
Noting the conditions of \Cref{res:coup:same-env}, this implies that
\[
	\prb[(x_0,\eta_0),(y_0,\xi_0)]{ \tauc > t } \le \eps^2
\Quad{when}
	t = \tfrac1\mu \log\rbb{3/\eps^2} \log_{(M+1)}\! n.
\]

Combining these three bounds we obtain, for these $s$ and $t$, that
\[
	\tvb{
		\prb[x_0,\eta_0]{ (X_{s+t},\eta_{s+t}) \in \cdot }
	-	\prb[y_0,\xi_0]{  (Y_{s+t},\xi_{s+t})  \in \cdot } }
\le
	\eps^2 + \oh1 + \oh1
\le
	\eps.
\]
Hence for all~$\eps \in (0,1)$ we have
\[
	\mu \cdot \tmix(\eps)
\le
	\tfrac32 \logn + \log_{(M)}\! n.
\]
This completes the proof of the upper bound.
\end{Proof}

It remains to prove \Cref{res:hyp}.

\begin{Proof}[Proof of \Cref{res:hyp}]
We prove the upper bound first.
We do this by relating the TV distance to the $L_\infty$ distance.
The probability an edge is in the same state $z \in \{0,1\}$ as initially is exactly
\[
	e^{-t} + \rbb{ 1-e^{-t} } \, \prb{ \Bern(p) = z }
=
	p^z (1-p)^{1-z} + e^{-t} \rbb{ 1 - p^z (1-p)^{1-z} }.
\]
Also, it is well-known that for a reversible transition kernel $P = (P_t)_{t\ge0}$ with invariant distribution $\pi$, writing $d_p(t)$ for the $p$-norm at time $t$ (for $p \in [1,\infty]$), we have
\[
	d_\TV(t) = \tfrac12 d_1(t) \le \tfrac12 d_2(t)
\Quad{and}
	d_\infty(2t)
=
	\rbb{ d_2(t) }^2
=
	\max_x \, P_{2t}(x,x)/\pi(x) - 1;
\]
see \cite[Exercise 4.5 and Proposition 4.15]{LPW:markov-mixing}.
We hence deduce that
\[
	d_\infty(2t)
=
	\max_{\eta_0} \, \normb{ \prb[\eta_0]{ \eta_{2t} \in \cdot } - \pi_\ER }_\infty
=
	\max_{\eta_0} \, \prb[\eta_0]{ \eta_{2t} = \eta_0 } / \pi_\ER(\eta_0) - 1.
\]
Calculating this directly, recalling that $N = \binom n2$ and $p = \lambda/n$ with $\lambda$ a constant, we see that
\[
	d_\infty(2t)
=
	\rbb{ 1 + e^{-2t}(1/p - 1) }^N - 1
\le
	\expb{ e^{-2t} N/p } - 1.
\]
Hence if we set $t = \tfrac12 \log(N/p) + \tfrac12 C_\eps$, for some large constant $C_\eps$, then we obtain
\[
	d_\infty(2t)
\le
	\expb{ 1/C_\eps } - 1
\le
	2/C_\eps.
\]
Finally we deduce that
\(
	d_\TV(t) \le \tfrac12 \sqrt{ d_\infty(2t) } \le 1/\sqrt{2C_\eps},
\)
proving the upper bound.

We now pursue the lower bound.
For this, we consider the statistic
\(
	N_t = \sum_{e=1}^N \one{ \eta_t(e) = 1 },
\)
ie the number of open edges at time $t$.
Observe that $N_t \sim \Bin(N,p)$ when $\eta_0 \sim \pi_\ER$.
Consider starting $\eta_0$ from the all-1 state, which we denote $\ONE \in \{0,1\}^N$.
Then write
\[
	q_t = e^{-t} + (1 - e^{-t})p = p + e^{-t}(1-p),
\]
and observe that $N_t \sim \Bin(N,q_t)$ when $\eta_0 = \ONE$.
Now define the set
\[ \textstyle
	A_t
=
	\brb{ \zeta \in \{0,1\}^N \mid \sum_{e=1}^N \oneb{ \zeta(e) = 1 } \ge \tfrac12(p+q_t) N }.
\]
This will be our distinguishing statistic/set.
Recall that
\[
	\exb{ \Bin(N,r) } = Nr
\Quad{and}
	\Var{ \Bin(N,r) } = Nr(1-r) \le Nr.
\]
Take $t = \tfrac12 \log(N/p) - \tfrac12 \log C_\eps$, for some large constant $C_\eps$.
Note that $q_t - p \ge \tfrac12 \sqrt{C_\eps p / N}$; also $q_t \le 2p$ since $\sqrt{p/N} \ll p$.
(This is where we use the condition $p \gg 1/N$.)
Hence, by Chebyshev,
\[
	\prb[\ONE]{ \eta_t \notin A }
\le
	\prb{ \absb{ \Bin(N,q_t) - q_t N } \ge \tfrac12(q_t - p) N }
\le
	\frac{ 4Nq_t }{ (q_t-p)^2 N^2 }
\le
	\frac{ 50p }{ (Cp/N) \cdot N }
=
	\frac{50}{C_\eps};
\]
similarly, $\pr[\pi_\ER]{ \eta_t \in A } \le 50/C_\eps$.
Hence $d_\TV(t) \ge 1 - 100/C_\eps$, proving the lower bound.
\end{Proof}

It remains to prove \Cref{res:coup:same-env}. To prove this, we carefully define a coupling, and use the result on dual-walker isolation, \Cref{res:2:tau>t;H}, that we proved in the previous section.

\subsection{Coupling Description and Proof of Tail Bound}
\label{sec:coup:description_proof}

Below, we write $(x,y)$ for the \emph{undirected} edge with endpoints $x$ and $y$; in particular, $(x,y) = (y,x)$.
We only use the coupling below once the environments have coupled and the two walks have subsequently become then jointly isolated.
We now define the coupling.

\begin{figure}
	\centering
	\includegraphics[width=0.9\textwidth]{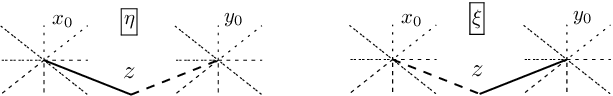}
	
	\vspace{4mm}
	
	\parbox{0.8\textwidth}{\caption{%
			The dotted lines represent the closed edges incident to $x$ and $y$ (recall that they are both isolated initially).
			The full   line indicates opening an edge to another vertex $z$;
			the dashed line indicates leaving it closed.
			The dotting/dashing is reversed in $\xi$ compared with $\eta$%
		}
		\label{fig:coup:opening_edges}
	}
\end{figure}

\begin{figure}
	\centering
	\includegraphics[width=0.9\textwidth]{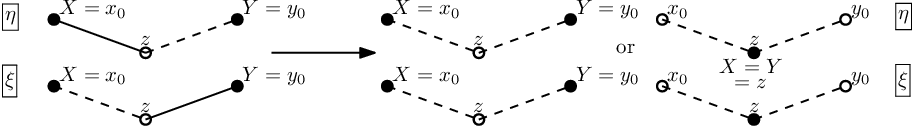}
	
	\vspace{4mm}
	
	\parbox{0.8\textwidth}{\caption{%
			The full   line indicates an open edge;
			the dashed line indicates a closed edge.
			The walkers move along the open edges, moving together.
			On the left-hand  side the filled dots represent where the walkers start;
			on the right-hand side the filled dots represent where the walkers end;
			the empty circles indicate empty sites%
		}
		\label{fig:coup:walker_moves}
	}
\end{figure}

\begin{defnt}
\label{defn:coupling}
Suppose that $(X,\eta)$ and $(Y,\xi)$ are in the states $(x,\eta_0)$ and $(y,\xi_0)$, respectively.
Assume that $\eta_0 = \xi_0$ and both $x$ and $y$ are isolated vertices in the environment $\eta_0 = \xi_0$.
Let $\eta$ evolve in the standard way.
Couple $\xi$ to $\eta$ as follows.
Suppose that edge $(u,v)$ refreshes in $\eta$:
\begin{itemize}[noitemsep,topsep=0pt,label={\ensuremath{\bm\cdot}}]
	\item 
	if $\{u,v\} \cap \{x,y\} = \emptyset$ or $\{u,v\} = \{x,y\}$, then perform the same update to $(u,v)$ in $\xi$ as in $\eta$;
	
	\item 
	if $u = x$ and $v \notin \{x,y\}$, then perform the same update to $(y,v)$ in $\xi$ as to $(u,v)$ in $\eta$;
	
	\item 
	if $u \notin \{x,y\}$ and $v = y$, then perform the same update to $(u,x)$ in $\xi$ as to $(u,v)$ in $\eta$.
\end{itemize}
This corresponds to a relabelling of $x$ and $y$ in $\xi$.
(See \Cref{fig:coup:opening_edges} for an illustration.)

%
%
%
%
%

While the environments are run like this, the environment $\eta$ from the perspective of the walk $X$ looks exactly the same as $\xi$ from the perspective of $Y$, modulo the label difference $x$--$y$.
This allows us to couple $X$ and $Y$, modulo the relabelling.
(See \Cref{fig:coup:walker_moves} for an illustration.)
So at every time, we have one of the following three situations:
\begin{itemize}[noitemsep,topsep=\smallskipamount,label={\ensuremath{\bm\cdot}}]
	\item 
	both $X$ and $Y$ are at some vertex $z \notin \{x,y\}$;
	
	\item 
	$X$ is at $x$ and $Y$ is at $y$;
	
	\item 
	$X$ is at $y$ and $Y$ is at $x$.
\end{itemize}

Observe that this defines a genuine Markovian coupling.
When the systems $(X,\eta)$ and $(Y,\xi)$ start from $(x,\eta_0)$ and $(y,\xi_0)$, respectively, we denote this coupling $\bp_{(x,\eta_0),(y,\xi_0)}$. 
\end{defnt}

%
%
%
%
%

We now describe how to couple two processes $(X,\eta)$ and $(Y,\xi)$, when the environments are initially the same, but the walks are not necessarily isolated.
(This is the set-up of \Cref{res:coup:same-env}.)
In the below algorithm, we define a time $\tauc$ at which $(X,\eta)$ and $(Y,\xi)$ agree.

\begin{enumerate}
	\item \label{i}
	Run the environments together (in the natural coupling, without any relabelling) and the walks independently until the two walks are jointly isolated, ie until time
	\[
		\tau_0
	=
		\inf\brb{ t \ge 0 \mid d_{\tau_0}^\eta(X_{\tau_0}) = 0 = d_{\tau_0}^\xi(Y_{\tau_0}) }.
	\]
	Note that $\eta_{\tau_0} = \xi_{\tau_0}$.
	Write $x = X_{\tau_0}$ and $y = Y_{\tau_0}$.
	
	\item \label{ii}
	Set $k = 1$.
	Use the coupling from \Cref{defn:coupling}: run until $x$ becomes non-isolated in $\eta$ (and hence $y$ becomes non-isolated in $\xi$), and then on until the first time after that both $x$ and $y$ are isolated (in both $\eta$ and $\xi$); call this first time $\sigma_k$ and the final time $\tau_k$.
	That is, set
	\[
		\sigma_k
	=
		\inf\brb{ t \ge \tau_0 \mid d_t^\eta(x) > 0 }
	\Quad{and}
		\tau_k
	=
		\inf\brb{ t \ge \sigma_k \mid d_t^\eta(x) = 0 = d_t^\eta(y) };
	\]
	by the relabelling of the coupling, in the above definition we could swap $(X,\eta) \leftrightarrow (Y,\xi)$ and the times would be the same.
	
	\item \label{iii}
	If $X_{\tau_k} \notin \{x,y\}$, then $X_{\tau_k} = Y_{\tau_k}$ (and vice versa).
	In this case, we have successfully coalesced the full processes.
	We then stop, setting $\tauc = \tau_k$. (Also set $K = k$.)
	
	Otherwise, we have $\{X_{\tau_k}, Y_{\tau_k}\} = \{x, y\}$. Then the walks are not at the same vertex, but the environments are in the same state and the walks are jointly isolated.
	Hence we can increment $k \to k+1$ and return to Step \ref{ii}.
	By symmetry, assume $X_{\tau_k} = x$ and $Y_{\tau_k} = y$.
\end{enumerate}

This means that, once we have the walks simultaneously isolated in the same environment, we can control their evolution very carefully.
It will be straightforward to see that the probability that Step \ref{ii} `succeeds', ie ends with $X_{\tau_k} = Y_{\tau_k}$ is $\tfrac12 + \oh1$; hence we run Step \ref{ii} an order 1 number of times to couple.
It remains to analyse how long Steps \ref{i} and \ref{ii} take (Step \ref{iii} only makes definitions, and so takes no time to `run').
Step \ref{i} is given precisely by \Cref{res:2:tau>t;H}.

\begin{lem}
	The probability Step \ref{ii} ends with $X_{\tau_k} = Y_{\tau_k}$ is at least $\tfrac13$, and hence $K \preccurlyeq \Geo(\tfrac13)$.
\end{lem}

The intuition behind this lemma is as follows.
Consider just $X$ and the vertex $x = X_{\tau_0}$.
In the (random) time interval between when $x$ has degree 1 and when it becomes isolated, the component of $x$ has size at least 2, and moreover $1/\mu \gg n$, so the walk takes a large number of steps if it is in this component, and so has probability $\tfrac12 + \oh1$ of being at $x$ when $x$ becomes isolated.
This heuristic is made rigorous in the proof of \Cref{res:ind:prob<2/3}, where an analogous claim is considered.
From this we deduce that $K \preccurlyeq \Geo(\tfrac13)$.

We still need to find the distribution of $\tau_1 - \tau_0$. (Note that $(\tau_k - \tau_{k-1})_{k\ge1}$ are iid.)

\begin{lem}
	There exists a constant $C$ so that,
	for all $L \ge 1$ and all $n$ sufficiently large,
	we have
	\[
		\prb{ \tau_1 - \tau_0 > C L / \mu } \le e^{-L}.
	\]
\end{lem}

The intuition behind this lemma is as follows.
There are $2n-3$ edges incident to $\bra{x,y}$, in the complete graph; in equilibrium, in expectation $p(2n-3) \approx 2\lambda$ will be open.
The number of edges open is well approximated by a birth and death chain with birth rate $\lambda \mu$ and death rate $\mu$. For this chain, the return time to 0 has mean order $1/\mu$ and an exponential tail.
This heuristic is made rigorous in the proof of \Cref{res:ind:exp-decay-d0}, where a similar claim is considered -- there only one vertex is considered, and some random number of edges are open initially; the same argument applies here.
To go from $\tau_1 - \tau_0$ to $\tau_k - \tau_0 = \sum_{\ell=1}^k (\tau_\ell - \tau_{\ell-1})$, we use the simple fact that if a random variable has an exponential tail, then so does a sum. This can be proved by applying Chernoff; cf the proof of \Cref{res:ind:exp-decay-return}.
(Note that these statements and proofs do not require the graph to be `good'.)

From these two lemmas, we immediately get the following corollary.

\begin{cor}
\label{res:tauk-tau0}
	There exists a constant $C$ so that, for all $L \ge 1$ and all $n$ sufficiently large, we have
	\[
		\prb{ \tauc - \tau_0 > CL/\mu } \le e^{-L}.
	\]
\end{cor}

Now that we know bounds on both $\tau_0$ and $\tauc - \tau_0$, given by \Cref{res:2:tau>t;H,res:tauk-tau0} respectively, the bound on $\tauc$ in \Cref{res:coup:same-env} follows immediately.
(For the application of \Cref{res:2:tau>t;H}, recall from \Cref{defn:very-good-graph} that $\pr[\eta_0]{ \Hbad{0,n/\mu} \le n^{-1} } \le n^{-1}$ for all $\eta_0 \in H$.)
Note also that we may assume $\mu t \ge 3 \log_{(M)}\! n$, else the claim holds trivially.

\section{Invariant Initial Environment}
\label{sec:iiie}

In this section we prove \Cref{res:iiie:statement}, which concerns the case where we draw $\eta_0$ according to $\pi_\ER$ and set $X_0 = 1$; the reader should recall the precise statement.
Throughout this entire section we consider the measure $\pr[1,\ER]{\cdot}$; for ease of notation, we drop the subscript and just write $\pr{\cdot}$.


\begin{Proof}[Proof of \Cref{res:iiie:statement} \normalfont{(\emph{lower bound})}]
Suppose the walk starts from an isolated vertex: it cannot have mixed before an incident edge opens. We make this idea precise and rigorous. We have
\[
	\prb{d_0(1) = 0} = (1-p)^n = e^{-\lambda} \rbb{1 - o(1)}.
\]

Let $\tau$ be the first time an edge incident to $X_0 = 1$ opens. By counting edges and their respective rates, we see that $\tau \sim \me(\lambda\mu(1-1/n)) \succcurlyeq \me(\lambda\mu)$.
Let $T = \lambda \mu t$, and observe that
\[
	\prb{ \tau > t \mid d_0(1) = 0 } \ge \prb{\me(\lambda\mu) > T/(\lambda\mu)} = e^{-T}.
\]

On the event $\{d_0(1)=0\} \cap \{\tau > t\}$, we have $X_t = 1$ (in fact $X_s = 1$ for all~$s \le t$), and hence
\[
	\tv{\prb{ X_t \in \cdot } - \pi_U}
\ge
	\prb{X_t = 1} - \pi_U(1)
\ge
	\prb{ d_0(1) = 0, \, \tau > t } - \pi_U(1)
\ge
	e^{-(T+\lambda)} - \tfrac1n.
\]
We desire $T$ so that $\tv{\prb{ X_t \in \cdot } - \pi_U} \ge \epsilon$. By the above, we may take
\(
	T = - \log(\eps+1/n) - \lambda.
\)
If $\epsilon < e^{-3\lambda}$, then $T \ge \tfrac12 \log(1/\epsilon)$. This proves the lower bound, as $t = T / (\lambda\mu)$.
\end{Proof}

The aim of the remainder of this section is to prove the upper bound in \Cref{res:iiie:statement}; herein we assume that $\mu \le \tfrac23(1+\lambda)^{-1}/n$. To this end, let $\tau$ be the first time our initial vertex is isolated and the walk is not there, ie
\[
	\tau = \inf\brb{ t \ge 0 \mid d_t(X_0) = 0, X_t \neq X_0 }.
\]

The \emph{idea} is that at time $\tau$ we are nearly uniform \emph{and} have lost information about where we started, and so our total variation does not become large in the future. We show this rigorously.

\begin{prop}
\label{res:ind:TV}
	For all~$n$ and all~$t$, we have
	\[
		\tvb{\prb{ X_t \in \cdot } - \pi_U} \le \prb{X_t = 1, \, \tau \le t \wedge (n/\mu)} + \prb{\tau > t \wedge (n/\mu)} + \tfrac1n.
	\]
\end{prop}

\begin{Proof}
Note that by construction and the symmetry of the graph, at all times $t \ge 0$ we must have that $\pr{X_t = x}$ is constant over $x \in V \setminus \{1\} = \{2,...,n\}$: define $\rho_t = \pr{X_t = 1}$; then $\pr{X_t = x} = (1-\rho_t)/(n-1)$ for all~$x \in \{2,...,n\}$. We then have
\begin{align*}
	2 \, \tvb{\prb{ X_t \in \cdot } - \pi_U}
&=
	(n-1) \absB{ \tfrac{1-\rho_t}{n-1} - \tfrac1n }+ \absB{ \rho_t - \tfrac1n }
=
	2 \absb{ \rho_t - \tfrac1n }
\le
	2 \rbb{ \rho_t + \tfrac1n }.
\end{align*}
Decomposing according to the event $\{\tau > t \wedge (n/\mu)\}$ completes the proof.
	%
\end{Proof}

\begin{prop}
\label{res:ind:tau-tail}
	There exists a constant $C$ so that, for all~$K \ge 2$ and all~$n$ sufficiently large, we have
	\[
		\prb{\tau > C K/\mu} \le e^{-K}.
	\]
\end{prop}

Before we prove this, we define some preliminary notation, and then state three claims. First, let $\sigma_0 = \sigma_1' = 0$, and let $\sigma_1$ be the first time the initial vertex, 1, becomes isolated, ie
\[
	\sigma_1 = \inf\brb{ t \ge 0 \mid d_t(1) = 0 },
\]
and for $i \ge 1$ define inductively
\begin{gather*}
	\sigma'_{i+1}  = \inf\brb{ t \ge \sigma_i      \mid d_t(1) > 0 },
\quad
	\sigma_{i+1}   = \inf\brb{ t \ge \sigma'_{i+1} \mid d_t(1) = 0 }
\\
	\text{and}
\quad
	\sigma''_i = \inf\brb{ t \le \sigma_i  \mid d_s(1) = 1 \, \forall  s \in [t,\sigma_i) }.
\end{gather*}
In words, $\sigma_i$ is the $i$-th time the vertex 1 becomes isolated, $\sigma_i'$ is the first time after this that it becomes non-isolated and $[\sigma_i'',\sigma_i)$ is the interval in which it is degree~1 immediately before becoming isolated for the $i$-th time. By the memoryless property, $\sigma'_i - \sigma_{i-1} \sim^\iid \me(\lambda\mu(1-1/n))$.

Define $\tau_i := \sigma_i - \sigma_{i-1}$ for $i \ge 1$; then $\tau_1$ is the time it takes to become isolated initially, and, for $i \ge 2$, $\tau_i$ is the time between the $(i-1)$-st and $i$-th times we become isolated. Note that the random variables $\{\tau_i\}_{i\ge2}$ are all independent and identically distributed.

We now state three lemmas which we use to deduce \Cref{res:ind:tau-tail}.

\begin{lem}
\label{res:ind:exp-decay-d0}
	There exists a constant $C$ so that, for all~$K \ge 1$ and all~$n$ sufficiently large, we have
	\[
		\prb{\tau_1 > C K /\mu} \le e^{-K}.
	\]
\end{lem}

\begin{lem}
\label{res:ind:exp-decay-return}
	There exists a constant $C$ so that, for all~$K \ge 2$ and all~$n$ sufficiently large, we have
	\[
		\prb{\textstyle \sum_{i=2}^K \tau_i > C K /\mu} \le e^{-K}.
	\]
\end{lem}

\begin{lem}
\label{res:ind:prob<2/3}
	For all~$n$ sufficiently large and all~$i \ge 1$, we have
	\[
		\prb{X_{\sigma_i} = 1 \mid X_{\sigma_j} = 1 \, \forall  j < i} \le \tfrac23.
	\]
\end{lem}

We now show how to conclude our tail bounds on $\tau$ from these three lemmas.

\begin{Proof}[Proof of \Cref{res:ind:tau-tail}]
	Consider an integer $K \ge 2$. 	\Cref{res:ind:prob<2/3} tells us that
	\[ \textstyle
		\prb{X_{\sigma_i} = 1 \, \forall  i \le K}
	=
		\prod_{i=1}^K \prb{X_{\sigma_i} = 1 \mid X_{\sigma_j} = 1 \, \forall  j < i}
	\le
		\rbb{2/3}^K.
	\]
	Combining this with \Cref{res:ind:exp-decay-d0} and \Cref{res:ind:exp-decay-return} tells us that
	\[
		\prb{\sigma_K \le C' K / \mu, \ \exists \ k \le K \ST X_{\sigma_k} \neq 1}
	\ge
		1 - e^{-K} - e^{-K} - \rbr{3/2}^{-K}
	\]
	for a suitably large constant $C'$. From this we deduce our claim.
\end{Proof}

To complete the proof of our tail bound, it remains only to prove our three lemmas; we do this at the end of the section. For now, we turn to upper bounding $\pr{X_t = 1, \, \tau \le t \wedge (n/\mu)}$.

\begin{lem}
\label{res:ind:rho-bound}
	There exists a constant $C$ so that, for all~$n$ sufficiently large and all~$t$, we have
	\[
		\prb{X_t = 1, \, \tau \le t \wedge (n/\mu)} \le C / n.
	\]
\end{lem}

\begin{Proof}
Write $\mi_t$ for the set of isolated vertices at time $t$. Write $s = t \wedge (n/\mu)$. First we lower bound the number of isolated vertices at time $\tau$ on the event $\{\tau \le s\}$.
From \Cref{res:ER_very-good},
\[
	\prb{ |\mi_\tau \setminus \{X_\tau\}| \le \tfrac12 c_* n, \, \tau \le s }
=
	\Ohb{n^{-2}}.
\]

By the symmetry of the complete graph, we must have that $\pr{X_t = x \mid \mf_\tau}$ is constant over $x \in \mi_\tau \setminus \{X_\tau\}$ on the event $\{\tau \le s\}$; let $\xi_t$ be this (random) value. ($\xi_t$ is an $\mf_\tau$-measurable random variable.) Now, by construction of $\tau$, we have that $X_0 = 1 \in \mi_\tau \setminus \{X_\tau\}$. This says that
\[
	\xi_t = \prb{X_t = 1 \mid \mf_\tau} \oneb{ \tau \le s }
	\Quad{and hence}
	\prb{X_t = 1, \, \tau \le s} = \exb{\xi_t}.
\]
It remains to bound $\ex{\xi_t}$, which we now do. Note that we have
\[
	1
\ge
	\prb{ \tau \le s }
\ge
	\exb{ |\mi_\tau \setminus \{X_\tau\}| \cdot \prb{ X_t = 1 \mid \mf_\tau } \oneb{ \tau \le s } }.
\]
Letting $A = \{|\mi_\tau \setminus \{X_\tau\}| \ge \tfrac12 c_* n\}$, we have $\pr{A^c, \tau \le s} = \Oh{n^{-2}}$, as above. Hence
\begin{align*}
	1
&\ge
	\exb{ |\mi_\tau \setminus \{X_\tau\}| \cdot \prb{ X_t = 1 \mid \mf_\tau } \oneb{ \tau \le s } \one{A} }
\\&
\ge
	\tfrac12 c_* n \exb{ \xi_t \one{A} }
\ge
	\tfrac12 c_* n \rbb{ \exb{ \xi_t} - \prb{A^c, \tau \le s } }.
\end{align*}
Rearranging completes the proof:
\[
	\prb{ X_t = 1, \, \tau \le s }
=
	\exb{\xi_t}
\le
	(\tfrac12 c_* n)^{-1} + \Oh{n^{-2}}
\le
	3c_*^{-1} / n.
\qedhere
\]
\end{Proof}

We can now give the proof of the upper bound in \Cref{res:iiie:statement}.
\begin{Proof}[Proof of \Cref{res:iiie:statement} \normalfont{(\emph{upper bound})}]
	\Cref{res:ind:rho-bound} says that, for all~$t$, we have
	\[
		\prb{ X_t = 1, \, \tau \le t \wedge (n/\mu) } \le C'/n,
	\]
	for a constant $C'$.
	Hence we have
	\[
		\tvb{ \prb{ X_t \in \cdot } - \pi_U }
	\le
		\prb{ \tau > t \wedge (n/\mu) } + (C'+1)/n.
	\]
	Observe that this upper bound is (weakly) monotone-decreasing in $t$, and \Cref{res:ind:tau-tail} gives us a constant $C$ so that
	\[
		\tvb{ \prb{ X_t \in \cdot } - \pi_U }
	\le
		\eps^2 + (C'+1)/n
	\le
		\eps
	\Quad{when}
		t = 2 C \logeps / \mu.
	\]
	Hence we deduce that $\tmix(\epsilon) \le 2 C \log(1/\epsilon)/\mu$.
\end{Proof}

\begin{Proof}[Proof of \Cref{res:ind:exp-decay-d0}] \label{pf:exp-decay-d0}
	Write $d_t = d_t(1)$. Also rescale time by $\mu$, so as to remove the $\mu$ factors from the workings. Observe that the jump-rates of $d$ are as follows:
	\begin{align*}
		&k \to k+1 \Quad{at rate} q_+(k) = (n-1-k)p = (\lambda - \lambda(1 + k)/n);
	\\
		&k \to k-1 \Quad{at rate} q_-(k) = k(1-p) = (k - \lambda k/n).
	\end{align*}
	Let $q(k) = q_+(k) + q_-(k)$, and observe that $q(k) \ge q(0) \ge 1$ for all~$k \ge 0$. We now couple $d$ with an auxiliary process $d'$, which has rate-1 jumps. Above $3\lambda$, $d'$ has probability $\tfrac23$ of going up and $\tfrac13$ of going down; below $3\lambda$, it has the same probabilities as $d$, ie $q_+(k)/q(k)$ for up and $q_-(k)/q(k)$ for down. Set $d'_0 = d_0$, and write $\tau'_1$ for the hitting time of 0 by $d'$. We then have $\tau_1 \preccurlyeq \tau'_1$.
	
	Note that once $d'$ reaches $\ceil{3\lambda}$, it moves directly to 0 (in $\ceil{3\lambda}$ steps) with probability bounded away from 0. The hitting time of $\ceil{3\lambda}$ is that of a biased simple random walk. Since $d_0 \sim \Bin(n-1,\lambda/n)$, we may assume that $d_0 \le C K$ for some sufficiently large constant $C$ with a penalty $e^{-K}$ to the probability. Given this, we see that the hitting time of $\ceil{3\lambda}$ has mean $\Th1$ and an exponential tail. Once $d'$ hits $\ceil{3\lambda}$, we perform a geometric number of excursions, the length of which have an exponential tail.
	Hence $\tau'_1$ has mean $\Th1$ and an exponential tail.
\end{Proof}

\begin{Proof}[Proof of \Cref{res:ind:exp-decay-return}] \label{pf:exp-decay-return}
	Again, drop the $\mu$ factors. Note that $\tau_i \preccurlyeq \tau_1$, and $\tau_1$ has mean $\Th1$ with an exponential tail. Since the $\tau_i$ are independent, we then apply the Chernoff bound to a sum of $K$ independent $\tau_1$ random variables to deduce the lemma.
\end{Proof}

\begin{Proof}[Proof of \Cref{res:ind:prob<2/3}] \label{pf:prob<2/3}
Fix $i \ge 1$. For $t \in (\sigma''_i,\sigma_i)$ we have that $d_t(1) = 1$; write $x_i$ for the neighbour of 1 in the interval $(\sigma''_i,\sigma_i)$. Note that all the $\sigma$-times depend only on the environment, not also on the walk. We describe a coupling between $X$ and an auxiliary walk $X'$ which is confined to the pair $\{1,x_i\}$. The coupling will have the property that
\[
	\prb{ X_{\sigma_i}  = 1 \mid X_{\sigma_j} = 1 \, \forall  j < i }
\le
	\prb{ X'_{\sigma_i} = 1 \mid X_{\sigma_j} = 1 \, \forall  j < i }.
\]
In particular, $X'$ will be the usual simple random walk on $\{1,x_i\}$, jumping at rate $1/(n-1)$. Thus we shall see that the probability on the right-hand side is `approximately' $\tfrac12$.

We now explicitly define the coupling. Start $X'$ from 1. If both $X$ and $X'$ are at 1, then move them together; if both $X$ and $X'$ are at $x_i$ and $X$ chooses vertex 1 to jump to, then move them together; otherwise let them evolve independently. Observe that, wherever $X$ is at time $\sigma''_i$, we always have for $t \in [\sigma''_i,\sigma_i]$ that $X_t = 1$ implies $X'_t = 1$. Hence our desired inequality holds.

Observe that $X'$ is at 1 if it has taken an even number of steps (and at $x_i$ if odd). Hence
\[
	\prb{ X'_{\sigma_i} = 1 \mid X_{\sigma_j} = 1 \, \forall  j < i }
=
	\prb{ \Po(r) \ \text{is even} }
=
	\tfrac12 \rbb{ 1 + \exb{ e^{-2r} } }
\Quad{where}
	r = (\sigma_i - \sigma''_i)/(n-1).
\]
We have $\sigma_i - \sigma''_i \sim \me( (\lambda+1-3p)\mu )$ by counting edges and rates, and so the lemma follows since
\[
	\exb{ e^{-2r} }
=
	\frac{(\lambda+1-3p)\mu}{ (\lambda+1-3p)\mu + 2/(n-1) }
\le
	\tfrac12(\lambda + 1) \mu n
\le
	\tfrac13.
\qedhere
\]
	%
\end{Proof}


\begin{quote}
	\normalsize
	\textbf{\textsf{Acknowledgements.}}
	The authors would like to thank Andrew Swan for many helpful discussions at the start of this project.
\end{quote}


\renewcommand{\bibfont}{\sffamily}
\printbibliography[heading=bibintoc]


\appendix
\section{Proofs of Erd\H{o}s--R\'enyi Structure Results}
\label{app}

In this section we give the proofs of the Erd\H{o}s--R\'enyi structure results given in \Cref{res:ER_good}.
The results will all be for supercritical Erd\H{o}s--R\'enyi random graphs.
As mentioned before, such results are usually shown to hold with probability $1-\oh1$ (see \cite{B:random-graphs}); here we quantify this $\oh1$.

\medskip

Properties \ref{gg:max-deg} and \ref{gg:num-edges} follow easily from properties of the Binomial distribution. We give these proofs immediately.
Applying Chernoff's bound, one can easily show, for $\lambda > 1$, that
\[
	\prb{ \Bin(n-1,\lambda/n) \ge 3 \lambda K } \le e^{-K}.
\]
Taking $K = 10 \logn$ and applying the union bound, we deduce \ref{gg:max-deg}.
The number of edges in the graph is precisely $\Bin(\binom{n}2,\lambda/n)$, and so we deduce \ref{gg:num-edges} from concentration of the Binomial.

\medskip

The remaining results take more to prove.
A method we shall use repeatedly is to find expectation of various quantities and then show concentration. In particular, we shall use a result on \textit{typical bounded difference} due to Warnke \cite{W:typ-bounded-diff}, which builds on work by McDiarmid \cite{McD:bounded-diff}.

Property \ref{gg:size} is standard, but the proof that we shall give is a nice, straightforward example of the above typical bounded differences approach, and so we give it as a type of `warm-up proof'.

In order to prove \ref{gg:num-deg1-giant}, \ref{gg:removal-edges} and \ref{gg:vtx-far-from-core}, we use \ref{gg:giant-expansion}. As such, we first prove \ref{gg:giant-expansion} directly (not using concentration, with one exception); we then formulate the concentration inequality of Warnke~\cite{W:typ-bounded-diff}, and apply it in various situations.

Write $G \sim \ER(n,p)$ if $G$ has the Erd\H{o}s--R\'enyi distribution on $n$ vertices with edge probability~$p$.
Direct calculation as at the start of the proof of \cite[Theorem 2.14]{FK:intro-random-graphs} tells us that if $G \sim \ER(n,\lambda/n)$ with $\lambda > 1$ fixed then there exist positive constants $c$ and $C$ so that the probability there is a component with size in $[C \logn, c n]$ is $\Oh{n^{-100}}$.

\subsection{Expansion Properties of the Giant}

We now consider \ref{gg:giant-expansion}. We first state precisely what it is that we prove.

\begin{prop}[Expansion Properties of the Giant, \ref{gg:giant-expansion}]
\label{app:giant-expansion}
	Suppose $G \sim \ER(n,\lambda/n)$, with $\lambda > 1$ fixed, and write $\mg$ for its largest component (breaking ties arbitrarily). Then there exists a positive constant $c$ so that $\Phi_\mg \ge c \logn[-2]$ and $\gamma_\mg \ge c \logn[-4]$ with probability $1 - \Oh{n^{-9}}$.
\end{prop}

We prove this via a sequence of lemmas, the proof of which are deferred until after the proof of \Cref{app:giant-expansion}. First, though, we make introductory observations and definitions.

\medskip

Given a graph $G$, write $\delta(G)$ for its minimum degree.
The \textit{core} of $G$, denoted $\mc_G = \mc(G)$, is obtained by removing all isolated vertices and then recursively removing all degree~1 vertices; note that $\delta(\mc_G) \ge 2$.
The \textit{kernel} of $G$, denoted $\mk_G = \mk(G)$, is obtained from the core by replacing each maximal 2-path that joins vertices $u$ and $v$ of degree at least 3 by an edge $(u,v)$, and deleting each isolated cycle. (Note that the kernel may have self-loops and/or multiple edges, and may not be connected.) A loop contributes 2 to the degree of its incident vertex, so $\delta(\mk_G) \ge 3$. Note that if~$G$ is connected, then so are $\mc_G$ and $\mk_G$.
Write $\mg_G = \mg(G)$ for the largest component of $G$.

Given a degree sequence $\bm d$, write $G(\bm d)$ for the set of all (multi-)graphs with degree sequence $\bm d$, and write $U(\bm d)$ for $\Unif(G(\bm d))$; also write $\bm d^*$ for the degree sequence $\bm d$ with the degree~2 vertices removed. Observe that if $G \sim U(\bm d)$, then $\mk_G \sim U(\bm d^*)$.

Also write $G'(\bm d)$ for the set of all simple graphs with degree sequence $\bm d$, and write $U'(\bm d)$ for $\Unif(G'(\bm d))$.
Pittel, Spencer and Wormald showed in \cite[Proposition 1(b)]{PSW:core-unif-dist}, for $G \sim \ER(n,p)$, that $\mc_G \sim U'(\bm d)$ conditional on having degree sequence $\bm d$.
Since this includes conditioning on being simple, we cannot deduce that $\mk_G \sim U'(\bm d^*)$; however, it's known that a realisation of $U(\bm d)$ is simple with order 1 probability under certain conditions (we state and reference this precisely below), and this will be sufficient for us.

The \textit{configuration model} was introduced by Bollob\'as in \cite{B:CM}, and provides a method for constructing a graph with a given degree sequence with distribution that is uniform over all graphs with the given degree sequence. (Bender and Canfield in \cite{BC:CM} also introduced a highly related object; see also Bollob\'as \cite{B:random-graphs}, or van der Hofstad \cite{H:random-graphs-book} for a modern description.)

As some notation, for a graph $H$, write $H = (V_H,E_H)$ and $d_H(x)$ for the degree of a vertex $x$ in $H$. For $S \subseteq H$, write $d_H(S) = \sum_{x \in S} d_H(x)$ for the total degree of~$S$, $\Phi_H(S) = e_H(S)/d_H(S)$, $I_H(S)$ for the number of induced total degree (ie the sum of the degrees in the subgraph induced by $S$, ie twice the number of edges in the subgraph induced by $S$) and $e_H(S)$ for the number of edges in $E_H$ between $S$ and $S^c$.
Then $\Phi_H(S) = 1 - I_H(S)/d(S)$ and
\[
	\Phi_H = \min\brb{ \Phi_H(S) \mid 0 < d_H(S) \le |E_H| }.
\]
Also, for a degree sequence $\bm d$, write $\delta_{\bm d}$ for its minimal degree and $\theta_{\bm d}$ for its average degree.

\begin{lem}
\label{app:U(d)-expander}
	Let $\bm d = (d_1, ..., d_m)$ be a degree sequence with $\delta_{\bm d} \ge 3$ and $\theta_{\bm d} \le m^{1/5}$.
	Suppose that $G \sim U(\bm d)$, and write $\mg$ for its largest component (breaking ties arbitrarily).
	We have
	\[
		\prb{ \Phi(\mg) < 10^{-5} } = \Ohb{ m^{-50} }.
	\]
\end{lem}

The proof of this lemma will closely follow \cite[Lemma 12]{ACF:kernel-expander}, but with a few slight changes to strengthen the probability. We want to convert this from a result about $U(\bm d^*)$ into a result about the kernel of the giant of an Erd\H{o}s--R\'enyi graph. Recall that if $F \sim \ER(n,p)$ then $\mc_F \sim U'(\bm d)$ conditional on having degree sequence $\bm d$, and that $\mk(U(\bm d)) \sim U(\bm d^*)$.


\begin{lem}
\label{app:kernel-size}
	Suppose $G \sim \ER(n,\lambda/n)$ with $\lambda > 1$ a constant, and write $\mg$ for its largest component (breaking ties arbitrarily). There exists a positive constant $c$ so that $|\mk_\mg| \ge c n$ with probability at least $1 - \Oh{n^{-10}}$.
\end{lem}

\Cref{app:kernel-size} is proved using a concentration result; its proof is deferred until the next section.

\medskip

We also need to know that if $G \sim U(\bm d)$ then $\pr{\simple{G}}$ is not too small.

\begin{lem}
\label{app:simplicity-ER}
	Suppose $F \sim \ER(n,\lambda/n)$ with $\lambda > 1$ fixed.
	Let $\bm d$ be a random degree sequence with the distribution of the degree sequence of $\mc_F$.
	Let $G \sim U(\bm d)$.
	Then there exists a positive constant $c$ so that, for all~$n$ sufficiently large, we have
	\[
		\prb{ \simple{G} } \ge c.
	\]
\end{lem}

The above results have been about the kernel, or the core. 
To convert this into a result about the giant, we consider the `decorations' to the kernel.
Consider a graph $G$ and its kernel $\mk_G$. An edge of $\mk_G$ is first expanded into a path in the core, and then trees are hung from each vertex of the path. For an edge $e$ of $\mk_G$, write $D_e$ for the number of vertices added in this `decorating' process; write $D = \max_e D_e$.

Suppose the edge $(x,y)$ in $\mk_G$ has decoration with vertex set $S$, where we do not include $x$ or $y$ in $S$. If $S = \emptyset$, then $(x,y)$ is an edge in $G$; suppose $S \neq \emptyset$. Then $G[S]$, the subgraph of $G$ induced by $S$, is a tree. Moreover, all the (potential) edges between $S$ and $S^c$ are closed in $G$, with the exception of one edge between $S$ and $x$ and one between $S$ and $y$.

\begin{lem}
\label{app:largest-decoration}
	Suppose that $G \sim \ER(n,\lambda/n)$ with $\lambda > 1$ fixed.
	Write $D$ for the size of the largest decoration (as above).
	There exists a constant $C$ so that
	\[
		\prb{ D > C \logn } = \Ohb{ n^{-100} }.
	\]
\end{lem}

Given these two lemmas, we can now prove \Cref{app:giant-expansion}.

\begin{Proof}[Proof of \Cref{app:giant-expansion}]
We first note that the inequality $\tfrac12 \Phi^2 \le \gamma \le 2 \Phi$ holds; this is due to Jerrum and Sinclair \cite{JS:cheeger} and Lawler and Sokal \cite{LS:cheeger}. Hence given $\Phi \ge c \logn[-2]$ it follows that $\gamma \ge c' \logn[-4]$ with $c' = \tfrac12 c^2$. So it suffices to just prove the statement about $\Phi_\mg$.

Let $F \sim \ER(n,\lambda/n)$, with $\lambda > 1$ fixed, and let $\bm d$ be a random degree sequence with the distribution of the degree sequence of $\mc_F$ (ie the core of an Erd\H{o}s--R\'enyi graph). Let $G' \sim U'(\bm d)$, $G \sim U(\bm d)$ and $H \sim U(\bm d^*)$.
We first show that
\[
	\prb{ \Phi_{\mg(\mk(F))} < 10^{-5} } = \Ohb{n^{-40}}.
\]
Let $c_1$ be the constant from \Cref{app:U(d)-expander}, $c_2$ from \Cref{app:kernel-size} and $c_3$ from \Cref{app:simplicity-ER}. Then
\GAP{3}
\begin{align*}
&	\prb{ \Phi_{\mg(\mk(F))} < 10^{-3} }
=
	\prb{ \Phi_{\mg(\mk(G'))} < 10^{-3} }
\\&\gap
\le
	\prb{ \Phi_{\mg(\mk(G))} < 10^{-3}, \, |\mk(G)| \ge c_2 n \mid \simple{G} }
+	\prb{ |\mk(G')| < c_2 n }
\\&\gap
\le
	\prb{ \Phi_{\mg(\mk(G))} < 10^{-3}, \, |\mk(G)| \ge c_2 n } / \prb{ \simple{G} }
+	\prb{ |\mk(G')| < c_2 n }
\\&\gap
\le
	(c_1 n)^{-50} / c_3 + n^{-100}
=
	\Ohb{ n^{-50} }.
\label{eq:kernel-expander}
\nt
\end{align*}

We now show how to move from the isoperimetric constant of the kernel to that of the graph.
To do this, we first set up some notation.
Let $G \sim \ER(n,\lambda/n)$ with $\lambda > 1$ fixed;
write $\mg$ for the giant, $\mc$ for its core and $\mk$ for its kernel.
For an edge $(x,y) \in E_\mk$, write $\me_{x,y}$ for the set of vertices in $\mg \setminus \mk$ that are the decoration of the edge $(x,y)$;
For a vertex $z \in \mc$, write $\mt_z$ for the dangling tree hung from $z$ and $\mt'_z = \mt_z \cup \{z\}$.
For $x \in \mk$, write $\md_x = \mt_x \cup (\cup_{y \sim^\mk x} \me_{x,y})$; note that $x \notin \md_x$.
Define the following properties for sets $S \subseteq \mg$, which we denote $(\dagger)$ and $(\ddagger)$, respectively:
\begin{gather*}
	\forall \, x \in S \cap \mk \text{ we have } \md_x \subseteq S.
\label{prop:fill-in}
\tag{\ensuremath\dagger}
\\
	\forall \, x \in S \setminus \mk \ \exists \, y \in S \cap \mk \text{ with } x \in \md_y;
\label{prop:non-isol}
\tag{\ensuremath\ddagger}
\end{gather*}

We do not consider $\mg$ to be a random graph, but rather consider each realisation of it. We assume that it has the following properties:
	there exists a constant $c$ so that $|\mg| \ge c n$;
	$|\mk| > 0$;
	writing $\Delta$ as an upper bound for the maximum degree and for the largest decoration, there exists a constant $C$ so that we may take $\Delta \le C \logn$.

First we show how to restrict attention to sets $S$ with $\mk \not\subseteq S$. Indeed, suppose $S$ has $\mk \subseteq S$ and $d(S) \le |E_\mg|$.
We claim that there exist (at least) $3 \Delta^2$ vertices $z \in \mc$ such that $\mt_z \not\subseteq S$ but $\mt_z \subseteq \md_x$ for some $x \in \mk \subseteq S$.
Indeed, $\cup_{x \in \mk} \md_x = \mg$ and, and so by the definition of $\md_x$ were this claim not the case then $S = \mg \setminus T$ where $T$ is a union of subsets of at most $2 \Delta^2$ trees $\mt_z$, and hence has size at most $2 \Delta^3$, and so $d(T) \le 6 \Delta^3$ (since the components are trees). But $d(S) \le |E_\mg|$, and so $d(T) \ge |E_\mg| \ge cn$, contradicting $d(T) \le 6 \Delta^2 \asymp \logn[3]$.

Given this property, define a new set, $S'$, from $S$ in the following way. First choose an element $x \in \mk$ and remove it from our set. Since $d_\mg(x) \le \Delta$ and the $|\md_x| \le \Delta$, this removes at most $\Delta^2$ vertices $z \in \mc$ with the above property; hence at least $2 \Delta^2$ remain. Choose $2 \Delta$ such $z \in \mc$, say $\{z_1, ..., z_{2\Delta}\}$. To each $z_i$, associate a vertex $x_i \in \mk$ with $z_i \in \md_{x_i}$; note that the $x_i$ need not be distinct. Now add to our set all the vertices of $\md_{x_i}$ for each $i$. So $S' = S \cup (\cup_{i=1}^{2\Delta} \md_{x_i}) \setminus x$.

We now compare the outer boundary and edge-count of $S$ and $S'$.
Removing the first vertex can remove at most $d_\mg(x) \le \Delta$ from the number of edges in $S$, but adding in the decorations adds at least $2 \Delta$ (as $2 \Delta$ vertices are added). Hence $d(S') - d(S) \ge 2 \Delta$ (since $d(S)$ is twice the number of edges in $S$).
Removing the first vertex can add at most $d_\mg(x) \le \Delta$ edges to the outer boundary, but adding the decorations removes at least $2 \Delta$ (since filling such a tree $\mt_z$ removes at least 1). Hence $|\partial_\mg S| - |\partial_\mg S'| \ge \Delta$.
In particular, $\Phi_\mg(S') \le \Phi_\mg(S)$.

Of course, it may be the case that $d(S') > |E_\mg|$. However, we note that $d(S') \le d(S) + 2 \Delta^2 \le \tfrac32 |E_\mg|$. To allow such $S'$, we define $\Phi_\mg'$ as follows:
\[
	\Phi_\mg' = \min\brb{ S \subseteq \mg \mid \rbr{0 < d(S) \le |E_\mg|} \textsc{ or } \rbr{|E_\mg| \le d(S) \le \tfrac32|E_\mg| \text{ and } \mk \not\subseteq S} }.
\]
Clearly $\Phi_\mg \ge \Phi_\mg'$.
The above considerations show that we can then exclude $S$ with $\mk \subset S$:
\[
	\Phi_\mg' = \min\brb{ S \subseteq \mg \mid 0 < d(S) \le |E_\mg|, \, \mk \not\subseteq S }.
\]
From now on we shall assume that sets $S$ in question do not have $\mk \subseteq S$, and hence $S^c \cap \mk \neq \emptyset$.

Now consider restricting attention to sets satisfying $(\dagger)$ or both $(\dagger)$ and $(\ddagger)$:
\begin{align*}
	\Phi_\mg^\dagger
&=
	\min\brb{ \Phi_\mg(S) \mid 0 < d(S) \le \tfrac32|E_\mg|, \, \mk \not\subseteq S, \, \text{$S$ satisfies $(\dagger)$} };
\\
	\Phi_\mg^\ddagger
&=
	\min\brb{ \Phi_\mg(S) \mid 0 < d(S) \le \tfrac32|E_\mg|, \, \mk \not\subseteq S, \, \text{$S$ satisfies $(\dagger,\ddagger)$} }.
\end{align*}
Since decorations are trees, since they are a (core-)path with dangling trees hung from each vertex of the path, we have $\Phi_\mg = \Phi_\mg^\dagger$. (Note that a set $S$ with $S \cap \mk = \emptyset$ satisfies $(\dagger)$ vacuously, and a set $S$ satisfying $(\ddagger)$ must have $S \cap \mk \neq \emptyset$.)

We now show that it suffices to consider `connected' $S$, ie $S \subseteq \mg$ so that the induced graph~$G[S]$ is connected. Indeed, suppose that $S = S_1 \cup S_2$ with $S_1 \cap S_2 = \emptyset$ and no edges between $S_1$ and $S_2$. Then $d(S) = d(S_1) + d(S_2)$ and $|\partial S| = |\partial S_1| + |\partial S_2|$.
Now, for any $a,b,c,d > 0$, we have
\[
	\frac{a+b}{c+d} \ge \min\brb{a/c, b/d},
\Quad{and hence}
	\Phi(S) \ge \min\brb{ \Phi(S_1), \Phi(S_2) },
\]
from which it follows that we may assume $S$ is connected.
Indeed, suppose $a/c < b/d$; then
\[
	\frac{a+b}{c+d}
\ge
	\frac{a+ad/c}{c+d}
=
	\frac{a(c+d)}{c(c+d)}
=
	\frac ac.
\]

Recall that the largest decoration is of size at most $\Delta$.
We claim that
\[
	\Phi_\mg = \Phi_\mg^\dagger = \min\brb{ \Phi_\mg^\ddagger, 1/\Delta }.
\label{eq:app:Phi-ddagger}
\nt
\]
Indeed, if $S$ is connected and does not satisfy $(\ddagger)$, then, since decorations are connected at kernel-vertices, it can contain vertices of at most one decoration and $S \cap \mk = \emptyset$. Since $|\partial S| \ge 1$ for all $S \notin \{\emptyset,\mg\}$ and decorations are trees, such $S$ have $\Phi_\mg(S) \ge 1/\Delta$.
This proves the claim.

Observe that for $S$ satisfying $(\dagger,\ddagger)$ we have
\[
	\abs{ \partial_\mg S^c }
=
	\abs{ \partial_\mg S }
=
	\abs{ \partial_\mk (S \cap \mk) }
=
	\abs{ \partial_\mk (S^c \cap \mk) }.
\label{eq:app:SK-out}
\nt
\]

We claim now that a set $S$ satisfying $(\ddagger)$ satisfies
\begin{subequations}
\[
	d_\mg(S) \le \Delta d_\mk(S \cap \mk) + \Delta |S \cap \mk| \le 2 \Delta \, d_\mk(S \cap \mk).
\label{eq:app:dGS}
\nt
\]
Indeed, by associating vertices $z \in \me_{x,y}$ with $(x,y) \in E_\mk$, we get the first term; with the second term coming from the dangling tree attached to vertices of the kernel. We also claim that if $S$ satisfies both $(\dagger,\ddagger)$, then $S^c$ satisfies $(\ddagger)$ unless $\mk \subseteq S$, and hence the same inequality, ie
\[
	d_\mg(S^c) \le \Delta d_\mk(S^c \cap \mk) + \Delta |S^c \cap \mk| \le 2 \Delta \, d_\mk(S^c \cap \mk).
\label{eq:app:dGSc}
\nt
\]
\end{subequations}
Indeed, take $x \in S^c \setminus \mk$.
First suppose that $x \in \md_{y,z}$ for some $(y,z) \in E_\mk$. We must then have $y,z \in S^c$, else by $(\dagger)$ we would have $x \in S$. Hence we have $x \in \md_y$ with $y \in S^c \cap \mk$.
The other case is to suppose that $x \in \mt_y$ for some $y \in \mk$. Again, we must have $y \in S^c$, else by $(\dagger)$ we would have $x \in S$. Hence we have $x \in \mt_y \subseteq \md_y$ with $y \in S^c \cap \mk$.
This proves the claim.

Instead of considering $\Phi_\mg^\ddagger$, for a graph $H$ we consider
\[
	\Psi_H = \min\brb{ \Psi_H(S) \mid S \notin \{\emptyset, H\} }
\Quad{where}
	\Psi_H(S) = \frac{ 2 |\partial_H S| |E_H| }{ d_H(S) d_H(S^c) }.
\]
Since any $S \subseteq H$ has $d_H(S^c) \le 2 |E_H|$, we have $\Psi_H \ge \Phi_H$.
If we define $\Phi_H$ by considering sets~$S$ with $0 < d(S) \le 2(1-c)|E_H|$, then $c \Psi_H \le \Phi_H$.
Using the above expression for $|\partial_\mg S|$ and inequalities for $d_\mg(S)$ and $d_\mg(S^c)$, namely \eqref{eq:app:SK-out}, \eqref{eq:app:dGS} and \eqref{eq:app:dGSc}, for $S$ satisfying $(\dagger,\ddagger)$, we have
\[
	\Psi_\mg(S)
=
	\frac{ 2 |\partial_\mg S| |E_\mg| }{ d_\mg(S) d_\mg(S^c) }
\ge
	\frac{ 2 |\partial_\mk (S \cap \mk)| |E_\mk| }{ 4\Delta^2 \, d_\mk(S \cap \mk) d_\mg(S^c \cap \mk) }
\ge
	2 \Psi_\mk(S \cap \mk) / \Delta^2.
\]
We have $S \cap \mk, S^c \cap \mk \notin \{\emptyset,\mk\}$ by recalling the following:
	we are assuming $\mk \not\subseteq S$, which implies $S \cap \mk \neq \mk$, which in turn implies $S^c \cap \mk \neq \emptyset$;
	also $S$ satisfies $(\ddagger)$, which implies $S \cap \mk \neq \emptyset$, which in turn implies $S^c \cap \mk \neq \mk$.
Hence we deduce that
\[
	\Psi_\mg^\ddagger = \min\brb{ \Psi_H(S) \mid S \notin \{\emptyset, \mg\}, \, \text{$S$ satisfies $(\dagger,\ddagger)$} }
\ge
	2 \Psi_\mk / \Delta^2
\ge
	2 \Phi_\mk / \Delta^2.
\]
In our calculation of $\Phi_\mg$, we are using the condition $0 < d(S) \le |E_\mg|$, and hence
\[
	\Phi_\mg^\ddagger \ge \tfrac12 \Psi_\mg^\ddagger \ge \Phi_\mk / \Delta^2.
\]
We have shown that $\pr{ \Phi_\mk < 10^{-3} } = \Oh{n^{-500}}$, and hence the proposition follows.
\end{Proof}

%
%
%

It remains to give the deferred proofs of \Cref{app:U(d)-expander,app:simplicity-ER,app:largest-decoration}. (Recall that the proof of \Cref{app:kernel-size} is deferred further, until the next section.)

\begin{Proof}[Proof of \Cref{app:U(d)-expander}]
For any connected graph $H$, any $S \notin \{\emptyset,H\}$ has $e_H(S) \ge 1$, and hence $\Phi_H(S) \ge 1/d_H(S)$.  Thus we have
\[
	\Phi_H \ge \min\brb{ \Phi_H', \, 10^{-5} }
\Quad{where}
	\Phi_H' = \min\brb{ \Phi_H(S) \mid S \subseteq H, \, 10^5 \le d(S) \le |E_H| }.
\]
Since $\mg$ is a connected component of $G \sim U(\bm d)$, we can lower bound
\[
	\Phi_\mg \ge \min\brb{ \Phi_\mg', \, 10^{-5} }
\Quad{where}
	\Phi_\mg' = \min\brb{ \Phi_\mg(S) \mid S \subseteq G, \, 10^5 \le d(S) \le |E_G| },
\]
where we define $\Phi_\mg(S) = 0$ if $S \supseteq \mg$. The advantage of doing this is that the set over which we are minimising is no longer a random set.
It will be easier to work with $G$, rather than $\mg$, since $G$ has a nice form.
Fix $\beta = 99/100$.

\medskip

\emph{Consider `small' sets.}
For $s \in \bn$, define the collection of sets
\begin{align*}
	\ms_0(s)
&=
	\brb{ S \mid |S| = s, \, 10^5 \le d(S) \le (\theta m)^{1/4} }.
\end{align*}
Let $N_0(s)$ be the expected number of sets $S \in \ms_0(s)$ with $I_G(S) \ge \beta d_G(S)$.
As in \cite[Lemma 12]{ACF:kernel-expander},
\begin{align*}
	N_S(s)
&
=
	\sum_{S \in \ms_0(s)}
	\binom{d(S)}{\beta d(S)}
	\binom{\theta m / 2}{\beta d(S)/2}
	\bigg/
	\binom{\theta m}{\beta d(S)}
\\&
\lesssim
	\sum_{S \in \ms_0(s)}
	\binom{d(S)}{\beta d(S)}
	\rbbb{ \frac{\beta d(S)}{\theta m} }^{\beta d(S) / 2}
	\rbbb{ 1 - \frac{\beta d(S)}{\theta m} }^{(\theta m - \beta d(S))/2}
\\&
\lesssim
	\sum_{S \in \ms_0(s)}
	\rbbb{ \frac{d(S) e}{\beta d(S)} }^{\beta d(S)}
	\rbbb{ \frac{\beta d(S)}{\theta m} }^{\beta d(S)/2}
\\&
=
	\sum_{S \in \ms_0(s)}
	\rbbb{ \frac{e^2 d(S)}{\beta \theta m} }^{\beta d(S) / 2}
\lesssim
	\sum_{S \in \ms_0(s)}
	\rbb{ 5 m^{-1} }^{(1.1) d(S) / 3},
\end{align*}
since $\theta \ge 3$ and $\beta = 99/100$.
Note that there are at most $m^s$ terms in the sum.
\begin{itemize}
	\item []
	\emph{Consider $s \ge 10^4$.}
	We use the bound $d(S) \ge 3|S| = 3s$, which is valid for all sets. Hence
	\[
		N_0(s)
	\lesssim
		m^s \, \rbb{\tfrac14 m}^{-(1.1)s}
	\lesssim
		m^{-(0.05)s}
	\lesssim
		m^{-50}.
	\]
	
	\item []
	\emph{Consider $s \le 10^4$.}
	We use the bound $d(S) \ge 10^5$ for $S \in \ms_S(s)$. Hence
	\[
		N_S(s)
	\lesssim
		m^s \, m^{-(1.1)10^5/3}
	\lesssim
		m^{-10^4}.
	\]
\end{itemize}
Combining these two then tells us that
\[ \textstyle
	\sum_{s \ge 1} N_S(s)
=
	\Ohb{ m^{-50} }.
\]

\medskip

\emph{Consider `large' sets.}
For $s \in \bn$, define the collection of sets
\[
	\ms_1(s) = \brb{ S \subseteq \ms \mid |S| = s, \, (\theta m)^{1/4} < d(S) \le \theta m/2 }.
\]
Let $N_1(s)$ be the expected number of sets $S \in \ms_1(s)$ with $I_G(S) \ge \beta d_G(S)$.
Following precisely the `large sets' part of the proof of \cite[Lemma 12]{ACF:kernel-expander}, it is shown that
\[ \textstyle
	\sum_{s \ge 1} N_1(s)
=
	\Ohb{ \rbr{9/10}^{m^{1/16}} }
=
	\Ohb{m^{-50}}.
\]

\medskip

\emph{Combing `small' and `large' sets.}
Recall that $\Phi_H(S) = 1 - I_H(S)/d_H(S)$, and so if $I_H(S) \le \beta d_H(S)$ then $\Phi_H(S) \ge 1/100$.
Hence, by Markov's inequality, we have
\[ \textstyle
	\prb{ \Phi_\mg' < 1/100 }
\le
	\sum_{s \ge 1}
	N(s)
=
	\Ohb{m^{-50}}.
\]
Since $\Phi_\mg \ge \min\{\Phi_\mg', 10^{-5}\}$, the lemma follows.
\end{Proof}

To prove \Cref{app:simplicity-ER}, we reference a result due to McKay and Wormald \cite{McKW:simple} on the probability of simplicity of a graph chosen uniformly at random given its degree sequences. The type of statement is well-known now, and similar statements can be found in \cite[Theorem 7.12]{H:random-graphs-book}, \cite{B:CM} and \cite[\S2.4]{B:random-graphs}.
The statement below is not as general as that given in \cite[Lemma 5.1]{McKW:simple}, but rather slightly weaker result that will be easier to apply in our situation.

\begin{thm}[{\cite[Lemma 5.1]{McKW:simple}}]
\label{app:simplicity-general}
	Let $\bm d = (d_1, ..., d_m)$ be a degree sequence, with $M = \sum_i d_i$ even and $\max_i d_i \le M^{1/4}$.
	Let $N = \sum_i d_i(d_i-1)$.
	Let $G \sim U(\bm d)$.
	Then
	\[
		\prb{ \simple{G} }
	\ge
		\expb{ - \tfrac12 N/M - \tfrac14 N^2/M^2 - \tfrac12 N^2 / M^3 + \Ohb{ 1/M^{1/4} } }.
	\]
\end{thm}

We show that the degree sequence of the core, with sufficient probability, satisfies $N \asymp M \asymp n$.

\begin{Proof}[Proof of \Cref{app:simplicity-ER}]
Write $\mc = \mc_F$, and let $\bm d = (d_1, ..., d_m)$ be as in the statement.
\Cref{app:kernel-size} says that the kernel at least $cn$ vertices with probability at least $1 - n^{-10}$, for a positive constant~$c$, and thus the same holds true for the core, and hence $\pr{M \le cn} \le n^{-10}$.
Observe that $d_\mc(x) \le d_F(x)$ for all $x \in \mc$. Hence, in the notation of \Cref{app:simplicity-general}, we have $N \le \sum_{i=1}^n d_F(i)^2$.
This formulation has the advantage that we know the distribution of $d_F(i)$, namely $\Bin(n-1,\lambda/n)$.

Moreover, we know the correlations between degrees:
for an edge $e$ of the complete graph~$K_n$, write $\eta_e = 1$ if the edge is open in $F$ and $\eta_e = 0$ if it is closed; then $d_F(i) = \sum_{y \neq x} \eta_{x,y}$, and the random variables $(\eta_e)_{e \in E(K_n)}$ are independent $\Bern(\lambda/n)$ random variables.
It is then straightforward to show that $\ex{N} \asymp n$. By Markov's inequality, we may then choose $C$ sufficiently large so that $\pr{N \ge Cn} \le \tfrac12$. The result now follows since
\[
	\prb{ \simple{G} } \ge \prb{ \simple{G} \mid N \le Cn, \, M \ge cn } \prb{ N \le Cn, \, M \ge cn }.
\qedhere
\]
\end{Proof}

We now give the proof of \Cref{app:largest-decoration}.

\begin{Proof}[Proof of \Cref{app:largest-decoration}]
This proof uses similar ideas to that of \cite[Lemma 4.2]{FR:giant-mixing}. 
Write $X_k$ for the number of sets $S$ of size $k$ such that $G[S]$ is a tree and and $e_G(S) \le 2$. Then  we have
\[ \textstyle
	\exb{ X_k }
\le
	\binom{n}{k} \rbb{ \tfrac{\lambda}{n} }^{k-1} k^{k-2} \rbb{1 - \tfrac\lambda n}^{ k(n-k) + \binom{k}2 - k - 1 }.
\label{eq:app:E(X)}
\nt
\]
For $k \le c n$, with positive constant $c$ sufficiently small, it is straightforward to see, using $\binom nk \le (ne/k)^k$, that there exists a positive constants $c'$ and $C'$ so that
\[
	\exb{ X_k } \le C' n e^{-c' k}.
\]

For $k \ge cn$, we are slightly more careful.
Note that (due to positive correlation) the number of isolated vertices is stochastically dominated from below by $\Bin(n,\tfrac12 e^{-\lambda})$, and hence, using concentration of the Binomial (Hoeffding) we have
\[
	\prb{ |\mg| \ge \rbr{1 - c} n } = \Ohb{ n^{-30} }
\Quad{whenever}
	c \le \tfrac14 e^{-\lambda}.
\]
Since a decoration is a subset of the giant, we need only consider $cn \le k \le (1-c)n$ with $c$ as small as we desire.
Write $k = \alpha n$ with $c \le \alpha \le 1-c$. Using Stirling's formula in \eqref{eq:app:E(X)}, we have
\begin{align*}
	\exb{ X_k }
&
\lesssim
	n^{-5/2} \rbb{ \alpha^\alpha(1-\alpha)^{1-\alpha} }^{-n} (\alpha n)^{\alpha n} \rbb{ \lambda/n }^{\alpha n} \expb{ - \lambda \alpha(1-\alpha) n - \tfrac12 \lambda \alpha^2 n + \tfrac32 \lambda \alpha }
\\&
\lesssim
	\expb{ - n h(\alpha,\lambda) }
\Quad{where}
	h(\alpha,\lambda) = (1-\alpha)\log(1-\alpha) - \alpha \log\lambda + \lambda \alpha(1-\alpha) + \tfrac12 \lambda \alpha^2.
\end{align*}
It is elementary to check that there exists a constant $\gamma$ so that $h(\alpha,\lambda) > \gamma$ for all $\alpha \in [c,1-c]$. Hence in this case also there exist positive constants $c'$ and $C'$ so that
\[
	\exb{ X_k } \le C' n e^{-c' k}.
\]
Taking $k = 101 \logn / c'$ we deduce the lemma by Markov's inequality.
\end{Proof}

\subsection{Typical Bounded Difference}

McDiarmid \cite{McD:bounded-diff} in 1989 used martingale methods to show concentration results for functions with bounded differences. Warnke \cite{W:typ-bounded-diff} in 2016 extended this to only consider differences on a high probability set (hence the name `typical' bounded differences).
We first give a general formulation, and then show how to apply it to Erd\H{o}s--R\'enyi graphs. We shall only formulate this for $\{0,1\}$-valued random variables, as this is all that we shall need; \cite{W:typ-bounded-diff} covers more general functions.

\medskip

Let $N \in \bn$, and let $Z = (Z_1, ..., Z_N)$ be a family of independent, $\{0,1\}$-valued random variables with $\pr{Z_k = 1} = p_k$. Let $\mz \subseteq \{0,1\}^N$, and suppose that the function $f : \{0,1\}^N \to \br$ satisfies the following \textit{typical Lipschitz condition}:
	there exist numbers $(\alpha_k)_{k=1}^N$ and $(\beta_k)_{k=1}^N$ such that, for each $k$, we have $\alpha_k \le \beta_k$ and, whenever $z,z' \in \{0,1\}^N$ with $z_j = z_j'$ for all $j \neq k$,
	\[
		\absb{ f(z) - f(z') }
	\le
		\begin{cases}
			\alpha_k	& \text{ if } z,z' \in \mz, \\
			\beta_k		& \text{ otherwise};
		\end{cases}
	\]
we say that $f : \{0,1\}^N \to \br$ satisfies $\TL(\mz,\alpha,\beta)$.

Warnke \cite{W:typ-bounded-diff} showed in essence that, providing we choose $\mz$ to be such that $\pr{Z \notin \mz}$ is sufficiently small, we can restrict our attention from worst-case differences to `typical' differences.

\begin{thm}[{\cite[Theorem 2 and Remark 10]{W:typ-bounded-diff}}]
\label{app:bdd-diff-general}
	Consider $Z$ and $\mz$ as above.
	Suppose $f : \{0,1\}^N \to \br$ satisfies $\TL(\mz,\alpha,\beta)$.
	For any numbers $(\gamma_k)_{k=1}^N \subseteq (0,1]^N$,
	we have
	\[
		\prb{ \absb{ f(Z) - \ex{ f(Z) } } \ge R }
	\le
		2 \expbb{ - \frac{R^2}{2 \sum_{k=1}^N p_k(1-p_k) (\alpha_k + e_k)^2 + 2 C R / 3} }
	+	\Gamma Q \, \prb{ Z \notin \mz },
	\]
	where $e_k = \gamma_k(\beta_k - \alpha_k)$, $C = \max_k (\alpha_k + e_k)$, $q_k \le \min_b \pr{Z_k = b}$, $\Gamma = \sum_k \gamma_k^{-1}$ and $Q = \sum_k q_k^{-1}$.
\end{thm}

We now formulate this in terms of an Erd\H{o}s--R\'enyi graph with $n$ vertices and edge probability $p = \lambda/n$ with $\lambda > 1$ fixed. We take $N = \binom{n}{2}$ to be the number of edges in the complete graph on $n$ vertices, and $Z_k$ to be the state of the $k$-th edge, so $q_k = \lambda/n$ for all $k$. Let us also restrict attention to bounded $f$, say $f : \{0,1\}^N \to [0,M]$; $M$ is allowed to depend on $n$. We then set $\beta_k = M$ and $\gamma_k = (Mn)^{-1}$ for all $k$. Since $\alpha_k \le \beta_k$, we then have $0 \le \alpha_k + e_k \le \alpha_k + 1/n$, and so $C \le \max_k \alpha_k + 1/n$. From this we obtain the following corollary.

\begin{cor}
\label{app:bdd-diff-ER}
	Consider an $\ER(n,\lambda/n)$ graph.
	Let $Z_k$ be the state of the $k$-th edge.
	Suppose $f : \{0,1\}^N \to [0,n]$ satisfies $\TL(\mz,\alpha,n)$.
	We have
	\[
		\prb{ \absb{ f(Z) - \ex{f(Z)} } \ge R }
	\le
		2 \expbb{ - \frac{R^2}{ 2 \lambda n^{-1} \sum_{k=1}^N (\alpha_k + 1/n)^2 + 2 C R / 3 } }
	+	n^6 \, \prb{ Z \notin \mz },
	\]
	where $C = \max_k \alpha_k + 1/n$.
\end{cor}

First we note some preliminary results which will go into the definition of $\mz$.
As proved at the start of this section (\S\ref{app}) there exist positive constants $c_0$ and $C_0$ so that the probability there exists a component with size in $[C_0 \logn, c_0 n]$ is $\Oh{n^{-100}}$.
A minor adaptation to the end of the proof of \cite[Theorem 2.14]{FK:intro-random-graphs} replacing $c_1 = c - \logn/n$ by $c_1 = c - \logn[2]/n$, shows that the probability there exists more than one component of size at least $c_0 n$ is $\Oh{n^{-100}}$.
Fix this $c_0$ and~$C_0$ for the rest of the paper.

\begin{lem}
\label{app:giant-size}
	Suppose $G \sim \ER(n,\lambda/n)$, with $\lambda > 1$ fixed, and write $\mg$ for its largest component (breaking ties arbitrarily).
	There exists a positive constant $c$ so that
	\[
		\prb{ |\mg| \le c n } = \Ohb{n^{-90}}.
	\]
\end{lem}

This result (and stronger versions) has been proved in the past; see, for example, \cite[Gap Theorem and Lemma 4]{P:tree-census-giant}.
We use the proof below as a warm-up for our method.

\begin{Proof}[Proof of \Cref{app:giant-size}]
By counting the expected number of subcritical components, we see that the expected size of the largest component is order $n$, say at least $2 c n$ for a positive constant $c$. For the details of this calculation, see the proof of \cite[Theorem 2.14]{FK:intro-random-graphs}.

Let $\mz$ be the set of all graphs with no components with size in $[C_0 \logn, c_0 n]$ and at most one component with size at least $c_0 n$.
Note also that $|\mg| \le n$.
Now consider a graph $G \in \mz$, and consider adding/removing an edge to obtain a new graph $G'$, and require that $G' \in \mz$.
	Adding an edge can only increase $|\mg|$, and by at most $C_0 \logn$, since no non-giant component has size larger than this.
	Removing an edge can only decrease $|\mg|$, and by at most $C_0 \logn$, since the new component must have size at most this.
Hence in the typical bounded difference formulation we may take $\alpha_k = C_0 \logn$ for all $k$; this gives $C \le 2 C_0 \logn$ and $\sum_{k=1}^N (\alpha_k + 1/n)^2 \le C_0^2 n^2 \logn[2]$. Take $R = c n$. Hence overall we have
\[
	\prb{ |\mg| \le cn }
\le
	2 \expbb{ - \frac{c^2 n^2}{2 \lambda C_0^2 n \logn[2] + 2 c C_0 n \logn} }
+	n^6 \, \Ohb{ n^{-100} }
=
	\Ohb{ n^{-90} }.
\qedhere
\]
\end{Proof}

\medskip

We now prove \Cref{app:kernel-size}. The proofs have two parts:
	first we find the expectation;
	then we show concentration.
To find the expectation, we use the contiguous model of Ding, Lubetzky and Peres, which was stated in \Cref{app:contiguous}; concentration will be proved using the typical bounded differences method above.

One part of the proof of \Cref{app:kernel-size} will rely on a short lemma, which we state and prove now.

\begin{lem}
\label{app:largest-dangling-component}
	Write $\md$ for the set of all edges of the giant whose removal causes the giant to split into two components (ie edges that are a \textit{bridge} for the giant).
	For $e \in \md$, write $R_e$ for the size of the smaller component resulting from the removal of $e$ (breaking ties arbitrarily). Write $R = \max_e R_e$.
	There exists a constant $C$ so that
	\[
		\prb{ R > C \logn } = \Ohb{ n^{-100} }.
	\]
\end{lem}

\begin{Proof}
	Let $F$ be a graph, and let $F \setminus e$ denote the graph $F$ with the edge $e$ closed. Let $G \sim \ER(n,\lambda/n)$ with $\lambda > 1$ fixed. Assuming $e$ is open in $F$, using independence of the edges we have
	\[
		\prb{ G = F } = \prb{ G = F \setminus e } \cdot p/(1-p).
	\]
	Now consider removing an edge $e \in \md$.
	Removing $e$ creates a new component. Since the probability that more than one component of size greater than $C_0 \logn$ exists is $\Oh{n^{-100}}$, using the above inequality we deduce the lemma in the following way.
	
	Let $\ma$ be the set of all graphs with the property that its largest component (breaking ties arbitrarily) has an edge whose removal disconnects the component leaving two components each of size greater than $C_0 \logn$;
	associate to each graph $F \in \ma$ an edge with the above removal property, and call this edge $e_F$.
	Let $\mb$ be the set of all graphs with two components of size greater than $C_0 \logn$.
	Then we have
	\[
		\brb{ F \setminus e_F \mid F \in \ma } \subseteq \mb.
	\]
	Hence using the first inequality we obtain
	\begin{align*}
		\prb{ G \in \ma }
	= \textstyle
		\sum_{F \in \ma}
		\prb{ G = F }
	&
	= \textstyle
		\sum_{F \in \ma}
		\prb{ G = F \setminus e_F } \cdot p/(1-p)
	\\&
	\le \textstyle
		\sum_{H \in \mb}
		\prb{ G = H } \cdot p/(1-p)
	=
		\prb{ G \in \mb } \cdot p/(1-p).
	\end{align*}
	Finally, if a graph $F$ has $R > C_0 \logn$, then $F \in \ma$. We know that $\pr{G \in \mb} = \Oh{n^{-100}}$, and hence we have proved the lemma.
\end{Proof}

An immediate corollary of this is that any $x \in \mg$ has at most $C \logn$ removal edges.

\begin{cor}[Maximum Removal Edges, \ref{gg:removal-edges}]
	There exists a constant $C$ so that
	\[
		\prb{ \exists \, x \in \mg \ST R(x) > C \logn } = \Ohb{n^{-100}}.
	\]
\end{cor}

\begin{Proof}
	Recall that we write $\mr(x)$ for the set of removal edges for $x$. Then $\mr(x)$ defines a `dangling component' in the sense of \Cref{app:largest-dangling-component}; the result follows.
\end{Proof}

\begin{rmkt}
\label{app:dangling-rmk}
	Observe that the proof of \Cref{app:largest-dangling-component} immediately extends in the following way: instead of considering sets $S \subseteq \mg$ which can be disconnected from $\mg$ by the removal of one edge, we require two edges to be removed.
\end{rmkt}


To prove this lemma, instead of considering the kernel directly, we construct a `contiguous model' of the giant, as given in \cite[Theorem 1]{DLP:anatomy-super}; we state the precise formulation now.

\begin{thm}[{\cite[Theorem 1]{DLP:anatomy-super}}]
\label{app:contiguous}
	Write $\mg$ for the (random) largest component of an $\ER(n,\lambda/n)$ graph (breaking ties arbitrarily), with $\lambda > 1$ fixed.
	Let $\theta$ be the unique solution in $[0,1]$ of $\theta e^{-\theta} = \lambda e^{-\lambda}$; write $\chi = \lambda - \theta$.
	Generate $\tilde\mg$ in the following way.
	\begin{enumerate}
		\item \textit{Constructing the kernel.}
		Let $\Chi \sim N(\chi,1/n)$ and, independently of $\Chi$, let $D_i \sim \Po(\Chi)$ for $i = 1,...,n$, conditioned that $\sum_i D_i \one{D_i \ge 3}$ is even (where $\Po(\xi) = 0$ if $\xi \le 0$). Let
		\[ \textstyle
			N_k = \absb{ \bra{i \mid D_i = k} }
		\Quad{and}
			N = \sum_{k \ge 3} N_k.
		\]
		Choose $\tilde\mk$ uniformly at random among all (multi-)graphs with $N_k$ vertices of degree $k$ for $k \ge 3$ (and hence with a total of $N$ vertices).
		
		\item \textit{Expanding to the core.}
		Obtain $\tilde\mc$ by replacing the edges of $\tilde\mk$ by iid paths of geometric length with parameter $1 - \theta$ (and mean $1/\theta$).
		
		\item \textit{Hanging dangling trees.}
		Obtain $\tilde\mg$ by attaching an independent Galton-Watson tree with $\Poisson(\theta)$ offspring to each vertex of $\tilde\mc$.
	\end{enumerate}
	Then, for any set $\ma$ of graphs, $\pr{\tilde\mg \in \ma} \to 0$ implies $\pr{\mg \in \ma} \to 0$.
\end{thm}


\begin{Proof}[Proof of \Cref{app:kernel-size}]
We start by finding the expectation via the contiguous model of \Cref{app:contiguous}.
In that notation, the size of the kernel is exactly $N$.
Note that $\Chi \sim N(\chi,1/n)$ and, independently of $\Chi$, $D_i \sim \Po(\Chi)$ for $i = 1,...,n$ conditioned that $\sum_i D_i \one{D_i \ge 3}$ is even, and $N = \sum_i \one{D_i \ge 3}$. Hence there exists a positive constant $c$ so that $\pr{ \Po(\Chi) \ge 3 } \ge 2c$. Then
\begin{align*}
	\prb{ N \le cn \mid \textstyle \sum_i D_i \one{D_i \ge 3} \ \text{even} }
&
\le
	\prb{ \textstyle \sum_i \one{D_i \ge 3} \le cn }
\big/
	\prb{ \textstyle \sum_i D_i \one{D_i \ge 3} \ \text{even} }
\\&
=
	\prb{ \Bin\rbb{n, \pr{D_i \ge 3}} \le cn }
\big/
	\prb{ \textstyle \sum_i D_i \one{D_i \ge 3} \ \text{even} }
\\&
\le
	\prb{ \Bin(n,2c) \le cn }
\big/
	\prb{ \textstyle \sum_i D_i \one{D_i \ge 3} \ \text{even} }
\\&
\le
	\expb{ - 2 c^2 n }
\big/
	\prb{ \textstyle \sum_i D_i \one{D_i \ge 3} \ \text{even} },
\end{align*}
with the final inequality holding by Hoeffding. We now show that the denominator is not too small. By conditioning that $\sum_{i=1}^{n-1} D_i \one{D_i \ge 3}$ is either odd or even and taking a minimum, we see that
\[
	\prb{ \textstyle \sum_{i=1}^n D_i \one{D_i \ge 3} \ \text{even} }
\ge
	\min\brb{ \pr{D_1 = 3}, \, \pr{D_1 = 4} }.
\]
Hence we see that there exists a positive constant $c'$ so that $\pr{ \textstyle \sum_i D_i \one{D_i \ge 3} \ \text{even} } \ge c'$.
From this we deduce that
\[
	\prb{ | \mk_{\tilde\mg} | \le c n } = \oh1,
\Quad{and hence}
	\exb{ | \mk_{\mg} |} \ge \tfrac12 c n,
\]
where $\mg$ is the largest component of an $\ER(n,\lambda/n)$ graph, with $\lambda > 1$ fixed.

\medskip

We now show concentration about the mean. To do this, we let $f(Z) = |\mk(\mg(G))|$ where we identify the vector $Z$ in $\{0,1\}^N$ with a graph $G$. Note that $f$ is increasing in the sense that if an edge is added to the largest component then (assuming there was no tie for largest component) the size of the kernel can only increase.
Let $\mz$ be the set of all graphs
	with no components with size in $[C_0 \logn, c_0 n]$
and
	precisely one component with size at least $c_0 n$.

We now consider how the size of the kernel can change when the state of an edge is changed. Write $K$ for the original kernel and $K'$ for the one obtained by the removal.
We show that $| |K| - |K'| | \lesssim \logn$;
the rest of the proof follows exactly as in \Cref{app:giant-size}.

\medskip

\newpage
\emph{Consider removing the edge $(x,y)$.}
\begin{itemize}
	\item []
	Removing an edge cannot increase the size of the kernel (using the uniqueness of the giant).
	
	\item []
	Suppose first that the removal of the edge does not disconnect the kernel (ie the edge is not a \textit{bridge} for the kernel).
	If the removed edge lies in a dangling tree, then this has no effect on the size of the kernel. Suppose the edge lies in a core-path, and so removing it removes the corresponding kernel-edge; let $z$ be an endpoint of the \emph{kernel-edge} corresponding to $(x,y)$.
		If $d_K(z) \ge 4$, then $z$ remains in the kernel.
%
		If $d_K(z) = 3$, then the other two core-paths incident to $z$ are joined into a single core-path, and the vertex $x$ is removed from the kernel but all other vertices remain.
	Hence we see that removing such an edge can remove at most 2 vertices from the kernel.
	
	Now suppose that the edge is a bridge for the kernel. Since the graph obtained by removing $(x,y)$ must be in $\mz$, the disconnected component must have size at most $C_0 \logn$.
	In addition to these, as in the above case, one further vertex is removed from the giant (through the joining of two core-paths at this vertex).
\end{itemize}

\emph{Consider adding the edge $(x,y)$.}
\begin{itemize}
	\item []
	Adding an edge cannot decrease the size of the kernel (using the uniqueness of the giant).
	
	\item []
	If the added edge has both ends in the same dangling tree, then this can add at most two vertices to the kernel, namely the most recent common ancestor of $x$ and $y$ in the dangling tree and its root.
	
	\item []
	Similarly, if $x$ and $y$ are in different dangling trees (allowing $x$ and $y$ to be the roots), then two vertices are added, namely the roots of the dangling trees in which $x$ and $y$ lie; that is, writing $r_z$ for the root of the dangling tree of $z$, we have $K' = K \cup \{r_x,r_y\}$.
	
	\item []
	Connecting a subcritical component to the giant can add at most $C_0 \logn + 1$, in the same way as removing a bridge can remove this many.
\qedhere
\end{itemize}
\end{Proof}

Having now completed the proof of the expansion properties \ref{gg:giant-expansion}, we move onto calculating the number of degree 1 vertices in the giant. (This result is proved independently of the expansion properties, except for quoting a small calculation that we did in \Cref{app:kernel-size}.)

\begin{lem}[Degree~1 Vertices in the Giant, \ref{gg:num-deg1-giant}]
\label{app:deg1}
	Suppose $G \sim \ER(n,\lambda/n)$ with $\lambda > 1$ fixed, and write $N$ for the number of degree~1 vertices in the giant.
	There exists a positive constant $c$ so that
	\[
		\prb{ N \le cn } = \Ohb{n^{-100}}.
	\]
\end{lem}

\begin{Proof}
First we lower bound the expected number of degree~1 vertices in the giant in the contiguous model; write $\tilde N$ for $N$ in the contiguous model. As in the proof of \Cref{app:kernel-size}, there exists a constant $c_1$ so that $\pr{|\mk(\tilde\mg)| \le c_1 n} = \oh1$. Note that $\PGW(\theta)$ trees are hung from the vertices of the kernel (as well as of the other vertices of the core). Since $\theta$ is a constant, the probability that the PGW-tree is just the root and one offspring is a constant. Hence there exists a positive constant~$c_2$ so that $\pr{\tilde N \le c_2 n} = \oh1$. Hence, as in the proof of \Cref{app:kernel-size}, we have $\ex{N} \ge \tfrac12 c_2 n$.

We now consider the effect on $N$ of adding/removing an edge. By recalling that both the original and new graph must have subcritical components of size at most $C_0 \logn$, we see that the largest change is at most $C_0 \logn$. The rest of the proof follows exactly as in \Cref{app:giant-size}.
\end{Proof}

Finally we prove \ref{gg:vtx-far-from-core} on the number of vertices `far away from the core'. Recall the definitions of $R$ and $\mw$ from \Cref{notation}, as well as that of $\log_{(M)}$.

\begin{prop}[Vertices Far from the Core, \ref{gg:vtx-far-from-core}]
\label{app:vtx-far-from-core}
	Suppose $G \sim \ER(n,\lambda/n)$ with $\lambda > 1$ fixed.
	There exists a constant $C$ so that, for all~$M$, for all~$n$ sufficiently large, we have
	\[
		\prb{ | \mg \setminus \mw^M | / |\mg| \ge (\log_{(M-1)}\! n)^{-4} }
	\le
		C n^{-50}.
	\]
\end{prop}


\begin{rmkt}
	The proposition needs $n$ to be sufficiently large depending on $M$, say $n \ge n_0(M)$. If we choose~$\omega_*$ such that $n_0(\omega_*(n)) \le n$, then we can apply the above proposition simultaneously for all~$M \le \omega_*(n)$. Further we can require $\omega_*(n) \le n$, and then deduce \ref{gg:vtx-far-from-core} by the union bound.
\end{rmkt}

We prove this by a sequence of lemmas. In essence, we split our consideration into dangling trees and parts of the giant that can be removed by removing a kernel-edge. The latter parts we think of as `bad' parts of the giant; we consider these first.

\medskip

For a graph $F$, write $\mP_F$ for the set of all $x \in \mg(F)$ with the following property: there exists an edge` of the giant whose removal disconnects the giant, leaving $x$ in the smaller component (breaking ties arbitrarily), and this smaller component is not a tree; that is, the removed edge disconnects the kernel.
We call $\mP_F$ the set of \textit{peninsular vertices} of the giant of $F$.

\begin{lem}
\label{app:peninsular}
	Suppose $G \sim \ER(n,\lambda/n)$ with $\lambda > 1$ fixed.
	There exists a constant $C$ so that
	\[
		\prb{ |\mP_G| > C \logn } = \Ohb{ n^{-50} }.
	\]
\end{lem}

\begin{Proof}
We first bound the expectation; we then use the typical bounded differences methodology.

We consider the kernel, and consider the expected number of vertices of the kernel that can be removed from the kernel by removing one kernel-edge. We showed in the proof of \Cref{app:giant-expansion}, specifically \eqref{eq:kernel-expander}, that the kernel is an expander, ie has isoperimetric constant $\Phi = \Th1$, with sufficiently high probability, and so any components that are `broken off' will have size in the kernel at most $1/\Phi = \Th1$. Moreover, any component of the kernel that is broken off cannot be a tree (since the kernel has minimal degree 3).
%
Using similar arguments as in \Cref{app:largest-dangling-component}, we now show that $\ex{|\mP_G|} = \Oh{\logn}$.

For a graph $F$, let $\me_F$ be the set of all edges $e$ with the property that removing $e$ disconnects the kernel of the largest component (breaking ties arbitrarily). Write $F \setminus \me_F$ for the graph with all edge $e \in \me_F$ closed. Then, with $p = \lambda/n$, we have
\[
	\prb{G = F} = \prb{G = F \setminus \me_F} \cdot \rbb{p/(1-p)}^{|\me_F|} \le \prb{G = F \setminus \me_F}.
\]
Write $|\mt_F'|$ for the number of vertices in non-tree subcritical components.
Then observe that $|\mP_F| \le D_F \cdot |\mt_{F \setminus \me_F}'|$ where $D_F$ is the size of the largest decoration. Write $\ma$ for the set of all graphs $F$ with $D_F \le C \logn$, with $C$ coming from \Cref{app:largest-decoration}. We then have
\begin{align*}
	\exb{ |\mP_G| \oneb{G \in \ma} }
&
= \textstyle
	\sum_{F \in \ma} |\mP_F| \prb{G = F}
\\&
\le \textstyle
	\sum_{F \in \ma} |\mP_F| \prb{G = F \setminus \me_F}
\\&
\le \textstyle
	\sum_{F \in \ma} D_F \absb{ \mt_{F \setminus \me_F}' } \prb{G = F \setminus \me_F}
\\&
\le \textstyle
	C \logn \sum_H \abs{ \mt_H' } \prb{G = H}
\\&
=
	C \logn \, \exb{ |\mt_G'| }.
\end{align*}
It is a standard calculation to show that $\ex{|\mt_G'|} = \Th1$ (see the proof of \cite[Theorem 2.4]{FK:intro-random-graphs}, for example).
\Cref{app:largest-decoration} says that $\pr{G \notin \ma} = \Oh{n^{-100}}$ and $|\mP_G| \le n$, hence $\ex{|\mP_G|} = \Oh{\logn}$.

We now apply the typical bounded differences methodology. Recall \Cref{app:largest-dangling-component}, which says that the largest set which can be disconnected from the giant by removing one edge has size at most $C \logn$, with probability $1 - \Oh{n^{-100}}$.

Let $\mz$ the set of all graphs with this connectivity property, with no components with size in $[C_0 \logn, c_0 n]$ and with precisely one with size at least $c_0 n$, which we call the giant. We now consider changing the state of a single edge. We show that the change is at most $(2C+C_0) \logn$; the rest of the proof follows exactly as in \Cref{app:giant-size}.

\medskip

\emph{Consider removing an edge.}
	A set which was previously connected to the giant via two edges can be reduced to only one edge, and so at most $C \logn$ can be added.
%
	A set which was previously connected to the giant via one edge can be removed, and so at most $C \logn$ can be removed.

\emph{Consider adding an edge.}
	A subcritical component can be connected to a vertex already in $\mP$, adding at most $C_0 \logn$ vertices to $\mP$. Also, the edge could make a dangling tree into a dangling non-tree; this component could then be disconnected by removing one edge, and so must have size at most $C \logn$.
%
	If the edge connects such a bad component to another vertex in the giant, then this can only remove vertices from $\mP$. It can remove two such bad components, each of which must have size at most $C \logn$, and so remove at most $2 C \logn$ vertices from $\mP$.
	%
\end{Proof}

We consider now the vertices in dangling trees that are `far from the core'.

\begin{lem}
\label{app:dangling-trees-decay}
	Suppose $G \sim \ER(n,\lambda/n)$ with $\lambda > 1$ fixed.
	There exist positive constants $c$ and~$C$ so that, writing $T_k$ for the number of dangling trees of size $k$, for $K \le C_0 \logn$, we have
	\[
		\exb{ \textstyle \sum_{k \ge K} k T_k } \le C e^{-cK}.
	\]
\end{lem}

\begin{Proof}
	Recall that the largest dangling tree is size order $\logn$ with high probability. Choose $c$ sufficiently small so that $c \cdot C_0 \le 10$; then $\exp{-c K} \ge n^{-10}$ for all $K \le C_0 \logn$.
	
	We upper bound $T_k$ by the number of tree components of size $k$ with exactly 1 edge coming out. Using standard combinatorial counting, we then see that
	\[
		\exb{ T_k } \lesssim n \rbb{\lambda e^{1-\lambda}}^k.
	\]
	The lemma then follows.
\end{Proof}

We now have the ingredients required to prove \Cref{app:vtx-far-from-core}.

\begin{Proof}[Proof of \Cref{app:vtx-far-from-core}]
As is often the case, we first look at expectation and then use the typical bounded differences methodology.

Write $\mt_k$ for the number of vertices in dangling trees of size at least $k$. (So $\mt_k = \sum_{k \ge K} k T_k$ in the notation of \Cref{app:dangling-trees-decay}.) Write
\[
	\mv^K = \brb{ x \in \mg \mid R(x) \ge K };
\Quad{so}
	\mv^{\log_{(K)}\!n} = \mw^K.
\]
For $K \le C_0 \logn$, we upper bound
\[
	\absb{ \mv^K }
\le \textstyle
	\mt_K + |\mP|.
\]
In \Cref{app:peninsular}, we showed that $|\mP| \le C \logn$ with probability $1 - \Oh{n^{-50}}$. Hence this combined with \Cref{app:dangling-trees-decay} says that
\[
	\exb{ \absb{ \mv^K } } \le C n e^{-c K} + C' \logn.
\]
In particular, if we take $K = 5 c^{-1} \log_{(M)}\!n$, then we see that
\[
	\exb{ \absb{ \mw^M } } \le \tfrac12 n \rbb{\log_{(M-1)}\!n}^{-4};
\]
moreover, we may allow $M$ to depend on $n$ if we require $M$ so that $\log_{(M)}\! n \ge 1$ (ie $M \le \log^*\!n - 1$).

We now use the typical bounded differences methodology, applied to $\mt_k$.
Let $\mz$ the set of all graphs with no components with size in $[C_0 \logn, c_0 n]$, with precisely one with size at least $c_0 n$, which we call the giant, and largest decoration to the giant of size at most $C_0 \logn$. We now consider changing the state of a single edge.

\medskip

\emph{Consider adding an edge.}
If both endpoints are in the giant, then this can only decrease the number of vertices in dangling trees (or leave unchanged); two full dangling trees can be removed, removing at most $2 C_0 \logn$ vertices.
If one endpoint is outside the giant, then a subcritical component can be added, adding at most $C_0 \logn$ vertices.

\emph{Consider removing an edge.}
If the core is unchanged, then a dangling tree has been `pruned', and hence at most $C_0 \logn$ vertices can be removed.
If the kernel is not disconnected, but a kernel-edge is removed, then all the vertices that were in this core-path then belong to dangling trees in the new graph, and hence at most $C_0 \logn$ can be added.
If the kernel is disconnected, then the vertices on the core-edge, of which there are at most $C_0 \logn$, can be added; all the vertices in the new subcritical component are removed, and so at most $C_0 \logn$ are removed.

\medskip

Hence $\mt_k$ can change by at most $3 C_0 \logn$. Taking $R = n / \logn[-6]$ in \Cref{app:bdd-diff-ER} we obtain
\[
	\prb{ \mt_k \ge n(C e^{-ck} + \logn[-6]) } = \Ohb{n^{-90}}.
\]
Note that $R \gg \log n$. This proves the proposition.
\end{Proof}

\end{document}

%% file: ER_preamble_biblio.tex
\usepackage[doi=false,isbn=false,url=false,
	hyperref=auto,
	sorting=nyt,
	maxnames=10,
	maxcitenames=3,
	backend=biber,
	texencoding=auto,
	giveninits=true,
	block=space,
	style=numeric,
]{biblatex}

\DefineBibliographyStrings{english}{%
	andothers = {et al}
} 
 
\DeclareFieldFormat{eprint:arxiv}{%
	\href{https://arxiv.org/abs/#1}{\emph{arXiv:\allowbreak #1}}}

\DeclareFieldFormat{eprint:mrnumber}{%
	\href{http://www.ams.org/mathscinet-getitem?mr=MR#1}{MR\allowbreak #1}}

\DeclareFieldFormat{eprint:jstor}{%
	\href{http://www.jstor.org/stable/#1}{JSTOR\allowbreak #1}}

\DeclareFieldFormat{eprint:customeprint}{%
	\em Available at \upshape \ttfamily \href{https://#1}{#1}}

\DeclareFieldFormat{eprint:online}{%
	Available \href{https://#1}{online}}

\DeclareFieldFormat{eprint:inprep}{%
	\emph{In preparation}}

\DeclareFieldFormat{eprint:manual}{%
	\emph{#1}}

\DeclareFieldFormat{eprint:onarxiv}{%
	\emph{#1}}

\DeclareFieldFormat{eprint:toappear}{%
	\mbox{\emph{To appear}}}

\DeclareFieldFormat{eprint:accepted}{%
	\emph{Accepted at {#1}; to appear}}

\DeclareSourcemap{
	\maps[datatype=bibtex]{
		\map{
			\step[fieldsource=mrnumber,		fieldtarget=eprint, final]
			\step[fieldset=eprinttype,		fieldvalue=mrnumber]
		}
		\map{
			\step[fieldsource=arxiv,		fieldtarget=eprint, final]
			\step[fieldset=eprinttype,		fieldvalue=arxiv]
		}
		\map{
			\step[fieldsource=jstor,		fieldtarget=eprint, final]
			\step[fieldset=eprinttype,		fieldvalue=jstor]
		}
		\map{
			\step[fieldsource=customeprint,	fieldtarget=eprint, final]
			\step[fieldset=eprinttype,		fieldvalue=customeprint]
		}
		\map{
			\step[fieldsource=online,		fieldtarget=eprint, final]
			\step[fieldset=eprinttype,		fieldvalue=online]
		}
		\map{
			\step[fieldsource=inprep,		fieldtarget=eprint, final]
			\step[fieldset=eprinttype,		fieldvalue=inprep]
		}
		\map{
			\step[fieldsource=manual,		fieldtarget=eprint, final]
			\step[fieldset=eprinttype,		fieldvalue=manual]
		}
		\map{
			\step[fieldsource=onarxiv,		fieldtarget=eprint, final]
			\step[fieldset=eprinttype,		fieldvalue=onarxiv]
		}
		\map{
			\step[fieldsource=toappear,		fieldtarget=eprint, final]
			\step[fieldset=eprinttype,		fieldvalue=toappear]
		}
		\map{
			\step[fieldsource=accepted,		fieldtarget=eprint, final]
			\step[fieldset=eprinttype,		fieldvalue=accepted]
		}
	}
}

\DeclareFieldFormat{doi}{%
	\href{https://doi.org/#1}{DOI}%
}

\DeclareFieldFormat{url}{%
	\href{https://#1}{\texttt{#1}}%
}

\DeclareFieldFormat{comments}{%
	\emph{#1}%
} 
 
\DeclareFieldFormat
	[article,inbook,incollection,inproceedings,patent,thesis,unpublished]
	{title}{{#1\isdot}}

\DeclareFieldFormat
	[article,book,inbook,incollection,inproceedings,patent,thesis,unpublished, online]
	{date}{\mkbibparens{#1}}
\DeclareFieldFormat
	[article,book,inbook,incollection,inproceedings,patent,thesis,unpublished]
	{volume}{\mkbibbold{#1}}
\DeclareFieldFormat
	{pages}{\mkbibparens{#1}}

\DeclareFieldFormat{edition}{%
	\ifinteger{#1}
	{\mkbibordedition{#1}~\bibstring{edition}\addcomma}
	{#1~\bibstring{edition}\addcomma}}

\renewbibmacro*{volume+number+eid}{%
	\printfield{volume}%
\setunit*{\adddot}%
	\printfield{number}%
\setunit{\space}%
	\printfield{eid}%
}

\DeclareBibliographyDriver{article}{%
	\usebibmacro{bibindex}%
	\usebibmacro{begentry}%
	\usebibmacro{author}%
\setunit{\addspace}%
	\usebibmacro{date}%
\newunit\newblock
	\usebibmacro{title}%
\newunit\newblock
	\printfield{journaltitle}\addperiod%
\setunit*{\addspace}%
	\usebibmacro{volume+number+eid}%
\setunit*{\addspace}\newblock
	\printfield{pages}
\setunit*{\addspace}
	\printfield{comments}%
	\usebibmacro{eprint}%
\setunit*{\addspace}%
	\printfield{doi}%
	\usebibmacro{finentry}%
}

\DeclareBibliographyDriver{online}{%
	\usebibmacro{bibindex}%
	\usebibmacro{begentry}%
	\usebibmacro{author}%
\setunit{\addspace}%
	\usebibmacro{date}%
\newunit\newblock
	\usebibmacro{title}%
\newunit\newblock
	\usebibmacro{eprint}%
	\usebibmacro{finentry}%
}

\DeclareBibliographyDriver{book}{%
	\usebibmacro{bibindex}%
	\usebibmacro{begentry}%
	\usebibmacro{author}%
\setunit{\addspace}\newblock
	\usebibmacro{date}%
\newunit\newblock
	\usebibmacro{title}%
\setunit*{\addspace}\newblock
	\printfield{edition}%
\newunit\newblock
	\printlist{publisher}
\setunit*{\addspace}%
	\printfield{volume}%
\setunit*{\addspace}\newblock%
	\printfield{comments}%
	\usebibmacro{eprint}%
\setunit*{\addspace}
	\printfield{doi}
	\usebibmacro{finentry}%
}

\DeclareBibliographyDriver{thesis}{%
	\usebibmacro{bibindex}%
	\usebibmacro{begentry}%
	\usebibmacro{author}%
\setunit{\addspace}\newblock
	\usebibmacro{date}%
\newunit\newblock
	\usebibmacro{title}.%
\setunit*{\addspace}\newblock
	\space Thesis, \printlist{institution}
	\usebibmacro{eprint}%
	\printfield{comments}%
\setunit*{\addspace}
	\printfield{doi}%
	\usebibmacro{finentry}%
}

\DeclareBibliographyDriver{inproceedings}{%
	\usebibmacro{bibindex}%
	\usebibmacro{begentry}%
	\usebibmacro{author}%
\setunit{\addspace}%
	\usebibmacro{date}%
\newunit\newblock
	\usebibmacro{title}%
\newunit\newblock
	\printfield{booktitle}\addcomma%
\setunit*{\addspace}%
	\printfield{series}\addcomma%
\setunit*{\addspace}\newblock
	\usebibmacro{editor}\addcomma
\setunit*{\addspace}\newblock
	\printlist{publisher}
\setunit*{\addspace}%
	\printfield{volume}%
\setunit*{\addspace}%
	\printfield{pages}
\setunit*{\addspace}
	\usebibmacro{eprint}%
	\printfield{comments}%
\setunit*{\addnnspace}
	\printfield{doi}
	\usebibmacro{finentry}%
}

\DeclareBibliographyDriver{inbook}{%
	\usebibmacro{bibindex}%
	\usebibmacro{begentry}%
	\usebibmacro{author}%
\setunit{\addspace}%
	\usebibmacro{date}%
\newunit\newblock
	\usebibmacro{title}%
\newunit\newblock
	\printfield{booktitle}\addcomma%
\setunit*{\addspace}%
	\printfield{series}\addcomma%
\setunit*{\addspace}\newblock
	\usebibmacro{editor}\addcomma
\setunit*{\addspace}\newblock
	\printlist{publisher}
\setunit*{\addspace}%
	\printfield{volume}%
\setunit*{\addspace}%
	\printfield{pages}
\setunit*{\addspace}
	\usebibmacro{eprint}%
	\printfield{comments}%
\setunit*{\addspace}
	\printfield{doi}
	\usebibmacro{finentry}%
}

\DeclareBibliographyDriver{incollection}{%
	\usebibmacro{bibindex}%
	\usebibmacro{begentry}%
	\usebibmacro{author}%
\setunit{\addspace}%
	\usebibmacro{date}%
\newunit\newblock
	\usebibmacro{title}%
\newunit\newblock
	\printfield{booktitle}\addcomma%
\setunit*{\addspace}%
	\printfield{series}\addcomma%
\setunit*{\addspace}\newblock
	\printlist{publisher}
\setunit*{\addspace}%
	\printfield{volume}%
\setunit*{\addspace}%
	\printfield{pages}
\setunit*{\addspace}
	\usebibmacro{eprint}%
	\printfield{comments}%
\setunit*{\addspace}
	\printfield{doi}
	\usebibmacro{finentry}%
}

\AtEveryBibitem{%
	\clearfield{day}%
	\clearfield{month}%
} 
 
\DeclareNameAlias{sortname}{given-family}
\DeclareNameAlias{default}{given-family}

%% file: ER_preamble_main.tex
\usepackage{indentfirst, environ}
\usepackage[utf8]{inputenc}

\usepackage{titlesec}
\usepackage{titletoc}

\titleformat{\title}
{\sffamily \Huge \bfseries \scshape \boldmath}{}{}{}
\titleformat{\section}
{\sffamily \Large \bfseries \scshape \boldmath}{\thesection}{1em}{}
\titleformat{\subsection}
{\sffamily \large \bfseries \scshape \boldmath}{\thesubsection}{1em}{}
\titleformat{\subsubsection}
{\sffamily \normalsize \bfseries \scshape \boldmath}{\thesubsubsection}{1em}{}
\titleformat{\paragraph}[runin]
{\sffamily \normalsize \bfseries \scshape \boldmath}{\theparagraph}{1em}{}
\titleformat{\subparagraph}[runin]
{\sffamily \normalsize \bfseries \scshape \boldmath}{\thesubparagraph}{1em}{}

\usepackage{authblk}

\usepackage{graphicx}
\usepackage[labelsep=period]{caption}

\usepackage{chngcntr}
\counterwithin{figure}{section}

\usepackage[UKenglish]{babel}
\usepackage{amsmath, amsthm, amssymb, mathrsfs, bm}
\usepackage{wasysym}

\usepackage{titlesec}
\usepackage{xifthen}

\usepackage{enumitem}
{\begin{itemize}[noitemsep,topsep=\smallskipamount,label={\bm\cdot}]}
	{\end{itemize}}

\newcommand{\romannumbering}{%
	\renewcommand{\labelenumi}{(\roman{enumi})}%
	\renewcommand{\theenumi}{(\roman{enumi})}}

\usepackage[dvipsnames,hyperref]{xcolor}

\usepackage{manyfoot}
\SetFootnoteHook{\hspace*{-1.8em}}
\DeclareNewFootnote{bl}[gobble]
\newcommand{\blfootnote}[1]{\footnotebl{\sffamily#1}}

\newcommand{\nt}{\addtocounter{equation}{1}\tag{\theequation}}

\newenvironment{Proof}[1][\proofname]{%
	\proof[\upshape\bfseries\sffamily\boldmath#1]
}{\endproof}

\newcommand{\Quad}[1]{\quad \text{#1} \quad}

\newcommand{\simple}[1]{\textnormal{\(#1\) simple}}

\newcommand{\tauc}{\tau_\textnormal{c}}

\newcommand{\XinG}[1]{
	X_{#1} \in \mathcal{G}_{#1}}

\newcommand{\existsXnotinG}[1]{
	\exists \, u \le #1 \textnormal{ s.t.\! } X_u \notin \mathcal{G}_u}

\renewcommand{\epsilon}{\varepsilon}
\newcommand{\logeps}{\log(1/\epsilon)}
\newcommand{\eps}{\epsilon}

\usepackage{xspace}

\let\originalexp\exp
\let\exp\relax

\DeclareRobustCommand{\exp} [1]{\originalexp(#1)}
\newcommand{\expb} [1]{\originalexp\bigl( #1 \bigr)}

\newcommand{\expbb}[1]{\originalexp\biggl( #1 \biggr)}

\newcommand{\gap}{}
\newcommand{\GAP}[1]{%
	\renewcommand{\gap}{\hspace{#1em}}}

\newcommand{\abs}  [1]{\lvert #1 \rvert}
\newcommand{\absb} [1]{\big\lvert #1 \bigr\rvert}
\newcommand{\absB} [1]{\Bigl\lvert #1 \Bigr\rvert}

\newcommand{\normb} [1]{\bigl\| #1 \bigr\|}

\newcommand{\rbr} [1]{ ( #1 ) }
\newcommand{\rbb} [1]{\bigl( #1 \bigr)}
\newcommand{\rbB} [1]{\Bigl( #1 \Bigr)}
\newcommand{\rbbb}[1]{\biggl( #1 \biggr)}

\newcommand{\sbb} [1]{\bigl[ #1 \bigr]}

\newcommand{\bra} [1]{ \{ #1 \} }
\newcommand{\brb} [1]{\bigl\{ #1 \bigr\}}
\newcommand{\brB} [1]{\Bigl\{ #1 \Bigr\}}
\newcommand{\brbb}[1]{\biggl\{ #1 \biggr\}}

\newcommand{\TV}{\textnormal{TV}}

\newcommand{\mix}{\textnormal{mix}}

\newcommand{\RW}{\textnormal{RW}}
\newcommand{\ST}{{\textnormal{ s.t.\! }}}

\newcommand{\TL}{\textnormal{TL}}

\newcommand{\Unif}{\textnormal{Unif}}
\newcommand{\iid}{\textnormal{iid}}
\newcommand{\ER}{\textnormal{ER}}

\newcommand{\Geo}{\textnormal{Geo}}
\newcommand{\Bin}{\textnormal{Bin}}
\newcommand{\Po}{\textnormal{Po}}
\newcommand{\Poisson}{\textnormal{Poisson}}

\newcommand{\Bern}{\textnormal{Bern}}

\newcommand{\tmix}{t_{\textnormal{mix}}}

\newcommand{\tunif}{t_{\textnormal{unif}}}

\newcommand{\pimin}{\pi_{\min}}
\newcommand{\isol}{\textnormal{isol}}
\newcommand{\tauisol}{\tau_{\textnormal{isol}}}

\newcommand{\PGW}{\textnormal{PGW}}

\newcommand{\Ggood}  [1]{
	{ \sg[#1] }				}

\newcommand{\Ggoodbb}[1]{
	{ \sg\biggl[#1\biggr] } }

\newcommand{\Gbad}  [1]{
	{ \sg[#1]^c }				}

\newcommand{\Hgood}[1]{
	{ \sh[#1] }				}
\newcommand{\Hgoodb}[1]{
	{ \sh\bigl[#1\bigr]	}	}

\newcommand{\Hbad}[1]{
	{ \sh[#1]^c }				}

\newcommand{\Xt}{(X_t)_{t\ge0}}

\newcommand{\Chi}{X}

\newcommand{\caps}{\cap \cdots \cap}

\newcommand{\floor}[1]{\lfloor #1 \rfloor}
\newcommand{\floorb}[1]{\bigl\lfloor #1 \bigr\rfloor}
\newcommand{\ceil}[1]{\lceil #1 \rceil}

\newcommand{\midb}{ \, \big| \, }

\newcommand{\midbb}{ \, \bigg| \, }

\newcommand{\ONE}  {\bm1}
\newcommand{\one}  [1]{\bm1( #1 )}
\newcommand{\oneb} [1]{\bm1\bigl( #1 \bigr)}

\newcommand{\tv} [1]{\| #1 \|_{\textnormal{TV}}}
\newcommand{\tvb}[1]{\bigl\| #1 \bigr\|_{\textnormal{TV}}}

\newcommand{\logn}[1][]{
	\ifthenelse{\equal{}{#1}}
	{\log n}
	{(\log n)^{#1}}
}

\newcommand{\loglogn}[1][]{
	\ifthenelse{\equal{}{#1}}
	{\log\log n}
	{(\log\log n)^{#1}}
}

\newcommand{\pr}[2][]{
	\ifthenelse{\equal{}{#1}}
	{\mathbb{P}(#2)}
	{\mathbb{P}_{#1}(#2)}
}
\newcommand{\prb}[2][]{
	\ifthenelse{\equal{}{#1}}
	{\mathbb{P}\bigl( #2 \bigr)}
	{\mathbb{P}_{#1}\bigl( #2 \bigr)}
}
\newcommand{\prB}[2][]{
	\ifthenelse{\equal{}{#1}}
	{\mathbb{P}\Bigl( #2 \Bigr)}
	{\mathbb{P}_{#1}\Bigl( #2 \Bigr)}
}
\newcommand{\prbb}[2][]{
	\ifthenelse{\equal{}{#1}}
	{\mathbb{P}\biggl( #2 \biggr)}
	{\mathbb{P}_{#1}\biggl( #2 \biggr)}
}
\newcommand{\prBB}[2][]{
	\ifthenelse{\equal{}{#1}}
	{\mathbb{P}\Biggl( #2 \Biggr)}
	{\mathbb{P}_{#1}\Biggl( #2 \Biggr)}
}
\newcommand{\prs}[2][]{
	\ifthenelse{\equal{}{#1}}
	{\mathbb{P}\left( #2 \right)}
	{\mathbb{P}_{#1}\left( #2 \right)}
}

\newcommand{\prg}[2][]{
	\ifthenelse{\equal{}{#1}}
	{\mathbb{P}^\mg(#2)}
	{\mathbb{P}^\mg_{#1}(#2)}
}
\newcommand{\prgb}[2][]{
	\ifthenelse{\equal{}{#1}}
	{\mathbb{P}^\mg\bigl( #2 \bigr)}
	{\mathbb{P}^\mg_{#1}\bigl( #2 \bigr)}
}
\newcommand{\prgB}[2][]{
	\ifthenelse{\equal{}{#1}}
	{\mathbb{P}^\mg\Bigl( #2 \Bigr)}
	{\mathbb{P}^\mg_{#1}\Bigl( #2 \Bigr)}
}
\newcommand{\prgbb}[2][]{
	\ifthenelse{\equal{}{#1}}
	{\mathbb{P}^\mg\biggl( #2 \biggr)}
	{\mathbb{P}^\mg_{#1}\biggl( #2 \biggr)}
}
\newcommand{\prgBB}[2][]{
	\ifthenelse{\equal{}{#1}}
	{\mathbb{P}^\mg\Biggl( #2 \Biggt)}
	{\mathbb{P}^\mg_{#1}\Biggl( #2 \Biggr)}
}
\newcommand{\prgs}[2][]{
	\ifthenelse{\equal{}{#1}}
	{\mathbb{P}^\mg\left( #2 \right)}
	{\mathbb{P}^\mg_{#1}\left( #2 \right)}
}

\newcommand{\ex}[2][]{
	\ifthenelse{\equal{}{#1}}
	{\mathbb{E}(#2)}
	{\mathbb{E}_{#1}(#2)}
}
\newcommand{\exb}[2][]{
	\ifthenelse{\equal{}{#1}}
	{\mathbb{E}\bigl( #2 \bigr)}
	{\mathbb{E}_{#1}\bigr( #2 \bigr)}
}
\newcommand{\exB}[2][]{
	\ifthenelse{\equal{}{#1}}
	{\mathbb{E}\Bigl( #2 \Bigr)}
	{\mathbb{E}_{#1}\Bigl( #2 \Bigr)}
}
\newcommand{\exbb}[2][]{
	\ifthenelse{\equal{}{#1}}
	{\mathbb{E}\biggl( #2 \biggr)}
	{\mathbb{E}_{#1}\biggl( #2 \biggr)}
}
\newcommand{\exBB}[2][]{
	\ifthenelse{\equal{}{#1}}
	{\mathbb{E}\Biggl( #2 \Biggr)}
	{\mathbb{E}_{#1}\Biggl( #2 \Biggr)}
}

\newcommand{\Var}[2][]{
	\ifthenelse{\equal{}{#1}}
	{\mathbb{V}\textnormal{ar} \bigl(#2\bigr)}
	{\mathbb{V}\textnormal{ar}_{#1} \bigl(#2\bigr)}
}
\newcommand{\VAR}[2][]{
	\ifthenelse{\equal{}{#1}}
	{\textnormal{Var}(#2)}
	{\textnormal{Var}_{#1}(#2)}
}

\newcommand{\Om}  [1]{\Omega( #1 )}

\newcommand{\om}  [1]{\omega( #1 )}

\newcommand{\Oh}  [1]{\mathcal{O}( #1 )}
\newcommand{\Ohb} [1]{\mathcal{O}\bigl( #1 \bigr)}

\newcommand{\oh}  [1]{o( #1 )}

\newcommand{\Th}  [1]{\Theta( #1 )}

\newcommand{\bn}{\mathbb{N}}

\newcommand{\bp}{\mathbb{P}}
\newcommand{\bq}{\mathbb{Q}}
\newcommand{\br}{\mathbb{R}}

\newcommand{\bz}{\mathbb{Z}}

\newcommand{\ma}{\mathcal{A}}
\newcommand{\mb}{\mathcal{B}}
\newcommand{\mc}{\mathcal{C}}
\newcommand{\md}{\mathcal{D}}
\newcommand{\me}{\mathcal{E}}
\newcommand{\mf}{\mathcal{F}}
\newcommand{\mg}{\mathcal{G}}

\newcommand{\mi}{\mathcal{I}}
\newcommand{\mj}{\mathcal{J}}
\newcommand{\mk}{\mathcal{K}}

\newcommand{\mn}{\mathcal{N}}

\newcommand{\mP}{\mathcal{P}}

\newcommand{\mr}{\mathcal{R}}
\newcommand{\ms}{\mathcal{S}}
\newcommand{\mt}{\mathcal{T}}
\newcommand{\mU}{\mathcal{U}}
\newcommand{\mv}{\mathcal{V}}
\newcommand{\mw}{\mathcal{W}}

\newcommand{\mz}{\mathcal{Z}}

\renewcommand{\sc}{\mathscr{C}}
\newcommand{\se}{\mathscr{E}}

\newcommand{\sg}{\mathscr{G}}
\newcommand{\sh}{\mathscr{H}}
\newcommand{\so}{\mathscr{O}}

\newcommand{\ninf}{{n\to\infty}}

\usepackage[
	colorlinks,
	pdfauthor = {Sam\ Olesker-Taylor\ and\ Perla\ Sousi},
	pdftitle = {Cutoff\ for\ Random\ Walk\ on\ Dynamical\ Erdos--Renyi\ Graph}
]{hyperref}
\hypersetup{
	colorlinks,
	linkcolor	= {blue},
	urlcolor	= {blue},
	citecolor	= {ForestGreen}
}

\usepackage{aliascnt}
\usepackage[capitalise]{cleveref}

\newtheoremstyle{sfsl}
{1\baselineskip}		
{1\baselineskip}		
{\slshape}				
{}						
{\bfseries\sffamily}	
{.}						
{0.5em}					
{\thmname{#1}\thmnumber{ #2}\thmnote{ \textnormal{\sffamily(#3)}}}

\newtheoremstyle{sfup}
{1\baselineskip}		
{1\baselineskip}		
{\upshape}				
{}						
{\bfseries\sffamily}	
{.}						
{0.5em}					
{\thmname{#1}\thmnumber{ #2}\thmnote{ \textnormal{\sffamily(#3)}}}

\theoremstyle{sfsl}

\newtheorem*{thm*}{Theorem}
\newtheorem{thm} {Theorem}[section]

\newtheorem*{lem*}    {Lemma}
\newtheorem{lem} [thm]{Lemma}

\newtheorem*{prop*}{Proposition}
\newtheorem{prop} [thm]{Proposition}

\newtheorem{cor} [thm]{Corollary}

\theoremstyle{sfup}

\newtheorem{defnT}[thm]{Definition}
	\newenvironment{defnt}
	{\pushQED{\qed}\defnT}
	{\popQED\enddefnT}

\newtheorem*{defn*}{Definition}
\newtheorem*{DEFNT}{Definition}
\newenvironment{defnt*}
	{\pushQED{\qed}\DEFNT}
	{\popQED\endDEFNT}

\newtheorem{rmkT}[thm]{Remark}
	\newenvironment{rmkt}
	{\pushQED{\qed}\rmkT}
	{\popQED\endrmkT}

\newtheorem*{rmk*} {Remark}
\newtheorem*{rmkTT}{Remark}
\newenvironment{rmkt*}
	{\pushQED{\qed}\rmkTT}
	{\popQED\endrmkTT}

\newtheorem*{rmks*} {Remarks}
\newtheorem*{rmksTT}{Remarks}
\newenvironment{rmkst*}
	{\pushQED{\qed}\rmksTT}
	{\popQED\endrmksTT}

\newtheorem*{clm*} {Claim}
\newtheorem*{clmTT}{Claim}
	\newenvironment{clmt*}
	{\pushQED{\qed}\clmTT}
	{\popQED\endclmTT}

\newtheorem*{exm*} {Example}
\newtheorem*{exmTT}{Lemma}
	\newenvironment{exmt*}
	{\pushQED{\qed}\exmTT}
	{\popQED\endexmTT}

\newtheorem{notationT}[thm]{Notation}
	\newenvironment{notationt}
	{\pushQED{\qed}\notationT}
	{\popQED\endnotationT}

\newtheorem*{notation*} {Notation}
\newtheorem*{notationTT}{Notation}
\newenvironment{notationt*}
	{\pushQED{\qed}\notationTT}
	{\popQED\endnotationTT}

\newtheorem{innercustomgeneric}{\customgenericname}
\providecommand{\customgenericname}{}
\newcommand{\newcustomtheorem}[2]{%
	\newenvironment{#1}[1]
	{%
		\renewcommand\customgenericname{#2}%
		\renewcommand\theinnercustomgeneric{##1}%
		\innercustomgeneric
	}
	{\endinnercustomgeneric}
}

\newcustomtheorem{customthm}{Theorem}

\Crefname{thm}{Theorem}{Theorems}
\Crefname{thmT}{Theorem}{Theorems}
\Crefname{lem}{Lemma}{Lemmas}
\Crefname{lemT}{Lemma}{Lemmas}
\Crefname{cor}{Corollary}{Corollaries}
\Crefname{corT}{Corollary}{Corollaries}
\Crefname{prop}{Proposition}{Propositions}
\Crefname{propT}{Proposition}{Propositions}

\newtheorem*{note*} {Note}
\newtheorem*{noteTT}{Note}
	\newenvironment{notet*}
	{\pushQED{\qed}\noteTT}
	{\popQED\endnoteTT}